\documentclass[upref]{amsart}
\usepackage{latexsym,amsmath,amstext,mathrsfs,amssymb,amsthm,amscd,amsopn,verbatim,graphics}

\usepackage{times}
\numberwithin{equation}{section}

\DeclareMathOperator{\laction}{L}
\DeclareMathOperator{\raction}{R}

\DeclareMathOperator{\id}{id}

\DeclareMathOperator{\partr}{H}

\DeclareMathOperator{\conj}{c}

\DeclareMathOperator{\tange}{T}

\DeclareMathOperator{\verti}{V}

\DeclareMathOperator{\GL}{GL}

\DeclareMathOperator{\coh}{H}
\newcommand{\tspace}[1]{\tange\! #1}

\newcommand{\tangent}[2]{\tspace_{#1}{#2}}

\DeclareMathOperator{\pr}{pr}

\DeclareMathOperator{\orient}{or}

\DeclareMathOperator{\ev}{ev}

\newcommand{\dd}{\mathrm{d}}

\DeclareMathOperator{\tr}{Tr}

\DeclareMathOperator{\ddtt}{\frac{d}{d\mathnormal{t}}}
\DeclareMathOperator{\ddss}{\frac{d}{d\mathnormal{s}}}
\newcommand{\holo}[3]{\partr(#1)\arrowvert_{#2}^{#3}}
\newcommand{\holon}[3]{\partr\!\left(#1;#2;#3\right)}
\newcommand{\holonom}[5]{\partr\!\left(#1;#2;#3;#4,#5\right)}
\newcommand{\holonomy}[6]{\partr\!\left(#1;#2;#3,#4;#5,#6\right)}

\DeclareMathOperator{\op}{\mathnormal{op}}

\DeclareMathOperator{\ad}{ad}

\DeclareMathOperator{\Ad}{Ad}

\newcommand{\lo}{\mathrm{L}}
\newcommand{\PM}{\mathrm{P}\!M}
\newcommand{\unint}{\mathrm{I}}

\newcommand{\Lg}{\mathfrak{g}}

\DeclareMathOperator{\End}{End}

\newtheorem{Thm}{Theorem}[section]
\newtheorem{Prop}[Thm]{Proposition}
\newtheorem{Lem}[Thm]{Lemma}
\newtheorem{Cor}[Thm]{Corollary}

\theoremstyle{remark}
\newtheorem{Rem}[Thm]{Remark}
\newtheorem*{Ack}{Acknowledgment}

\theoremstyle{definition}

\newtheorem{Def}[Thm]{Definition}

\newtheorem*{Reg*}{Regularization procedure}

\newcommand{\Lie}[2]{{\left[{\,{#1}\,,\,{#2}\,}\right]}}

\newcommand{\bbR}{{\mathbb{R}}}

\newcommand{\calA}{\mathcal{A}}

\newcommand{\calH}{\mathcal{H}}

\newcommand{\calG}{\mathcal{G}}

\newcommand{\calM}{\mathcal{M}}

\newcommand{\calW}{\mathcal{W}}

\newcommand{\tildep}{\widetilde{p}}

\DeclareMathOperator{\diff}{Diff}

\newcommand\qq{\rm}
\newcommand\cmp[1]{{\qq Commun.\ Math.\ Phys.\ \bf #1}}

\newcommand\anm[1]{{\qq Ann.\ Math.\ \bf #1}}

\begin{document} 

\title[Holonomy and parallel transport\ldots]{Holonomy and parallel transport in the differential geometry of the space of loops and the groupoid of generalized gauge transformations}

\author[C.~A.~Rossi]{Carlo~A.~Rossi}
\address{Dept.\ of Mathematics --- Technion ---  
32000 Haifa --- Israel}  
\email{crossi@techunix.technion.ac.il}
\thanks{C.~A.~Rossi acknowledges partial support of Aly Kaufman Fellowship}

\begin{abstract}
The motivation for this paper stems~\cite{CR} from the need to construct explicit isomorphisms of (possibly nontrivial) principal $G$-bundles on the space of loops or, more generally, of paths in some manifold $M$, over which I consider a fixed principal bundle $P$; the aforementioned bundles are then pull-backs of $P$ w.r.t.\ evaluation maps at different points.

The explicit construction of these isomorphisms between pulled-back bundles relies on the notion of {\em parallel transport}.
I introduce and discuss extensively at this point the notion of {\em generalized gauge transformation between (a priori) distinct principal $G$-bundles over the same base $M$}; one can see immediately that the parallel transport can be viewed as a generalized gauge transformation for two special kind of bundles on the space of loops or paths; at this point, it is possible to generalize the previous arguments for more general pulled-back bundles.
 
Finally, I discuss how flatness of the reference connection, w.r.t.\ which I consider holonomy and parallel transport, is related to horizontality of the associated generalized gauge transformation.
\end{abstract}

\maketitle

\tableofcontents

\section{Introduction}
In the last section of~\cite{CR}, where the Topological Quantum Field theoretical background behind the higher-order cohomology groups of the space of knots in $\bbR^m$, for $m\geq 3$, of Cattaneo--Cotta-Ramusino--Longoni~\cite{CCL} was explicitly constructed, some hints were mentioned towards possible generalizations of the computations in~\cite{CR}, dealing with iterated Chen-type integrals, to a {\em nontrivial principal $G$-bundle $P$}; in fact, the main object throughout the paper was the so-called {\em generalized holonomy}, viewed as an iterated integral, on some fixed principal $G$-bundle $P$, which we assumed to be trivial in order to simplify the computations.

In the last section, we addressed informally the problem of defining iterated integrals, as our expression for the parallel transport, for a nontrivial bundle $P$.
The same problem arises also in Section 2 of~\cite{CCR}: the authors discuss the notion of iterated Chen integrals in relationship with so-called ``special connections'' on the space of horizontal paths in a (not necessarily trivial) principal bundle over a $4$-manifold $M$, namely they need iterated Chen integrals of forms of the adjoint type on $M$ in order to compute the curvature of special connections.
An explicit formula is displayed there for the Chen bracket; still, as the authors themselves pointed out, it is without details, to which they planned to dedicate some further paper.
One of the main subtleties of the task in~\cite{CCR} and~\cite{CR} lies therein, that we need to identify pulled-back bundles of $P$ on the space of loops (or, more generally, of paths) in a general manifold $M$ w.r.t.\ evaluation maps at different points.
We sketched without details therein some arguments leading to the answer; in particular, we pointed out the importance in this task of the {\em holonomy} and, more generally, of the {\em parallel transport}.
In fact, the parallel transport, which depends explicitly by its very construction on a connection $A$ on $P$, defines an isomorphism between two particular type of pulled-back bundles on the space of loops or paths in $M$.

In the present paper, I explain all details of this construction, performing the complete computations we hinted at in the last section of~\cite{CR}.
The paper is organised as follows: in Section~\ref{sec-backgr}, I recall the main notions, namely connections on principal $G$-bundles over the manifold $M$, and I discuss the two equivalent ways of defining a connection, since I will make use of them both in the subsequent computations.   

In Section~\ref{sec-holpartransp}, I recall the notions of {\em holonomy w.r.t.\ a chosen connection $A$ of a loop $\gamma$ in $M$} and {\em parallel transport w.r.t.\ $A$ of a path $\gamma$}; I then state and prove two technical lemmata containing information about some sort of equivariance displayed by holonomy and parallel transport w.r.t.\ the actions of the structure group $G$ of $P$ and the gauge group $\calG_P$of $P$.
The contents of the first two sections are standard facts in gauge theory; I will nonetheless review them in some detail to fix conventions and notations.

In Section~\ref{sec-isombun}, I discuss the notion of isomorphisms of $G$-principal bundles (a priori distinct) over $M$.It is well known that fibre-preserving, $G$-equivariant automorphisms of a {\em fixed} $G$-principal bundles are in one-to-one correspondence to maps from $P$ to $G$, equivariant w.r.t.\ conjugation on $G$; similarly, fibre-preserving, $G$-equivariant (iso)morphisms between (a priori distinct) principal $G$-bundles over the same base space $M$ are in one-to-one correspondence to what I call {\em generalized gauge transformations}.
These are maps from the {\em fibred product} of the aforementioned bundles (a concept somehow mimicking the notion of Whitney sum of vector bundles) with values in $G$, equivariant w.r.t.\ an action of the product group $G\times G$ on $G$, which, when restricted to the diagonal subgroup $G$, restricts to the conjugation on $G$.
The key element of the above correspondence lies in the canonical map from the fibred product of a bundle $P$ with itself associated to the identity map of $P$; I also discuss its properties.
I then proceed with the discussion of the structure of the set of generalized gauge transformations; it turns out to be (obviously) a {\em groupoid} over the category of principal $G$-bundles over $M$; this notion generalizes that of gauge group, and I therefore speak of the {\em groupoid of generalized gauge transformations} of $M$ and $G$.

In Section~\ref{sec-groupoidconstr}, I partly review, partly explain some general constructions for groupoids; in particular, I recall the concept of {\em actions of groupoids on sets} and construct what I call the {\em generalized conjugation for groupoids}, which is an action of the product of a groupoid with itself on the groupoid itself, which mimics somehow the conjugation of a group on itself, a notion which no longer makes sense for a groupoid.
I also discuss the notion of {\em equivariant morphisms} from (left or right) groupoid actions to other (left or right) groupoid actions.

In Section~\ref{sec-applpullbun}, I interpret holonomy and parallel transport as generalized gauge transformations for a particular type of pulled-back bundles over the space of loops or paths in $M$ of a fixed principal bundle $P$ over $M$ w.r.t.\ evaluation maps.
Thus, there is a map associating to a connection on $P$ a generalized gauge transformation; interpreting the gauge group of $P$ as a groupoid, hence the space of connections on $P$ as a $\calG_P$-space, the technical lemmata of Section~\ref{sec-holpartransp} may be interpreted in the sense that there is an equivariant map of groupoid actions from the space of connections to the groupoid of generalized gauge transformations.

In Section~\ref{sec-flatness}, I discuss the important consequences of the restriction of the map discussed in Section~\ref{sec-applpullbun} between the space of connections and the groupoid of generalized gauge transformations to the space of {\em flat connections}.
Namely, the flatness (which has a non-abelian cohomological interpretation) implies in a highly nontrivial way the {\em horizontality} of a the corresponding generalized gauge transformation; the key step for this achieving this result lies in the well-known fact that flat connections induce representations of the {\em homotopy group of $M$}.

\begin{Ack}
I thank A.\ S.\ Cattaneo and G.\ Felder for many inspiring suggestions, mainly on the subject of groupoids and for constant support; I also acknowledge the pleasant atmosphere at the Department of Mathematics of the Technion, where this work was accomplished.
\end{Ack}

\section{Background definitions}\label{sec-backgr}
In this section, I introduce the main notions and notations I use throughout the paper; among them, I want to discuss in detail the notion of {\em connection on principal bundles}.
Let me notice that I work on a general principal bundle $P\overset{\pi}\mapsto M$ over a real manifold $M$ (if not otherwise stated, $M$ is assumed to be connected and paracompact); I do not assume any particular property on the structure group $G$.
By $\Lg$ I denote the corresponding Lie algebra. 

\begin{Def}\label{def-connform}
A {\em connection $1$-form} on the principal bundle $P$ is a $1$-form $A$ on $P$ with values in $\Lg$, satisfying the two following requirements:
\begin{itemize}
\item[i)] ({\em Equivariance}) 
\begin{equation}\label{eq-equiv}
\raction_g^*\!A=\Ad\!(g^{-1})A,
\end{equation} 
where by $\raction_g$ I have denoted the (free) right action of $G$ on $P$. 
\item [ii)] ({\em Verticality}) For any $\xi\in\Lg$,
\begin{equation}\label{eq-verticality}
A_p\!\left(\tangent{e}{\laction_p}\!(\xi)\right)=\xi,\quad \forall p\in P;
\end{equation} 
here, I have denoted by $\raction_p$ the fibre injection $G\to P$ given by $g\mapsto pg$.
\end{itemize}
\end{Def}

For our purposes, it is better to introduce a slightly different characterisation of connections; later, I will discuss the relationship between them.

First of all, a tangent vector $X$ to $P$ at the point $p$ is said to be {\em vertical}, if it satisfies the equation 
\[
\tangent{p}{\pi}(X)=0.
\]
The {\em vertical space $\verti_p\!P$}, consisting of all vertical tangent vectors at $p$, for any $p\in P$, is isomorphic to the Lie-algebra $\Lg$ via
\[
\xi\mapsto \tangent{e}{\laction_p}(\xi),\quad \xi\in \Lg.
\] 
It turns out that the vertical spaces $\verti_p\!P$ can be glued together to give a smooth vector bundle, the so-called {\em vertical bundle $\verti\!P$} (whose typical fibre is isomorphic to the Lie algebra $\Lg$); it is clearly a subbundle of $\tspace\!P$.

Alas, there is in principle no canonical way to define a complement of $\verti\! P$ w.r.t.\ Whitney sum, i.e.\ there is no {\em canonical} bundle $\coh\!P$ such that $\verti\!P\oplus\coh\!P=\tspace\!P$.
It turns out that the choice of such a complementary bundle relies on the choice of a connection $1$-form $A$ on $P$, as it is motivated by the following definition:
\begin{Def}\label{def-horiz}
Given a connection $1$-form $A$ on $P$, a tangent vector $X_p$ to $P$ at $p$ is said to be {\em $A$-horizontal} if the following equation holds:
\begin{equation*}
A_p\!\left(X_p\right)=0. 
\end{equation*}
(For the sake of brevity, when the connection $1$-form $A$ is clear from the context, I simply speak of horizontal vectors.)
\end{Def}

\subsection{The bridge between connection $1$-forms and horizontal bundle}
I discuss, for the sake of completeness, how a connection $1$-form $A$ on $P$ gives rise to a {\em smooth} assignment to any $p\in P$ of a subspace $\coh_p\!P\subset \tangent{p}{P}$, such that the following two requirements hold:
\begin{itemize}
\item[i)] $\tangent{p}{P}=\verti_p\!P\oplus\coh_p\!P$;
\item[ii)] $\coh_{pg}\!P=\tangent{p}{\raction_g}\left(\coh_p\!P\right)$, for any $g\in G$,
\end{itemize} 
and vice versa. 
In fact, the two requirements listed above are also an alternative definition of connection as a way of splitting the tangent bundle of $P$ into the Whitney sum of the vertical bundle, which is clearly $G$-invariant as a consequence of the identity
\[
\raction_g\circ\!\laction_p=\laction_{pg}\circ\!\conj\!\left(g^{-1}\right),\quad \forall p\in P, g\in G,
\]
and some horizontal bundle.

Given such an assignment, it is possible to define a corresponding connection $1$-form via
\begin{equation}\label{eq-distconn}
A_p\!\left(X_p\right)\colon=\xi_{X_p},
\end{equation}
where $\xi_{X_p}$ is the unique element of $\Lg$ which corresponds to the vertical part of $X_p$ w.r.t.\ the above splitting, i.e.\
\[
\tangent{e}{\laction_p}\!\left(\xi_{X_p}\right)=X_p^v,
\]
where $X_p^v$ is the vertical part of $X_p$ w.r.t.\ the above splitting.
It follows immediately that the horizontal space at $p\in P$ is exactly the kernel of $A_p$. 
The very definition of vertical space and Equation (\ref{eq-distconn}) imply together that $A$ satisfies (\ref{eq-verticality}). 
The invariance given by Condition ii) above ensures that both projections into vertical and horizontal subspace are $G$-invariant, which in turn leads to (\ref{eq-equiv}).

On the other hand, given a connection $1$-form $A$ on $P$, the corresponding splitting of $\tspace\!P$ at some point $p\in P$ is given by
\[
X_p=\left(X_p-X_{p}^A\right)+X_{p}^A,\quad X_p^A\colon=\tangent{e}{\laction_p}\!\left(A_p(X_p)\right),\quad p\in P.
\]
In fact, by (\ref{eq-verticality}), it is easy to prove that the linear operator $\tangent{e}{\laction_p}\!\circ\! A_p$ is the projection onto the vertical subspace; its kernel is therefore the horizontal subspace corresponding to the choice of $A$.
Condition (\ref{eq-equiv}), in turn, alongside the identity
\[
\laction_{pg}\circ \conj\!\left(g^{-1}\right)=\raction_g\circ \laction_p,\quad \forall p\in P, g\in G,
\]
ensures that the corresponding distribution is $G$-equivariant ($\conj$ denotes the conjugation on $G$).
Finally, let me spend some words on the concept of {\em gauge transformation of a principal $G$-bundle $P$ over $M$}.
A gauge transformation $\sigma$ of $P$ is a (smooth) map from $P$ to $P$, enjoying the following two properties:
\[
\pi\circ \sigma =\pi,\quad \sigma(pg)=\sigma(p)g,\quad \forall p\in P,g\in G.
\]
The first condition means that any gauge transformation respects the fibres of $P$; the second one means that $\sigma$ is equivariant w.r.t.\ the (right) $G$-action on $P$.
Later, we will see that there is another way of defining gauge transformations of a bundle $P$, but, for the moment, let me skip the problem.

Another problem I will address to later is that, in fact, any gauge transformation is an {\em isomorphism}; the proof of this fact is easy, but I prefer to postpone it to Section~\ref{sec-isombun}, deserving it a treatment in a more general context.

Let me notice that this last fact that means that the set of gauge transformations is a group w.r.t.\ the product operation given by composition; hence, it makes sense to speak of the {\em gauge group of a principal bundle $P$}, which I denote by $\calG_P$.

Finally, the gauge group $\calG_P$ operates on the space of connections on $P$, which I denote by $\calA=\calA_P$, to make manifest the dependence on the chosen bundle $P$; the (right) action is given explicitly by\[
A^\sigma\colon=\sigma^*A,\quad A\in\calA_P,\sigma\in \calG_P.
\]

\section{Holonomy and parallel transport}\label{sec-holpartransp}
Let me consider a curve $\gamma$ on $M$; by the word ``curve'', I mean in this context (if not otherwise stated) a {\em piecewise smooth map from the unit interval $\unint$ to $M$}.
\begin{Def}\label{eq-horlift}
Given a connection $1$-form $A$ on $P$, a {\em horizontal lift of $\gamma$ based at $p\in P$} is a curve on $P$, lying over $\gamma$, based at the point $p$ and such that all its tangent directions are ($A$-)horizontal.
\end{Def}

I quote without proof from \cite{KN1}, Chapter 2, Section 3, the following Theorem, which is the main ingredient of many of the subsequent constructions:
\begin{Thm}\label{thm-horlift}
Given a connection $1$-form $A$ on $P$ and a curve $\gamma$ in $M$, there is a unique horizontal lift of $\gamma$ based at the point $p$, which I denote by $\widetilde{\gamma}_{A,p}$.
\end{Thm}
Now, it is possible to display two important consequences of Theorem~\ref{thm-horlift}.

\subsection{Holonomy: definition and main properties}\label{ssec-holonom}
I consider a loop $\gamma$, i.e.\ a curve in $M$ satisfying $\gamma(0)=\gamma(1)$; I choose additionally a point $p\in P$ lying over $\gamma(0)=\gamma(1)$. 
Given a connection $A$, by Theorem~\ref{thm-horlift}, there is a unique horizontal lift $\widetilde{\gamma}_{A,p}$ of $\gamma$ based at $p$.
Since $\widetilde{\gamma}_{A,p}$ is a lift of $\gamma$, $\widetilde{\gamma}_{A,p}(1)$ also lies over $\gamma(0)$, and, as $G$ acts transitively on each fibre of $P$, it makes sense to propose the following
\begin{Def}\label{def-holonom}
The {\em holonomy of $\gamma$ w.r.t.\ the connection $A$ and base point $p$ over $\gamma(0)$} is the unique element of $G$, usually denoted by $\holon{A}{\gamma}{p}$, satisfying
\begin{equation}\label{eq-holonom}
\widetilde{\gamma}_{A,p}(1)=p\!\holon{A}{\gamma}{p},
\end{equation}
where $\widetilde{\gamma}_{A,p}$ is the unique horizontal lift based at $p\in P$ over $\gamma(0)$, of the loop $\gamma$. 
\end{Def}
The structure group $G$ and the gauge group $\calG_P$ operate from the right on $P$, resp.\ on the space of connections $\calA$. 
The holonomy depends by its very construction on a given loop $\gamma$, a connection $A$ on $P$ and a base point $p\in P_{\gamma(0)}$; hence, it makes sense to consider it as a function on the cartesian product of the space of connections $\calA$ on $P$, the space of loops with values in the structure group $G$ and the bundle $P$ itself; later, we will see that it is more precisely a section of some principal bundle over it.

The next Lemma shows how the holonomy behaves w.r.t.\ the action of $G$, resp.\ of $\calG$, on $P$, resp.\ $\calA$.

\begin{Lem}\label{lem-propholonom}
I assume $g$, resp.\ $\sigma$, to be an element of $G$, resp.\ a gauge transformation; I denote by $g_\sigma$ the function from $P$ with values in $G$ canonically associated to $\sigma$ via the formula 
\[
\sigma(p)=p g_\sigma(p).
\]
(This formula makes sense by the transitivity of the action of $G$ on any fibre of $P$.)

Then , the following formulae hold:
\begin{align}
\label{eq-liehol} \holon{A}{\gamma}{pg}&=\conj\!\left(g^{-1}\right)\holon{A}{\gamma}{p};\\ 
\label{eq-gaugehol} \holon{A^\sigma}{\gamma}{p}&=\conj\!\left(g_{\sigma}(p)^{-1}\right)\holon{A}{\gamma}{p}.
\end{align}
\end{Lem}
\begin{proof}
By Definition~\ref{def-holonom}, the identity holds
\begin{align*}
\widetilde{\gamma}_{A,pg}(1)=pg\ \holon{A}{\gamma}{pg},
\end{align*}
where $\widetilde{\gamma}_{A,pg}$ is the unique horizontal lift of $\gamma$ based at $pg$ given by Theorem~\ref{thm-horlift}. 
I consider the curve $\raction_g\left(\widetilde{\gamma}_{A,p}\right)$: it is clearly based at $pg$, and (\ref{eq-equiv}) implies that
\begin{align*}
A_{\raction_g\left(\widetilde{\gamma}_{A,p}(t)\right)}\left(\ddtt \raction_g\left(\widetilde{\gamma}_{A,p}(t)\right)\right)&=A_{\raction_g\left(\widetilde{\gamma}_{A,p}(t)\right)}\left[\tangent{\widetilde{\gamma}_{A,p}(t)}{\raction_g}\left(\ddtt \widetilde{\gamma}_{A,p}(t)\right)\right]=\\
&=\Ad\left(g^{-1}\right) \left[A_{\widetilde{\gamma}_{A,p}(t)}\left(\ddtt \widetilde{\gamma}_{A,p}(t)\right)\right]=0.
\end{align*}
Hence, the curve $\raction_g\left(\widetilde{\gamma}_{A,p}\right)$ is also horizontal, it lies over $\gamma$ (by the $G$-invariance of $\pi$) and is based at $pg$. 
By the uniqueness in Theorem~\ref{thm-horlift}, it follows $\raction_g\left(\widetilde{\gamma}_{A,p}\right)=\widetilde{\gamma}_{A,pg}$, whence
\[
\widetilde{\gamma}_{A,pg}(1)=\widetilde{\gamma}_{A,p}(1)g=p\ \holon{A}{\gamma}{p} g.
\]
Finally, the freeness of the action of $G$ implies (\ref{eq-liehol}).

Similarly, again by Definition~\ref{def-holonom}, it holds
\[
\widetilde{\gamma}_{A^\sigma,p}(1)=p\ \holon{A^\sigma}{\gamma}{p}.
\]
I claim now 
\begin{equation}\label{eq-gaugehor}
\sigma\left(\widetilde{\gamma}_{A^\sigma,p}\right)=\widetilde{\gamma}_{A,\sigma(p)}.
\end{equation}
Both curves $\sigma\left(\widetilde{\gamma}_{A^\sigma,p}\right)$ and $\widetilde{\gamma}_{A,\sigma(p)}$ lie over $\gamma$, as $\sigma$ is a gauge transformation, and both are clearly based at $\sigma(p)$.
To prove equation (\ref{eq-gaugehor}), it suffices to show that both curves are $A$-horizontal by uniqueness of horizontal lifts.
A direct computation shows:
\begin{align*}
A_{\sigma\left(\widetilde{\gamma}_{A^\sigma,p}(t)\right)}\left(\ddtt \sigma\left(\widetilde{\gamma}_{A^\sigma,p}(t)\right) \right)&=\left(A^\sigma\right)_{\widetilde{\gamma}_{A^\sigma,p}(t)}\left(\ddtt \widetilde{\gamma}_{A^\sigma,p}(t)\right)=\\
&=0,
\end{align*}
where the last identity is a consequence of the $A^\sigma$-horizontality of $\widetilde{\gamma}_{A^\sigma,p}$.
Hence, Identity (\ref{eq-gaugehor}) holds true.
Therefore, one gets
\[
\sigma(p)\ \holon{A^\sigma}{\gamma}{p}\overset{\text{By (\ref{eq-gaugehor})}}=\sigma(p)\ \holon{A}{\gamma}{\sigma(p)}
\]
Since the action of $G$ is free, it follows
\begin{equation*}
\begin{aligned}
\holon{A^\sigma}{\gamma}{p}&=\holon{A}{\gamma}{\sigma(p)}=\holon{A}{\gamma}{p g_\sigma(p)}\overset{\text{By (\ref{eq-liehol})}}=\\
&=\conj\left(g_\sigma(p)^{-1}\right)\ \holon{A}{\gamma}{p}.
\end{aligned}
\end{equation*}
\end{proof}

\subsection{Parallel transport: definition and main properties}\label{ssec-partransp}
The next object I want to define is the parallel transport w.r.t.\ a connection $A$ along a general curve $\gamma$.
\begin{Def}\label{def-partransp}
Let $\gamma$ be a curve in $M$, not necessarily closed, and let $t\in I$ and $p\in P$, resp.\ $q\in P$, be such that
\[
\pi(p)=\gamma(0),\quad\text{resp.}\quad \pi\left(q\right)=\gamma(t).
\] 
I define the {\em parallel transport from $p$ to $q$ along $\gamma$ from $0$ to $t$ w.r.t.\ $A$}, usually denoted by $\holonom{A}{\gamma}{t}{p}{q}$, as the unique element of $G$ obeying the rule
\begin{equation}\label{eq-partransp}
\widetilde{\gamma}_{A,p}(t)=q\ \holonom{A}{\gamma}{t}{p}{q},
\end{equation}
where $\widetilde{\gamma}_{A,p}$ is the horizontal lift of $\gamma$.
\end{Def}

The parallel transport defined by equation (\ref{eq-partransp}) satisfies two identities similar in spirit to (\ref{eq-liehol}) and (\ref{eq-gaugehol}).
\begin{Lem}\label{lem-proppartr}
Given a connection $A$ on $P$, a curve $\gamma$ in $M$, $p\in P_{\gamma(0)}$ and $q\in P_{\gamma(t)}$, for some $t$ in the unit interval, general elements $h$ and $k$ of $G$ and a gauge transformation $\sigma\in \calG$, the parallel transport of $\gamma$ from $p$ to $q$ along $\gamma$ from $0$ to $t$ satisfies the two following identities:
\begin{align}
\label{eq-liepartr} \holonom{A}{\gamma}{t}{pg}{qh}&=h^{-1} \holonom{A}{\gamma}{t}{p}{q} g;\\
\label{eq-gaugepartr} \holonom{A^\sigma}{\gamma}{t}{p}{q}&=g_{\sigma}\!\left(q\right)^{-1} \holonom{A}{\gamma}{t}{p}{q} g_\sigma(p).
\end{align}
\end{Lem} 
\begin{proof}
By Definition~\ref{def-partransp}, one has
\[
\widetilde{\gamma}_{A,p}(t)=q h\ \holonom{A}{\gamma}{t}{p}{q h}=q\ \holonom{A}{\gamma}{t}{p}{q}.
\]
By the freeness of the action of $G$ on $P$, one gets
\begin{equation}\label{eq-liepartr1}
\holonom{A}{\gamma}{t}{p}{q h}=h^{-1} \holonom{A}{\gamma}{t}{p}{q},\quad \forall h\in G.
\end{equation}
On the other hand, the identity holds
\begin{align*}
\widetilde{\gamma}_{A,pg}(t)&=q\ \holonom{A}{\gamma}{t}{pg}{q}=\\
&=\widetilde{\gamma}_{A,p}(t) g=\\
&=q\ \holonom{A}{\gamma}{t}{p}{q} g,
\end{align*}
where the second identity was shown in the proof of Lemma~\ref{lem-propholonom}.
It follows therefore
\begin{equation}\label{eq-liepartr2}
\holonom{A}{\gamma}{t}{pg}{q}=\holonom{A}{\gamma}{t}{p}{q} g,\quad \forall g\in G.
\end{equation}
Combining (\ref{eq-liepartr1}) and (\ref{eq-liepartr2}), one gets (\ref{eq-liepartr}).

As for the second identity, I make use again of (\ref{eq-gaugehor}): it then holds
\begin{align*}
\widetilde{\gamma}_{A,\sigma(p)}\!(t)&=\sigma\!\left(q\right) \holonom{A}{\gamma}{t}{\sigma(p)}{\sigma\left(q\right)}=\\
&=\sigma\!\left(\widetilde{\gamma}_{A^\sigma,p}(t)\right)=\\
&=\sigma\left(q\holonom{A^\sigma}{\gamma}{t}{p}{q}\right)=\\
&=\sigma\left(q\right) \holonom{A}{\gamma}{t}{\sigma(p)}{\sigma\left(q\right)}.
\end{align*}
The freeness of the action of $G$ implies then
\begin{align*}
\holonom{A^\sigma}{\gamma}{t}{p}{q}&=\holonom{A}{\gamma}{t}{\sigma(p)}{\sigma\left(q\right)}=\\
&=\holonom{A}{\gamma}{t}{p g_{\sigma(p)}}{q g_{\sigma}\left(q\right)}=\\
&=g_\sigma\!\left(q\right)^{-1} \holonom{A}{\gamma}{t}{p}{q} g_\sigma(p),
\end{align*}
where the last identity is a consequence of (\ref{eq-liepartr}).
\end{proof}

After having introduced the parallel transport of a curve $\gamma$ from $0$ to a point $t$ w.r.t.\ $A$ and having discussed some of its properties, I may also introduce another object, namely the parallel transport of the curve $\gamma$ from $s$ to $t$, where $s,t\in\unint$ satisfy $s<t$.

First of all, given a curve $\gamma$ from the unit interval $\unint$ to $M$, and given $s\in\unint$, let me define a new curve $\gamma_s$ from the interval $\left[0,1-s\right]$ by the assignment
\[
\gamma_s\!(t)\colon=\gamma(t+s),\quad \forall t\in\left[0,1-s\right].
\]

\begin{Def}
Let $\gamma$ be a curve in $M$, defined on the unit interval $\unint$; let then $s<t$ two points in $\unint$, and $p$ and $q$ two points in $P$ satisfying
\[
\pi(p)=\gamma(s),\quad \pi(q)=\gamma(t).
\]

The {\em parallel transport of $\gamma$ from $p$ to $q$ w.r.t.\ the connection $A$}, which I denote by $\holonomy{A}{\gamma}{s}{t}{p}{q}$, is the unique element of $G$ obeying the rule
\begin{equation}\label{eq-partrtransl}
\widetilde{\gamma_s}_{A,p}(t-s)=q\holonomy{A}{\gamma}{s}{t}{p}{q}.
\end{equation}
Recalling Equation (\ref{eq-partransp}), Equation (\ref{eq-partrtransl}) is equivalent to 
\begin{equation}\label{eq-partrtransl1}
\holonomy{A}{\gamma}{s}{t}{p}{q}=\holonom{A}{\gamma_s}{t-s}{p}{q}.
\end{equation}
\end{Def}

Lemma~\ref{lem-proppartr} together with Equation (\ref{eq-partrtransl1}) implies the useful
\begin{Cor}\label{cor-proppartr}
Given a connection $A$ on $P$, a curve $\gamma$ in $M$, two points $s,t\in\unint$ obeying $s<t$, $p\in P_{\gamma(s)}$ and $q\in P_{\gamma(t)}$, general elements $h$ and $k$ of $G$ and a gauge transformation $\sigma\in \calG$, the parallel transport of $\gamma$ from $p$ to $q$ along $\gamma$ from $s$ to $t$ satisfies the following identities:
\begin{align*}
\holonomy{A}{\gamma}{s}{t}{pg}{qh}&=h^{-1} \holonomy{A}{\gamma}{s}{t}{p}{q} g;\\
\holonomy{A^\sigma}{\gamma}{s}{t}{p}{q}&=g_{\sigma}\!\left(q\right)^{-1} \holonomy{A}{\gamma}{s}{t}{p}{q} g_\sigma(p).
\end{align*}
\end{Cor}

\section{Isomorphisms of principal bundles and generalized gauge transformations}\label{sec-isombun}
Let $P\overset{\pi}\to M$, $\widetilde{P}\overset{\widetilde{\pi}}\to M$ two principal bundles over the same manifold $M$ and with the same structure group $G$; the right action of $G$ on $P$, resp.\ $\widetilde{P}$, is denoted by $\raction_{\bullet}$, resp.\ $\widetilde{\raction}_{\bullet}$, or simply by $(p,g)\mapsto pg$, resp.\ $(\tildep,g)\mapsto \tildep g$, given the case.

\begin{Def}\label{def-bunisom}
Given two principal bundles $P$, $\widetilde{P}$ as above, a {\em $G$-equivariant map $\sigma$ from $P$ to $\widetilde{P}$} is a smooth map from $P$ to $\widetilde{P}$ satisfying the two properties:
\begin{equation}\label{eq-princbunisom}
\begin{aligned}
\widetilde{\pi}\circ \sigma&=\pi,\\
\sigma\circ\raction_g&=\widetilde{\raction}_g\circ\sigma,\quad \forall g\in G.
\end{aligned}
\end{equation}
The set of all such maps is denoted by $\calG_{P,\widetilde{P}}$.
\end{Def}

\begin{Rem}
The set $\calG_{P,P}$, for a given principal bundle $P$, is the set of gauge transformations of $P$.
\end{Rem}

Given two principal bundles $P$, $\widetilde{P}$ over $M$ as above, it is possible to form out of them a manifold as follows:
\begin{Def}\label{def-fibprod}
Given two principal $G$-bundles $P$ and $\widetilde{P}$ over the same base space $M$, their {\em fibred product}, denoted usually by $P\odot \widetilde{P}$, is defined as
\begin{equation*}
P\odot \widetilde{P}\colon=\left\{\left(p,\tildep\right)\in P\times \widetilde{P}\colon \pi(p)=\widetilde{\pi}\left(\tildep\right)\right\}.
\end{equation*}
\end{Def}
There is a natural map $\overline{\pi}$ from the fibred product $P\odot \widetilde{P}$ to $M$, which is simply
\begin{equation*}
\overline{\pi}\left(p,\widetilde{p}\right)\colon=\pi(p)=\widetilde{\pi}\left(\tildep\right),\quad \left(p,\tildep\right)\in P\odot \widetilde{P}.
\end{equation*}
Additionally, $P\odot \widetilde{P}$ receives a right $G\times G$-action:
\[
\left(p,\tildep;(g,h)\right)\mapsto \left(pg,\tildep h\right),\quad \forall \left(p,\tildep\right)\in P\odot \widetilde{P}, (g,h)\in G\times G.
\]
It is clear that the above action is free, as both actions of $G$ on $P$ and $\widetilde{P}$ are free; moreover, considering the fibre
\[
\left(P\odot\widetilde{P}\right)_x\colon=\overline{\pi}^{-1}\left(\left\{x\right\}\right),
\] 
it follows immediately that the action of $G\times G$ is transitive on it.

These two facts are not incidental, because of the following 
\begin{Prop}\label{prop-fibprod}
The fibred product $P\odot \widetilde{P}$ is a principal $G\times G$-bundle over $M$, with projection $\overline{\pi}$.
\end{Prop}
\begin{proof}
If $U$ is an open set of $M$, let me denote by $\varphi_U$, resp.\ $\widetilde{\varphi}_U$, the trivialization of $P$ over $U$, resp.\ of $\widetilde{P}$ over $U$.
A trivialization $\overline{\varphi}_U$ of $P\odot\widetilde{P}$ over $U$ may thus be defined via
\begin{equation*}
\begin{aligned}
\overline{\varphi}\colon \overline{\pi}^{-1}\!\left(U\right)&\longrightarrow U\times \left(G\times G\right)\\
\left(p,\tildep\right)&\longmapsto \left(\pi(p);\left(\pr_{2}\circ \varphi_U\right)(p),\left(\pr_{2}\circ \widetilde{\varphi}_U\right) \left(\tildep\right)\right).\,
\end{aligned}
\end{equation*}
where $\pr_2$ denotes the projection from $U\times G$ onto $G$. 
These maps are invertible, their inverses being given by
\begin{equation*}
\begin{aligned}
\overline{\varphi}_U^{-1}\colon U\times \left(G\times G\right)&\longrightarrow \overline{\pi}^{-1}\!\left(U\right)\\
\left(x;g,h\right)&\longmapsto \left(\varphi_U^{-1}\left(x,g\right),\widetilde{\varphi}_U^{-1}\left(x,h\right)\right).
\end{aligned}
\end{equation*}
It is clear from their definition that the maps $\overline{\varphi}_U$ and their inverses are smooth.
Hence, I have obtained a trivialization of the fibred product $P\odot \widetilde{P}$.

For the sake of completeness, let me write down explicitly the transition maps of the $G\times G$-principal bundle $P\odot \widetilde{P}$:
\begin{equation*}
\begin{aligned}
\overline{\varphi}_{U,V}\colon U\cap V&\longrightarrow \diff\!(G\times G)\\
x&\longmapsto \raction_{\varphi_{U,V}(x)}\times \widetilde{\raction}_{\widetilde{\varphi}_{U,V}(x)},
\end{aligned}
\end{equation*}
where $\varphi_{U,V}$, resp.\ $\widetilde{\varphi}_{U,V}$, are the transition maps of $P$, resp.\ $\widetilde{P}$ w.r.t.\ the trivializations $\varphi_U$, $\varphi_V$, resp.\ $\widetilde{\varphi}_U$, $\widetilde{\varphi}_V$, for $U$, $V$ any two open subsets of $M$ with nontrivial intersection.
\end{proof}

\begin{Rem}
Let me notice that the fibred product may be also seen as a principal $G$-bundle over $P$: in fact, by its very definition,
\[
P\odot\widetilde{P}=\pi^*\!(\widetilde{P}),
\]
and the latter manifold inherits clearly a principal $G$-bundle structure over $P$.
(Equivalently, the fibred product may be also viewed as a $G$-bundle over $\widetilde{P}$, as $P\odot \widetilde{P}=\widetilde{\pi}^*\!(P)$.)
\end{Rem}

\begin{Rem}
The fibred product of two principal bundles over the same base space and with the same structure group may be seen as an analogue of the Whitney sum of vector bundles for principal bundles.

Furthermore, it is clear that there is a canonical isomorphism
\[
P\odot \widetilde{P}\cong\widetilde{P}\odot P,
\]
for any two principal bundles $P$, $\widetilde{P}$ over the same base space $M$.
\end{Rem}

Now, there is a left action of $G\times G$ on $G$ specified by the rule:
\begin{equation*}
\begin{aligned}
\conj\colon G\times G&\longrightarrow \diff\!(G)\\
(g,h)&\longmapsto \conj\!\left(g,h\right)k\colon= h k g^{-1}.
\end{aligned}
\end{equation*}
It is quite evident that the restriction of $\conj$ to the diagonal subgroup $G$ of $G\times G$ gives the usual conjugation of $G$ on itself; therefore, one can speak of the above action as of the {\em generalized conjugation in $G$}. 

\begin{Def}\label{def-equivmaps}
Under the same hypotheses as in Definition~\ref{def-fibprod}, the set of {\em smooth $G\times G$-equivariant maps from $P\odot \widetilde{P}$ to $G$}, denoted by $C^{\infty}\!\!\left(P\odot \widetilde{P},G\right)^{G\times G}$, is the subset of $C^{\infty}\!\!\left(P\odot \widetilde{P},G\right)$ of those maps $K$ satisfying the equivariance w.r.t.\ the generalized conjugation 
\begin{equation}\label{def-eqisom}
K\!\left(pg,\tildep h\right)=\conj\!\left(g^{-1},h^{-1}\right)K\!\left(p,\tildep\right),\quad \forall \left(p,\tildep\right)\in P\odot \widetilde{P}, (g,h)\in G\times G.
\end{equation} 
\end{Def}

To a given principal $G$-bundle $P$ over $M$, one can associate the canonical map $\phi_P$ on $P\odot P$ with values in $G$ by the rule
\begin{equation}\label{eq-canbundle}
q=p \phi_P\!(p,q),
\end{equation}
for any pair $(p,q)$ in $P\odot P$. 
Notice that $i)$ the definition makes sense, because $p$ and $q$ lie in the same fibre and $G$ operates transitively on each fibre, and $ii)$ the element $\phi_P\!(p,q)$ is uniquely defined since the action of $G$ on $P$ is free.

The following proposition illustrates the main properties of $\phi_P$.
\begin{Prop}\label{prop-canprop}
The map $\phi_P$ defined by (\ref{eq-canbundle}) is smooth, and satisfies the two identities
\begin{align}
\label{eq-diagonal} \phi_P\!(p,p)&=e,\forall p\in P;\\
\label{eq-inverse} \phi_P\!\left(p,q\right)&=\phi_P\!\left(q,p\right)^{-1},\quad \forall p,q\in P\odot P;\\
\label{eq-canequiv} \phi_P\!(pg,qh)&=g^{-1}\phi_P\!(p,q) h,\quad \forall (p,q)\in P\odot P,\quad \forall (g,h)\in G\times G.
\end{align}
\end{Prop}
\begin{proof}
By considering a local trivialization $\overline{\varphi}_U$ of $P\odot P$ over $U\subset M$ open, one can find a local expression for the map $\phi_P$:
\[
(x;g,h)\overset{\phi_{P,U}}\leadsto g^{-1}h, 
\]
whence it follows immediately that $\phi_P$ is a smooth map.

Moreover, by its very definition, $\phi_P$ is uniquely determined by (\ref{eq-canbundle}).
Taking therefore the pair $(p,p)$ in $P\odot P$, by the freeness of the action of $G$ on $P$, both (\ref{eq-diagonal}) and (\ref{eq-inverse}) follow immediately.

On the other hand, again by (\ref{eq-canbundle}), one gets
\begin{align*}
\widetilde{p}h&=pg\phi_P\!(pg,\widetilde{p}h)=\\
&=p\phi_{P}\!(p,\widetilde{p})h,
\end{align*} 
whence (\ref{eq-canequiv}) follows again by the freeness of the action of $G$ on any fibre.
\end{proof}
\begin{Rem}
The canonical map $\phi_P$ was called by MacKenzie in~\cite{McK} the {\em division map of $P$}; the origin of the name lies obviously in its local description, or, equivalently, in its shape for the trivial bundle.
On the other hand, M{\oe}rdijk~\cite{Moer} called the division map $\phi_P$ (although working in the more general context of principal bundles with structure groupoids) a {\em cocycle on $P$ with values in $G$}; I will use the former name for $\phi_P$.
\end{Rem}

I have all elements now to state and prove the following Theorem, which gives an alternative characterisation of the maps introduced in Definition~\ref{def-bunisom}
\begin{Thm}\label{thm-isomequiv}
There is a one-to-one correspondence between the set $\calG_{P,\widetilde{P}}$ of bundle morphisms from $P$ to $\widetilde{P}$, as in Definition~\ref{def-bunisom}, and the set $C^{\infty}\!\!\left(P\odot \widetilde{P},G\right)^{G\times G}$ of smooth $G\times G$-equivariant maps from the fibred product $P\odot \widetilde{P}$ to $G$, as in Definition~\ref{def-equivmaps}.
\end{Thm}
\begin{proof}
Consider a bundle morphism $\sigma$ from $P$ to $\widetilde{P}$ and take a pair $\left(p,\tildep\right)\in P\odot \widetilde{P}$.
By the first property in (\ref{eq-princbunisom}), it follows that the pair $(\tildep,\sigma(p))$ belongs to the fibred product of $\widetilde{P}$ with itself. 
Therefore, one can consider the function on $P\odot \widetilde{P}$ with values in $G$ given by
\[
K_\sigma\!\left(p,\tildep\right)\colon=\phi_{\widetilde{P}}\!\left(\tildep,\sigma(p)\right).
\]
By its very definition, $K_\sigma\!\left(p,\tildep\right)$ is uniquely determined by the equation
\[
\sigma(p)=\tildep K_\sigma\!\left(p,\tildep\right).
\]
By Proposition~\ref{prop-canprop}, it follows immediately that $K_\sigma$ is smooth, since $\phi_{\widetilde{P}}$ is smooth and $\sigma$ also.
Equation (\ref{eq-canequiv}) implies
\begin{align*}
K_\sigma\!\left(pg,\tildep h\right)&=\phi_{\widetilde{P}}\!\left(\tildep h,\sigma(pg)\right)=\\
&\overset{\text{by (\ref{eq-princbunisom})}}=\phi_{\widetilde{P}}\!\left(\tildep h,\sigma(p)g\right)=\\
&=h^{-1}\phi_{\widetilde{P}}\!\left(\tildep,\sigma(p)\right)g=\\
&=h^{-1} K_\sigma\!\left(p,\tildep\right)g.
\end{align*}
Conversely, given an element $K$ of $C^{\infty}\!\left(P\odot \widetilde{P},G\right)^{G\times G}$, it is possible to construct a map $\sigma_K$ from $P$ to $\widetilde{P}$ as follows:
\[
\sigma_K\!\left(p\right)\colon=\tildep\ K\!\left(p,\tildep\right),
\]
$\tildep$ being any element of $\widetilde{P}$, such that the pair $\left(p,\tildep\right)$ lies in $P\odot \widetilde{P}$.

First of all, one has to show that the map $\sigma_K$ is well-defined: taking $\widetilde{q}$ to be another point in $\widetilde{P}$, such that $\widetilde{\pi}\left(\tildep\right)=\widetilde{\pi}\left(\widetilde{q}\right)$.
Since $\widetilde{p}$ and $\widetilde{q}$ belong to the same fibre, the identity holds
\[
\widetilde{q}=\tildep \phi_{\widetilde{P}}\!\left(\tildep,\widetilde{q}\right).
\]
It then follows
\begin{align*}
\sigma_K\left(p\right)&=\widetilde{q}\ K\!\left(p,\widetilde{q}\right)=\\
&=\tildep \phi_{\widetilde{P}}\!\left(\tildep,\widetilde{q}\right) \ K\!\left(p,\tildep\phi_{\widetilde{P}}\!\left(\tildep,\widetilde{q}\right)\right)=\\
&\overset{\text{By equivariance}}=\tildep\ K\!\left(p,\tildep\right).
\end{align*}
Therefore, $\sigma_K$ is well-defined.
It also satisfies $\widetilde{\pi}\circ \sigma_K=\pi$ by its very definition, and the $G$-equivariance follows from
\begin{align*}
\sigma_K(p g)&=\tildep K\!\left(pg,\tildep\right)=\\
&\overset{\text{By equivariance}}=\tildep K\!\left(p,\tildep\right)\ g=\\
&=\sigma_K(p) g,\quad \forall g\in G,
\end{align*}
hence $\sigma_K$ is a bundle morphism from $P$ to $\widetilde{P}$.
It remains to show that the assignments 
\[
\sigma\leadsto K_\sigma\quad\text{and}\quad K\leadsto \sigma_K 
\]
are inverse to each other.
Namely, one has
\begin{align*}
\sigma_{K_\sigma}(p)&=\tildep K_\sigma\!\left(p,\tildep\right)=\\
&=\tildep \phi_{\widetilde{P}}\!\left(\tildep, \sigma(p)\right)=\\
&\overset{\text{by (\ref{eq-canbundle})}}=\sigma(p),
\end{align*}
where $\tildep\in \widetilde{P}$ is chosen so that the pair $(p,\tildep)\in P\odot \widetilde{P}$.
On the other hand,
\begin{align*}
K_{\sigma_K}\!\left(p,\tildep\right)&=\phi_{\widetilde{P}}\!\left(\tildep,\sigma_K(p)\right)=\\
&=\phi_{\widetilde{P}}\!\left(\tildep,\tildep K\!\left(p,\tildep\right)\right)=\\
&\overset{\text{by (\ref{eq-canequiv})}}=\phi_{\widetilde{P}}\!\left(\tildep,\tildep\right)K\!\left(p,\tildep\right)=\\
&\overset{\text{by (\ref{eq-diagonal})}}=K\!\left(p,\tildep\right).
\end{align*}
\end{proof}
Notice that the bundle morphism $\sigma_K$ canonically associated to $K\in C^{\infty}\!\left(P\odot \widetilde{P},G\right)^{G\times G}$ is invertible; namely, an explicit inverse map can be defined via
\[
\sigma_K^{-1}(\tildep)\colon=p\ K\!\left(p,\tildep\right)^{-1},
\]
where $p\in P$ is chosen so that the pair $(p,\tildep)$ belongs to $P\odot \widetilde{P}$, and $K^{-1}$ denotes the inverse in $G$ of $K$.
As in the proof of Theorem~\ref{thm-isomequiv}, one can show that $\sigma_K^{-1}$ is well-defined, that it is $G$-equivariant and satisfies $\pi\circ \sigma_K^{-1}=\widetilde{\pi}$; in fact,
\begin{align*}
\left(\sigma_K^{-1}\circ \sigma_K\right)(p)&=\sigma_K^{-1}\!\left(\tildep K\!\left(p,\tildep\right)\right)=\\
&=\sigma_K^{-1}\! \left(\tildep\right) K\!\left(p,\tildep\right)=\\
&=p K\!\left(p,\tildep\right)^{-1} \!K\!\left(p,\tildep\right)=\\
&=p.
\end{align*}
It follows analogously that $\sigma_K\circ \sigma_K^{-1}=\id_{\widetilde{P}}$.
Notice that I made use of the fact that both definitions of $\sigma_K$ and $\sigma_K^{-1}$ do not depend on the choice of points in $\widetilde{P}$, resp.\ $P$.

This fact finds its motivation in the following
\begin{Lem}\label{lem-gaugeisom}
Any element $\sigma$ of $\calG_{P,\widetilde{P}}$ is a bijection.
\end{Lem}
\begin{proof}
First of all, any $\sigma$ in $\calG_{P,\widetilde{P}}$ is injective: in fact, assuming 
\[
\sigma(p)=\sigma(q)
\] 
for $p$ and $q$ in $P$, one gets immediately
\[
(\widetilde{\pi}\circ \sigma)(p)=\pi(p)=\pi(q)=(\widetilde{\pi}\circ \sigma)(q)\Rightarrow q=p \phi_P\!\left(p,q\right).
\]
Hence, it follows
\begin{align*}
\sigma(q)&=\sigma(p \phi_P\!\left(p,q\right))=\\
&=\sigma(p)\phi_P\!\left(p,q\right)=\\
&\overset{!}=\sigma(p).
\end{align*}
Therefore, by the freeness of the action of $G$, it follows
\[
\phi_P\!\left(p,q\right)=e\Rightarrow p=q,
\]
yielding injectivity of $\sigma$.

It remains to prove surjectivity of $\sigma$; for this purpose, it suffices to show that $\sigma$ is surjective on any fibre, since the equation $\sigma(p)=\tildep$ implies, by the first equation in (\ref{eq-princbunisom}), that $(p,\tildep)\in P\odot \widetilde{P}$.
Hence, taking in any fibre $P_x$ of $\pi$ a point $p_x\in P$ (whose existence is guaranteed by the surjectivity of $\pi$), consider its image $\sigma(p_x)$.
Consequently, taking a general $\tildep\in\widetilde{P}$, consider $p_{\widetilde{\pi}(\tildep)}\in P$.
Defining 
\[
q\colon=p_{\widetilde{\pi}(\tildep)}\phi_{\widetilde{P}}\!\left(\sigma(p_{\widetilde{\pi}(\tildep)}),\tildep\right),
\]   
a direct computation shows
\begin{align*}
\sigma(q)&=\sigma(p_{\widetilde{\pi}(\tildep)}\phi_{\widetilde{P}}\!\left(\sigma(p_{\widetilde{\pi}(\tildep)}),\tildep\right))=\\
&\overset{\text{by equivariance}}=\sigma(p_{\widetilde{\pi}(\tildep)})\phi_{\widetilde{P}}\!\left(\sigma(p_{\widetilde{\pi}(\tildep)}),\tildep\right)=\\
&\overset{\text{by (\ref{eq-canbundle})}}=\tildep.
\end{align*}
\end{proof}
The inverse of any bundle morphism $\sigma\in\calG_{P,\widetilde{P}}$ is automatically smooth; this can be checked locally or can be viewed as a consequence of Theorem~\ref{thm-isomequiv}.

\begin{Rem}
If one considers the case $P=\widetilde{P}$, the proof of Lemma~\ref{lem-gaugeisom} implies the well-known fact that a gauge transformation, i.e.\ a bundle endomorphism of $P$, is invertible, whence it follows that the set of gauge transformations is a group w.r.t.\ composition.
\end{Rem}

\subsection{The groupoid of generalized gauge transformations}\label{ssec-groupoid}
I consider now the category $\mathsf{Bun}_{G,M}$ of (smooth) principal bundles over $M$ with structure group $G$.
In order to avoid cumbersome notations, a general object of the category $\mathsf{Bun}_{G,M}$ is labelled by $P_i$, for some index $i$; the surjective submersion from $P_i$ to $M$ is labelled by $\pi_i$, and general elements of $P_i$ are always labelled by the index $i$, e.g.\ $p_i$, $\tildep_i$, $q_i$.
Theorem~\ref{thm-isomequiv} is a general version of the well-known correspondence between bundle endomorphisms of $P$ (also called {\em gauge transformations}) and (smooth) maps on $P$ with values in $G$, equivariant w.r.t.\ the conjugation action of $G$ on itself.
Moreover, this correspondence can be proved to be an isomorphism of groups, because the set of $G$-equivariant maps from $P$ to $G$ w.r.t.\ conjugation inherits a group structure from the multiplication in $G$.
Obviously, such a group structure is in principle not anymore available, if one considers, instead of a fixed principal $G$-bundle $P$ over $M$ and its gauge group, all objects of $\mathsf{Bun}_{G,M}$ and its set of morphisms, i.e.\ $G$-equivariant, fibre-preserving maps from an object $P$ to another object $\widetilde{P}$, which must not be a priori identical.        
Lemma~\ref{lem-gaugeisom} implies that any morphism of the category $\mathsf{Bun}_{G,M}$ in the above sense is invertible; hence, the morphisms of the category $\mathsf{Bun}_{G,M}$ build a {\em groupoid} over the objects of the same category.

Let me deal with the following replacement:
\begin{align*}
P\in\mathsf{Ob}(\mathsf{Bun}_{G,M})&\leadsto \left(P_1,P_2\right)\in\mathsf{Ob}(\mathsf{Bun}_{G,M}^2),\\
C^{\infty}\!(P,G)^G&\leadsto C^{\infty}\!(P_1\odot P_2,G)^{G\times G}.
\end{align*}
(Here, I have denoted by $\mathsf{Bun}_{G,M}^2$ the product category of $\mathsf{Bun}_{G,M}$ with itself.)
It is clear that one can consider the ``diagonal'' in $\mathsf{Ob}\left(\mathsf{Bun}_{G,M}^2\right)$, consisting of pairs of the form $\left(P,P\right)$, for $P\in\mathsf{Bun}_{G,M}$.
The $G$-equivariant, fibre-preserving isomorphisms of $P$ (gauge transformations) correspond, on the one hand, to the group $C^{\infty}\!(P,G)^G$, on the other hand, by Theorem~\ref{thm-isomequiv}, to the set $C^{\infty}\!(P\odot P,G)^{G\times G}$.

It can be therefore expected that both sets are in bijection; this is in fact true in virtue of the following
\begin{Lem}\label{lem-gaugetrsf}
The set $C^{\infty}\!(P,G)^G$ of $G$-equivariant maps on $P$ with values in $G$ is canonically bijective to the set $C^{\infty}\!(P\odot P,G)^{G\times G}$ of $G\times G$-equivariant maps on the fibred product of $P$ with itself with values in $G$.
\end{Lem}
\begin{proof}
First of all, define a canonical map from $C^{\infty}\!(P\odot P,G)^{G\times G}$ to $C^{\infty}\!(P,G)$ via the assignment
\begin{equation*}
K\in C^{\infty}\!(P\odot P,G)^{G\times G}\leadsto C^{\infty}\!(P,G)\ni g_K\!(p)\colon=K\!(p,p),
\end{equation*}
i.e., consider simply the {\em restriction to the diagonal of $P\odot P$} of any element of $C^{\infty}\!(P\odot P,G)^{G\times G}$.
One has to prove that such a map $g_K$, which is smooth by construction, belongs to $C^{\infty}\!(P,G)^G$, i.e.\ one has to prove its $G$-equivariance.
This is an immediate consequence of the fact that the generalized conjugation of $G\times G$ restricts to the usual conjugation, when taking its restriction from $G\times G$ to its diagonal subgroup.

Consider on the other hand the map $\phi_P$ in (\ref{eq-canbundle}), and to a general element $k$ of $C^{\infty}\!(P,G)^G$, associate the map $K$ defined via the rule
\begin{equation}\label{eq-invertwist}
K\!(p,\widetilde{p})\colon=\phi_K\!\left(p,\widetilde{p}\right)^{-1} k(p),
\end{equation}
for any pair $(p,\widetilde{p})$ in $P\odot P$.
It is clear that the map $K$ in (\ref{eq-invertwist}) is smooth by Proposition~\ref{prop-canprop} and the smoothness of $k$; hence, it remains to show that it is $G\times G$-equivariant.
This is a consequence of the following computation, which follows again from Proposition~\ref{prop-canprop}:
\begin{align*}
K\!(pg,\widetilde{p}h)&=\phi_P\!(pg,\widetilde{p}h)^{-1}k(pg)=\\
&\overset{\text{By equivariance of $k$ and $\phi_K$}}=\left(g^{-1}\phi_K(p,\widetilde{p})h\right)^{-1}g^{-1}k(p)g=\\
&=h^{-1}\phi_K(p,\widetilde{p})^{-1}g g^{-1}k(p)g=\\
&=\widetilde{\conj}\!\left(g^{-1},h^{-1}\right)K\!(p,\widetilde{p}).
\end{align*}

It remains to show that the assignments
\[
k\mapsto K\quad \text{and}\quad K\mapsto g_K
\]
are inverse to each other.
In fact, one gets
\begin{align*}
&k(p)\leadsto K(p,\widetilde{p})=\phi_K\!(p,\widetilde{p})^{-1}k(p)\leadsto \\
&\leadsto g_K(p)=K(p,p)=\phi_K\!(p,p)^{-1}k(p)\overset{\text{by the properties of $\phi_K$}}=k(p).
\end{align*}
On the other hand, 
\begin{align*}
&K\!(p,\tildep)\leadsto g_K(p)=K\!(p,p)\leadsto \phi_K\!(p,\widetilde{p})^{-1}g_K(p)=\phi_K\!(p,\widetilde{p})^{-1}K\!(p,p)=\\
&\overset{\text{by $G\times G$-equivariance of $K$}}=K\!(p,p\phi_K\!(p,\widetilde{p}))\overset{\text{By the definition of $\phi_K$}}=K(p,\widetilde{p}).
\end{align*}
\end{proof}
\begin{Rem}
It is easy to prove that the inverse in $G$ of the map $\phi_P$ yields an element of $C^{\infty}\!(P\odot P,G)^{G\times G}$, which, by Theorem~\ref{thm-isomequiv}, corresponds uniquely to a fibre-preserving isomorphism of $P$.
By the very arguments of the proof of Lemma~\ref{lem-gaugetrsf}, it follows immediately that $\phi_P$ corresponds to the identity isomorphism of $P$; namely,
\[
\id_P\leadsto g_{\id_P}(p)=e,\quad \forall p\in P\leadsto \phi_P\!(p,\tildep)^{-1}e=\phi_P\!(p,\tildep)^{-1},\quad \forall (p,\tildep)\in P\odot P.
\]
\end{Rem}

Motivated by the fact that the correspondence $\calG_P\leftrightarrow C^{\infty}\!(P,G)^G$ is an isomorphism of groups, I construct a ``product'' on any set of morphisms from $P_1$ to $P_2$ in the category $\mathsf{Bun}_{G,M}$.

Let me just begin with a notational remark:
\fbox{\parbox{12cm}{\bf A general bundle morphism from $P_i$ to $P_j$ is denoted by $\sigma_{ij}$.}}
Given three principal bundles $P_1$, $P_2$ and $P_3$, and bundle morphisms $\sigma_{12}$ and $\sigma_{23}$, since $\sigma_{23}\circ\sigma_{12}$ is obviously a bundle morphism from $P_1$ to $P_3$, there is a unique element in $C^\infty\!(P_1\odot P_3,G)^{G\times G}$, given by
\[
\left(\sigma_{23}\circ\sigma_{12}\right)\leadsto K_{\sigma_{23}\circ\sigma_{12}}\!\left(p_1,p_3\right)=\phi_{P_3}\!\left(p_3,\left(\sigma_{23}\circ\sigma_{12}\right)(p_1)\right).
\]
A direct computation gives
\begin{align*}
\left(\sigma_{23}\circ\sigma_{12}\right)(p_1)&=\sigma_{23}\left(p_2K_{\sigma_{12}}\!\left(p_1,p_2\right)\right)=\\
&=\sigma_{23}(p_2)K_{\sigma_{12}}\!\left(p_1,p_2\right)=\\
&=p_3 K_{\sigma_{23}}\!\left(p_2,p_3\right)K_{\sigma_{12}}\!\left(p_1,p_2\right),
\end{align*}
where $p_2\in P_2$ such that $\pi_2(p_2)=\pi_1(p_1)=\pi_3(p_3)$.
The freeness of the action of $G$ on $P_3$ and (\ref{eq-diagonal}) imply finally
\begin{equation}\label{eq-composisom}
K_{13}\!\left(p_1,p_3\right)\colon=K_{\sigma_{23}\circ\sigma_{12}}\!\left(p_1,p_3\right)=K_{\sigma_{23}}\!\left(p_2,p_3\right)K_{\sigma_{12}}\!\left(p_1,p_2\right).
\end{equation}
(In order to avoid cumbersome notations, I simply abbreviate $K_{\sigma_{12}}$ by $K_{12}$ and so on.)

Let me make some comments on Equation (\ref{eq-composisom}).
Since $K_{12}$ belongs to $C^\infty\!(P_1\odot P_2,G)^{G\times G}$ and $K_{23}$ to $C^\infty\!(P_2\odot P_3,G)^{G\times G}$, $K_{13}$ as defined in (\ref{eq-composisom}) does not depend on the choice of $p_2\in P_2$ such that $\pi_1(p_1)=\pi_2(p_2)=\pi_3(p_3)$: namely, for another representative $q_2=p_2 \phi_{P_2}\!\left(p_2,q_2\right)$, one gets
\begin{align*}
K_{13}\!(p_1,p_3)&=K_{23}\!(p_2\phi_{P_2}\!\left(p_2,q_2\right),p_3)K_{12}\!(p_1,p_2 \phi_{P_2}\!\left(p_2,q_2\right))=\\
&=K_{23}\!(p_2,p_3)\phi_{P_2}\!\left(p_2,q_2\right) \phi_{P_2}\!\left(p_2,q_2\right)^{-1}K_{12}\!(p_1,p_2)=\\
&=K_{23}\!(p_2,p_3)K_{12}\!(p_1,p_2).
\end{align*} 
On the other hand, $K_{13}$ as in (\ref{eq-composisom}) is also $G\times G$-equivariant w.r.t.\ the action $\widetilde{\conj}$
\begin{align*}
K_{13}\!(p_1 g,p_3 h)&=K_{23}\!(p_2,p_3 h)K_{12}\!(p_1 g,p_2)=\\
&=h^{-1}K_{23}\!(p_2,p_3)K_{12}\!(p_1,p_2)g=\\
&=\conj\!\left(g^{-1},h^{-1}\right)K_{13}\!(p_1,p_3),\quad \forall g,h\in G.
\end{align*}

For any two objects $P_1$, $P_2$ of the category $\mathsf{Bun}_{G,M}$, consider the set $C^{\infty}\!(P_1\odot P_2,G)^{G\times G}$ of {\em generalized gauge transformations from $P_1$ to $P_2$}.
For any triple $(P_1,P_2,P_3)$ of objects of $\mathsf{Bun}_{G,M}$, consider then the following operation 
\begin{equation*}
\begin{aligned}
\star \colon C^{\infty}\!(P_2\odot P_3,G)^{G\times G}\times C^{\infty}\!(P_1\odot P_2,G)^{G\times G}&\to C^{\infty}\!(P_1\odot P_3,G)^{G\times G}\\
\left(K_{23},K_{12}\right)&\mapsto K_{23}\star K_{12},
\end{aligned}
\end{equation*}
where $K_{23}\star K_{12}$ is defined by (\ref{eq-composisom}).
Consider then a $4$-tuple $\left(P_1,P_2,P_3,P_4\right)$ of objects of the category $\mathsf{Bun}_{G,M}$, and the three respective sets of generalized gauge transformations:
\[
C^{\infty}\!(P_1\odot P_2,G)^{G\times G},\quad C^{\infty}\!(P_2\odot P_3,G)^{G\times G}\quad \text{and}\quad C^{\infty}\!(P_3\odot P_4,G)^{G\times G}.
\]
It thus make sense to compute iterated operations of the map $\star$ as follows
\[
K_{34}\star\left(K_{23}\star K_{12}\right)\quad\text{and}\quad \left(K_{34}\star K_{23}\right)\star K_{12},
\]
for any $K_{12}\in C^{\infty}\!(P_1\odot P_2,G)^{G\times G}$, $K_{23}\in C^{\infty}\!(P_2\odot P_3,G)^{G\times G}$ and $K_{34}\in C^{\infty}\!(P_3\odot P_4,G)^{G\times G}$.
Then, explicit computations give
\begin{align*}
\left(K_{34}\star\left(K_{23}\star K_{12}\right)\right)\!(p_1,p_4)&=K_{34}\!(p_3,p_4)\left(K_{23}\star K_{12}\right)\!(p_1,p_3)=\quad \text{$(\pi_3(p_3)=\pi_1(p_1))$}\\
&=K_{34}\!(p_3,p_4)K_{23}\!(p_2,p_3) K_{12}\!(p_1,p_2)=\quad \text{$(\pi_2(p_2)=\pi_1(p_1))$}\\
&=\left(K_{34}\star K_{23}\right)\!(p_2,p_4)K_{12}\!(p_1,p_2)=\\
&=\left(\left(K_{34}\star K_{23}\right)\star K_{12}\right)\!(p_1,p_4),
\end{align*}  
which proves associativity of the operation $\star$, whenever it makes sense.

On the other hand, one may consider a pair $\left(P_1,P_2\right)$ of objects of the category $\mathsf{Bun}_{G,M}$.
As was already proved in Lemma~\ref{lem-gaugetrsf}, any gauge transformation $\sigma\in \calG_{P_1}=\calG_{P_1,P_1}$ corresponds uniquely to an element of $C^{\infty}\!(P_1\odot P_1,G)^{G\times G}$.
In particular, the unique element associated to the identity map on $P_1$ is simply $\phi_{P_1}^{-1}$.
It is interesting to compute an explicit expression for $K_{12} \star \phi_{P_1}^{-1}$, for any $K_{12}\in C^{\infty}\!(P_1\odot P_2,G)^{G\times G}$:
\begin{align*}
\left(K_{12} \star \phi_{P_1}^{-1}\right)\!(p_1,p_2)&=K_{12}\!(q_1,p_2)\phi_{P_1}\!(p_1,q_1)^{-1}=\\
&\overset{\text{by independence of the choice of $q_1$}}=K_{12}\!(p_1,p_2)\phi_{P_1}\!(p_1,p_1)^{-1}=\\
&=K_{12}\!(p_1,p_2),
\end{align*} 
where $\pi_1(p_1)=\pi_1(q_1)$.
On the other hand, using the same notations as before, it is possible to compute explicitly $\phi_{P_2}^{-1}\star K_{12}$:
\begin{align*}
\left(\phi_{P_2}^{-1}\star K_{12}\right)\!(p_1,p_2)&=\phi_{P_2}\!(q_2,p_2)^{-1}K_{12}\!(p_1,q_2)=\\
&\overset{\text{by independence of the choice of $q_2$}}=\phi_{P_2}\!(p_2,p_2)^{-1}K_{12}\!(p_1,p_2)=\\
&=K_{12}\!(p_1,p_2),
\end{align*}
where $\pi_2(p_2)=\pi_2(q_2)$.
Hence, for any object $P_1$ of the category $\mathsf{Bun}_{G,M}$, there is an element $\phi_{P_1}^{-1}$, which corresponds to the identity for the operation $\star$.

Finally, for any pair $\left(P_1,P_2\right)$ of objects of the category $\mathsf{Bun}_{G,M}$ and any bundle morphism between them represented by the generalized gauge transformation $K_{12}\in C^{\infty}\!(P_1\odot P_2,G)^{G\times G}$, one defines (see also the last part of the proof of Theorem~\ref{thm-isomequiv}) an element $\widetilde{K}_{12}$ as follows
\[
\widetilde{K}_{12}\!(p_2,p_1)\colon=K_{12}\!(p_1,p_2)^{-1},\quad (p_2,p_1)\in P_2\odot P_1.
\] 
It is almost immediate to check that $\widetilde{K}_{12}$ belongs to $C^{\infty}\!(P_2\odot P_1,G)^{G\times G}$.
Let us compute explicitly the product $\widetilde{K}_{12}\star K_{12}$:
\begin{align*}
\left(\widetilde{K}_{12}\star K_{12}\right)\!(p_1,q_1)&=\widetilde{K}_{12}\!(p_2,q_1)K_{12}\!(p_1,p_2)=\\
&=K_{12}\!(q_1,p_2)^{-1}K_{12}\!(p_1,p_2)=\\
&\overset{\text{by definition of $\phi_{P_1}$}}=K_{12}\!(p_1\phi_{P_1}\!(p_1,q_1),p_2)^{-1}K_{12}\!(p_1,p_2)=\\
&\overset{\text{by $G\times G$-equivariance of $K_{12}$}}=\phi_{P_1}\!(p_1,q_1)^{-1}K_{12}\!(p_1,p_2)^{-1}K_{12}\!(p_1,p_2)=\\
&=\phi_{P_1}\!(p_1,q_1)^{-1},
\end{align*}
where $p_2\in P_2$ is such that $\pi_2(p_2)=\pi_1(p_1)$.
On the other hand, similar computations yield
\[
K_{12}\star \widetilde{K}_{12}=\phi_{P_2},
\]
whence the assignment
\[
K_{12}\in C^{\infty}\!(P_1\odot P_2,G)^{G\times G}\leadsto C^{\infty}\!(P_2\odot P_1,G)^{G\times G}\ni \widetilde{K}_{12}
\]
gives an inverse for the operation $\star$.

To any pair of objects $\left(P_1,P_2\right)$ of $\mathsf{Bun}_{G,M}$, one can associate the set
\[
\left(P_1,P_2\right)\leadsto C^{\infty}\!(P_1\odot P_2,G)^{G\times G}\cong \calG_{P_1,P_2}.
\]
There are then maps $s$, $t$, the {\em source} and {\em target} respectively, from the set of generalized gauge transformations $C^{\infty}\!(P_1\odot P_2,G)^{G\times G}$, for any two objects $P_1$, $P_2$ of $\mathsf{Bun}_{G,M}$, to the objects of $\mathsf{Bun}_{G,M}$; the identity map $\iota$ from the objects of the category $\mathsf{Bun}_{G,M}$, to sets of generalized gauge transformations of the form $C^{\infty}\!(P\odot P,G)^{G\times G}$, for some object $P$ of $\mathsf{Bun}_{G,M}$, defined respectively via
\begin{align*}
s\left(K_{12}\right)&\colon=P_1,\quad K_{12}\in C^{\infty}\!(P_1\odot P_2,G)^{G\times G},\\ 
t\left(K_{12}\right)&\colon=P_2,\quad K_{12}\in C^{\infty}\!(P_1\odot P_2,G)^{G\times G},\\ 
i\left(P\right)&\colon=\phi_P^{-1}\in C^{\infty}\!(P\odot P,G)^{G\times G}\cong C^{\infty}\!(P,G)^{G} .
\end{align*}
There is a partially defined, associative product of the set of sets of the form $C^{\infty}\!(P_1\odot P_2,G)^{G\times G}$:
\begin{align*}
\star :C^{\infty}\!(P_2\odot P_3,G)^{G\times G}\times C^{\infty}\!(P_1\odot P_2,G)^{G\times G}&\to C^{\infty}\!(P_1\odot P_3,G)^{G\times G}\\
\left(K_{23},K_{12}\right)&\mapsto K_{23}\star K_{12}.
\end{align*}
It is obvious that 
\begin{align*}
s\left(K_{23}\star K_{12}\right)&=P_1=s\left(K_{12}\right),\\
t\left(K_{23}\star K_{12}\right)&=P_3=t\left(K_{23}\right),\\
s\left(i\left(P\right)\right)&=P=t\left(i\left(P\right)\right);\\
i\left(t\left(K_{12}\right)\right)\star K_{12}&=K_{12},\quad K_{12}\star i\left(s\left(K_{12}\right)\right)=K_{12},\quad \forall K_{12}\in C^{\infty}\!\left(P_1\odot P_2,G\right)^{G\times G}.
\end{align*}
It was also proved that there exists, for any $K_{12}\in C^{\infty}\!(P_1\odot P_2,G)^{G\times G}$, a unique element, previously denoted by $\widetilde{K}_{12}\in C^{\infty}\!(P_2\odot P_1,G)^{G\times G}$, which satisfies the property
\[
K_{12}\star \widetilde{K}_{12}=\phi_{P_2}=i\left(t\left(K_{12}\right)\right),\quad \widetilde{K}_{12}\star K_{12}=\phi_{P_1}=i\left(s\left(K_{12}\right)\right).
\]
Hence, from now on, I switch to the notation $\widetilde{K}_{12}\colon=K^{-1}_{12}$, meaning {\em not} the inverse of $K_{12}$ w.r.t.\ the multiplication in $G$, but w.r.t.\ the operation $\star$.

\section{Some general constructions for groupoids}\label{sec-groupoidconstr}
In this Section, I discuss some general facts about groupoids; I briefly discuss the concept of product groupoid, and, in more details, the concept of left and right $\calG$-spaces, for a general groupoid $\calG$.
In particular, I introduce the notion of {\em generalized conjugation for groupoids}, since a natural notion of conjugation action for groupoids is in some sense elusive.
I discuss later the concept of morphisms of groupoids, and later, of equivariant maps between left (and right) groupoid spaces, where the actions may come from distinct groupoids. 

To begin with, a general groupoid $\calG$, without any assumption on smooth structures, consists of a six-tuple $\left(\calG,X_\calG, s_\calG, t_\calG,\iota_\calG,j_\calG\right)$, whose elements are respectively two sets $\calG$ (the set of arrows) and $X_\calG$ (the set of points or the set of arrows), two surjective maps $s_\calG$ and $t_\calG$ from $\calG$ to $X$, the {\em source} and {\em target} of the groupoid, a map $\iota_\calG$ from $X$ to $\calG$, the {\em identity}, and a map $j_\calG$ from $\calG$ to itself, the {\em inversion of $\calG$}; these maps satisfy a series of axioms, for which see e.g.~\cite{M}.
Moreover, there is a well-defined partial product on the set of arrows, which is defined as follows: the set $\calG_2\subset \calG\times \calG$ consists of all pairs of elements $\left(g_1,g_2\right)$ of arrows, such that 
\[
s_\calG(g_1)=t_\calG(g_2);
\]  
then, there is an operation from $\calG_2$ to $\calG$, which is denoted as follows
\begin{align*}
\calG_2&\to\calG\\
(g_1,g_2)&\mapsto g_1g_2.
\end{align*}
The product has to be associative, whenever it makes sense:
\[
(g_1g_2)g_3=g_1(g_2g_3),\quad s_\calG(g_1)=t_\calG(g_2),\quad s_\calG(g_2)=t_\calG(g_3).
\]
The identity $\iota_\calG$ and the inversion $j_\calG$ have to satisfy additional identities involving the product, for which I refer to~\cite{M}.

\subsection{The product groupoid of two groupoids $\calG$, $\calH$}\label{ssec-prodgroupoid}
Let me introduce briefly the concept of product groupoid of two groupoids $\calG$ and $\calH$, which plays a central r{\^o}le in the definition of the generalized conjugation for a groupoid $\calG$ (see Subsection~\ref{ssec-genconjgroupoid}).
Given two groupoids $\calG$ and $\calH$, with respective sources, targets, identities and inversions, one may form the {\em product groupoid of $\calG$ and $\calH$} by putting
\begin{itemize}
\item[i)] the product set $\calG\times \calH$ as the set of arrows of the product groupoid;
\item[ii)] the product set $X_\calG\times X_\calH$ as the set of objects of the product groupoid; 
\item[iii)] the map
\[
s_{\calG\times \calH}\!\left(g,h\right)\colon=\left(s_\calG(g),s_\calH(h)\right),\quad \forall (g,h)\in \calG\times\calH,
\] 
as the source map of the product groupoid;
\item[iv)] the map
\[
t_{\calG\times \calH}\!\left(g,h\right)\colon=\left(t_\calG(g),t_\calH(h)\right),\quad \forall (g,h)\in \calG\times\calH,
\] 
as the target map of the product groupoid; 
\item[v)] the map 
\[
\iota_{\calG\times \calH}\!(x,y)\colon=\left(\iota_\calG(x),\iota_{\calH}(y)\right),\quad \forall (x,y)\in X_{\calG}\times X_\calH,
\] 
as the identity of the product groupoid;
\item[vi)] the map 
\[
j_{\calG\times \calH}\!(g,h)\colon=\left(j_\calG(g),j_\calH(h)\right),\quad \forall (g,h)\in \calG\times\calH,
\] 
as the inversion of the product groupoid;
\item[vii)] the partial product of the product groupoid is defined by the assignment
\[
(g_1,h_1)(g_2,h_2)\colon=(g_1g_2,h_1h_2),\quad s_{\calG\times\calH}\!\left(g_1,h_1\right)=t_{\calG\times\calH}\!\left(g_2,h_2\right).
\] 
\end{itemize}
Notice that the definition of product makes sense by the very definition of the source and target map in the product groupoid.
It is immediate to check that all axioms of groupoid are satisfied.

\subsection{Morphisms of groupoids}\label{ssec-morgroupoid}
Given two groupoids $\calG$ and $\calH$, I want to clarify the notion of {\em morphisms of groupoids}.
\begin{Def}\label{def-morgroupoid}
A morphism from the groupoid $\calG$ to the groupoid $\calH$ is a pair of maps $\left(\Phi,\varphi\right)$, such that the map $\Phi:\calG\to\calH$ is a map between the respective sets of arrows and the map $\varphi:X_\calG\to X_\calH$ is a map between the respective sets of objects.

Moreover, the following two conditions must hold:
\begin{itemize}
\item[i)] the three diagrams must commute   
\begin{equation}\label{eq-commdiagmorph}  
\begin{CD} 
\calG  @>\Phi>> \calH\\
@Vs_\calG VV              @VV s_\calH V\\
X_\calG @>\varphi>> X_\calH
\end{CD}\quad,\quad 
\begin{CD} 
\calG  @> \Phi >> \calH\\
@V t_\calG VV              @VV t_\calH V\\
X_\calH          @>\varphi>> X_\calH
\end{CD}\quad\text{and}\quad
\begin{CD} 
X_\calG  @> \varphi >> X_\calH\\
@V \iota_\calG VV              @VV \iota_\calH V\\
\calH          @>\Phi>> \calH
\end{CD}\quad;
\end{equation}
\item[ii)] the identity holds
\begin{equation}\label{eq-homomorph}
\Phi(g_1g_2)=\Phi(g_1)\Phi(g_2),\quad s_\calG(g_1)=t_\calG(g_2).
\end{equation} 
\end{itemize}
\end{Def}

\begin{Rem}
Notice that the commutativity of the first two diagrams in (\ref{eq-commdiagmorph}) implies that both sides of (\ref{eq-homomorph}) are well-defined.

Moreover, it is immediate to check the identity
\[
\Phi\circ j_\calG=j_\calH\circ \Phi.
\]
\end{Rem}

\subsection{The opposite groupoid}\label{ssec-oppgroupoid}
I discuss now the notion of {\em opposite groupoid} of a groupoid $\calG$, which is denoted by $\calG^{\op}$.

Set 
\begin{itemize}
\item[i)] the set of arrows $\calG$ of the initial groupoid as the set of arrows of the opposite groupoid;
\item[ii)] the set of objects $X_\calG$ of the initial groupoid as the set of objects of the opposite groupoid;  
\item[iii)] the map
\[
s_{\calG}^{\op} \colon=t_\calG
\]  
as the source map of the opposite groupoid;
\item[iv)] the map
\[
t_{\calG}^{\op}\colon=s_\calG
\]  
as the target map of the opposite groupoid; 
\item[v)] the map
\[
\iota_{\calG}^{\op}\colon=\iota_\calG
\]  
as the identity of the opposite groupoid;  
\item[vi)] the map
\[
j_{\calG}^{\op} \colon=j_\calG
\]  
as the inversion of the opposite groupoid;  
\item[vii)] the product $\cdot_{\op}$, defined via
\[
g_1\cdot_{\op} g_2\colon= g_2g_1,\quad s_{\calG}^{\op}(g_1)=t_{\calG}^{\op}(g_2), 
\] 
as the product on the opposite groupoid.
\end{itemize}
It is easy to check that the opposite groupoid, as defined above, satisfies all groupoid axioms; moreover, the pair $\left(j_\calG,\id_{X_\calG}\right)$ is an isomorphisms of groupoids from $\calG$ to the opposite groupoid $\calG^{\op}$.

\subsection{Left- and right $\calG$ actions for the groupoid $\calG$}\label{ssec-leftactgroupoid}
Given a groupoid $\calG$ and a set $M$, the next task is to define a natural notion of left $\calG$-action on the set $M$.

\begin{Def}\label{def-groupoidaction}
A left $\calG$-action of the groupoid $\calG$ on the set $M$ consists of a $3$-tuple $\left(M,J_M,\Psi_M\right)$, where $i)$ $J_M$ is a map from $M$ to the set of objects $X_\calG$ of $\calG$ (which is called the {\em momentum map} or shortly the {\em momentum} of the action) and $ii)$ a map $\Psi_M$ from $\calG\times_{J_M} M$ to $M$, where
\[
\calG\times_{J_M} M\colon=\left\{(g,m)\in\calG\times M\colon s_\calG(g)=J_M(m)\right\}.
\] 
It is customary to write
\[
\Psi_M\!\left(g,m\right)\colon= gm
\]

Moreover, the following requirements must hold
\begin{itemize}
\item[i)] \[
J_M\!\left(gm\right)=t_\calG(g),\quad \forall (g,m)\in \calG\times_{J_M}M;
\] 
\item[ii)] \[
g_1\left(g_2 m\right)=\left(g_1 g_2\right)m,\quad \forall (g_1,g_2)\in\calG_2,(g_1g_2,m)\in\calG\times_{J_M}M;
\] 
observe that Condition $\mathnormal{i})$ implies that the previous identity is well-defined.
\item[iii)] \[
\iota_\calG\left(J_M(m)\right)m=m,\quad \forall m\in M.
\] 
\end{itemize}
\end{Def}

\begin{Rem}
The definition of right $\calG$-action is similar, the only difference being that one has to switch the r{\^o}les of the source and target maps; consequently, the map $\Psi_M$ goes from the product $M\times_{J_M}\calG$ to $M$, and is denoted by
\[
\Psi_M(m,g)\colon=mg.
\]
Equivalently, a right $\calG$-action is a left $\calG^{\op}$-action, and the switch between the two actions is simply given by the inversion $j_\calG$. 
\end{Rem}

\subsection{The generalized conjugation of $\calG$}\label{ssec-genconjgroupoid}
It is well-known that there is in general no natural notion of conjugation action on a groupoid.
Nonetheless, as it is discussed in~\cite{Rgroupoid}, the conjugation action on an abstract group $G$ may be viewed as a specialization of the so-called generalized conjugation discussed at the beginning of Section~\ref{sec-isombun}; the latter action admits an easy generalization to the context of groupoids, which I am now going to discuss in detail.

First of all, consider the product groupoid $\calG^2$ of $\calG$ with itself; the aim is to construct a right action of $\calG^2$ on $\calG$.
The first ingredient one needs is a map from $\calG$ to the product of the set of objects $X_\calG$ with itself, the momentum of the (right) generalized conjugation: namely, I define
\[
J_{\conj}\!\left(g\right)\colon=\left(s_\calG(g),t_\calG(g)\right),\quad \forall g\in\calG.
\]
Consequently, the set $\calG\times_{J_{\conj}}\calG^2$ has the form
\[
\calG\times_{J_{\conj}}\calG^2=\left\{\left(g_3;g_1,g_2\right)\in\calG^3\colon \begin{cases}
t_\calG(g_1)&=s_\calG(g_3)\\
t_\calG(g_2)&=t_\calG(g_3)
\end{cases}\right\}.
\]
It makes thus sense to define the map $\Psi_{\conj}$ from $\calG\times_{J_{\conj}}\calG^2$ to $\calG$ by the rule
\begin{equation}\label{eq-genconjgroup}
\Psi_{\conj}\!\left(g_3;g_1,g_2\right)\colon= g_2^{-1}g_3 g_1,
\end{equation}
where, for the sake of simplicity, $g_2^{-1}\colon=j_\calG(g_2)$.

These data fit into a right $\calG^2$-action on $\calG$ in the sense of Definition~\ref{def-groupoidaction}, as shown in the following
\begin{Prop}\label{prop-genconj}
The triple $\left(\calG^2,J_{\conj},\Psi_{\conj}\right)$ defines a right $\calG^2$-action on $\calG$, which is called the {\em generalized conjugation of $\calG$}. 
\end{Prop}
\begin{proof}
First of all, let me compute, for any triple $\left(g_3;g_1,g_2\right)$ in $\calG\times_{J_{\conj}}\calG^2$, the following expression:
\begin{align*}
\left(J_{\conj}\circ\Psi_{\conj}\right)\!\left(g_3;g_1,g_2\right)&=J_{\conj}\!\left(g_2^{-1}g_3g_1\right)=\\
&\overset{\text{by definition of $J_{\conj}$}}=\left(s_{\calG}\!\left(g_2^{-1}g_3g_1\right),t_{\calG}\!\left(g_2^{-1}g_3g_1\right)\right)=\\
&=\left(s_{\calG}(g_1),t_{\calG}(g_2^{-1})\right)=\\
&=\left(s_{\calG}(g_1),s_{\calG}(g_2)\right)=\\
&=s_{\calG^2}\!\left(g_1,g_2\right), 
\end{align*}
which proves the first requirement for $\left(\calG^2,J_{\conj},\Psi_{\conj}\right)$ to be a right $\calG^2$ action.

Second, I compute explicitly
\begin{align*}
\Psi_{\conj}\!\left(\Psi_{\conj}\!\left(g_3;g_1,g_2\right);h_1,h_2\right)&=\Psi_{\conj}\!\left(g_2^{-1}g_3g_1;h_1,h_2\right)=\\
&=h_2^{-1}\left(g_2^{-1}g_3g_1\right)h_1=\\
&=\left(g_2h_2\right)^{-1}g_3 (g_1h_1)=\\
&=\Psi_{\conj}\!\left(g_3;g_1h_1,g_2h_2\right),
\end{align*}
whenever the identity makes sense.

Finally, one has to compute
\begin{align*}
\Psi_{\conj}\!\left(g_;\iota_{\calG^2}\!\left(J_{\conj}(g)\right)\right)&=\Psi_{\conj}\!\left(g;\iota_\calG(s_\calG(g)),\iota_\calG(t_\calG(g))\right)=\\
&=\iota_\calG(t_\calG(g))^{-1}g\iota_\calG(s_\calG(g))=\\
&=g,\quad \forall g\in \calG,
\end{align*}
whence the claim follows.
\end{proof}

\begin{Rem}\label{rem-genconj}
There is another right $\calG^2$-action on $\calG$; in fact, one can consider the map $\overline{J}_{\conj}$ from $\calG$ to $X_\calG\times X_\calG$ given by
\[
\overline{J}_{\conj}(g)\colon=\left(t_{\calG}(g),s_{\calG}(g)\right),
\]
and the map $\overline{\Psi}_{\conj}$ from $\calG\times_{\overline{J}_{\conj}}\calG^2$ to $\calG$ via
\[
\overline{\Psi}_{\conj}\!\left(g_3;g_1,g_2\right)\colon=g_1^{-1}g_3g_2.
\]
\end{Rem}
The action of Proposition~\ref{prop-genconj}, as well as the action of Remark~\ref{rem-genconj}, are called both generalized conjugations of $\calG$, because they are the natural analogues of the generalized conjugation of a group $G$, which I have briefly discussed at the beginning of Section~\ref{sec-isombun}.

Notice that there are also two left generalized conjugations of $\calG$: namely, take
\begin{equation*}
\begin{cases}
J_{\conj}^L(g)&\colon=J_{\conj}(g),\quad \Psi_{\conj}^L\!\left(g_1,g_2;g_3\right)\colon=g_2g_3g_1^{-1},\\
\overline{J}_{\conj}^L(g)&\colon=\overline{J}_{\conj}(g),\quad \overline{\Psi}_{\conj}^L\!\left(g_1,g_2;g_3\right)\colon=g_1g_3g_2^{-1},
\end{cases}
\end{equation*}
for $\left(g_1,g_2;g_3\right)$ in $\calG^2\times_{J_{\conj}^L}\calG$, resp.\ $\calG^2\times_{\overline{J}_{\conj}^L}\calG$.
The proof that both triples $\left(\calG,J_{\conj}^L,\Psi_{\conj}^L\right)$ and $\left(\calG,\overline{J}_{\conj}^L,\overline{\Psi}_{\conj}^L\right)$ define both left $\calG^2$-actions is immediate.

\subsection{Equivariant maps between groupoid-spaces}\label{ssec-equivgroupoid}
I want to define and discuss the concept of {\em equivariant map between groupoid-spaces}; here, by groupoid-space, I mean a set $M$ acted on (from the left or from the right) by a groupoid $\calG$.

For our purposes, let me consider the following general situation: given a left $\calG$-space $\left(M,J_M,\Psi_M\right)$ and a left $\calH$-space $\left(N,J_N,\Psi_N\right)$ respectively, where $\calG$, $\calH$ are any two groupoids and $M$, $N$ any two sets.
\begin{Def}\label{def-equivgroupoid}
A {\em (twisted) equivariant map} between the left $\calG$-space $M$ and the left $\calH$-space $N$ is a triple $\left(\Theta,\Phi,\varphi\right)$, where $\Theta$ is a map from the set $M$ to the set $N$, and the pair $\left(\Phi,\varphi\right)$ is a morphism from the groupoid $\calG$ to the groupoid $\calH$ in the sense of Definition~\ref{def-morgroupoid}.

Moreover, one imposes the commutativity of the following diagrams:
\begin{itemize}
\item[i)] \[
\begin{CD} 
M  @>\Theta>> N\\
@VJ_M VV              @VV J_N V\\
X_\calG @>\varphi>> X_\calH
\end{CD}\quad;
\] 
\item[ii)] \[ 
\begin{CD}
\calG\times_{J_M}M  @>\Phi\times \Theta  >> \calH\times_{J_N}N\\
@V\Psi_M VV              @VV \Psi_N V\\
M @>\Theta>> N 
\end{CD}
\]
\end{itemize}
\end{Def}  

\begin{Rem}
The first commutative diagram in Definition~\ref{def-equivgroupoid} implies that $\Phi\times \Theta$ maps $\calG\times_{J_M}M$ to $\calH\times_{J_N}N$, as the following explicit computation shows:
\begin{align*}
J_N\left(\Theta(m)\right)&=\varphi(J_M(m))=\\
&=\varphi(s_\calG(g))=\\
&=s_{\calH}(\Phi(g)),\quad \forall (g,m)\in \calG\times_{J_M} M.
\end{align*} 
Usually, the second diagram is encoded in the following identity:
\[
\Theta(gm)=\Phi(g)\Theta(m),\quad \forall (g,m)\in \calG\times_{J_M} M,
\]
which corresponds clearly to the usual definition of (twisted by $\Phi$) equivariance of a map $\Theta$ from a left $G$-space to a left $H$-space, for $G$, $H$ usual groups.
\end{Rem}

\subsection{Some explicit computations for the groupoid of generalized gauge transformations}\label{ssec-explcompgengauge}
Finally, I want to compute an explicit expression for the generalized conjugation of the groupoid of generalized gauge transformations.

For a general element $K_{12}\in C^{\infty}\!\left(P_1\odot P_2,G\right)^{G\times G}$, the generalized conjugation of it makes sense, whenever one considers the action of pairs $\left(K_{i1},K_{j2}\right)$, where $K_{ik}\in C^{\infty}\!\left(P_i\odot P_k,G\right)^{G\times G}$, for some principal $G$-bundle $P_i$ over $M$, and $k=1,2$.
The explicit formula takes the form
\begin{align*}
\left(K_{j2}^{-1}\star K_{12}\star K_{i1}\right)\!(p_i,p_j)&=K_{j2}\!(p_j,p_2)^{-1}K_{12}\!(p_1,p_2)K_{i1}\!(p_i,p_1), 
\end{align*}                      
where $p_1\in P_1$ and $p_2\in P_2$ are such that $\pi_{1}(p_1)=\pi_2(p_2)=\pi_i(p_i)=\pi_j(p_j)$.
The previous equation simplifies considerably, when considering the generalized conjugation of $K_{12}$ by a special pair $(g,h)$, where $g\in C^{\infty}\!(P_1,G)^G$ and $h\in C^{\infty}\!(P_2,G)^G$, i.e.\ when considering the action of a pair of gauge transformations on the source and on the target of $K_{12}$.
A direct computation shows immediately
\begin{align*}
\left(K_h^{-1}\star K_{12}\star K_g\right)\!(p_1,p_2)&=K_h\!(q_2,p_2)^{-1}K_{12}\!(q_1,q_2)K_g\!(p_1,q_1)=\\
&=h(p_2)^{-1}\phi_{P_2}\!(p_2,q_2)K_{12}\!(q_1,q_2)\phi_{P_1}\!(p_1,q_1)^{-1}g(p_1)=\\
&=h(p_2)^{-1}K_{12}\!(q_1 \phi_{P_1}\!(p_1,q_1)^{-1},q_2 \phi_{P_2}\!(p_2,q_2)^{-1})g(p_1)=\\
&\overset{\text{by (\ref{eq-canbundle})}}=h(p_2)^{-1} K_{12}\!(p_1,p_2) g(p_1),
\end{align*} 
whenever $\pi_1(p_1)=\pi_1(q_1)$ and $\pi_2(p_2)=\pi_2(q_2)$.

\section{Applications to bundles on the space of loops or paths}\label{sec-applpullbun}
In the present section, I interpret the holonomy as a gauge transformation of a bundle over the space of loops in the manifold $M$, and similarly the parallel transport as a generalized gauge transformation on bundles over the space of loops or, more generally, paths in $M$.

Consider first the space of free loops in $M$, which is denoted by $\lo M$; by the word ``free'', one means loops in $M$ without a specified base point.
There are canonical, smooth maps from the product manifold $\lo M\times \unint$ to $M$, the {\em evaluation map} and the {\em evaluation map at $0$}, defined respectively via
\begin{align}
\label{eq-evmap} \ev\left(\gamma;t\right)\colon&=\gamma(t);\\
\label{eq-evinit} \ev_0\left(\gamma;t\right)\colon&=\gamma(0).
\end{align}
Using the maps (\ref{eq-evmap}) and (\ref{eq-evinit}), one can construct two principal $G$-bundles on the cylinder $\lo M\times\unint$ as follows: consider a fixed principal $G$-bundle $P$ over $M$, and take the pull-backs
\begin{align*}
\ev^* P&=\left\{\left(\gamma;t;p\right)\in \lo M\times \unint\times P\colon \pi(p)=\gamma(t)\right\};\\
\ev_0^* P&=\left\{\left(\gamma;t;p\right)\in \lo M\times \unint\times P\colon \pi(p)=\gamma(0)\right\}.
\end{align*}
The canonical projections from both bundles onto $\lo M\times \unint$ map any $3$-tuple to the first two arguments; this allows us to write the fibred product of $\ev^*P$ and $\ev_0^*P$ as
\begin{equation*}
\ev_0^* P\odot \ev^* P=\left\{\left(\gamma;t;p,\widetilde{p}\right)\in \lo M\times \unint\times P\times P\colon \pi(p)=\gamma(0),\quad \pi\left(\tildep\right)=\gamma(t)\right\}.
\end{equation*}

Choosing a background connection $A$ on $P$, one can define (see Section~\ref{sec-holpartransp}) the holonomy w.r.t.\ $A$ and the parallel transport w.r.t.\ $A$.
The holonomy is an element of $G$, depending on a loop $\gamma$ and on a chosen base point $p\in P$, satisfying $\pi(p)=\gamma(0)=\ev_0(\gamma,t)$. 
Hence, it may be viewed as a map on the pulled-back bundle $\ev_0^*P$ with values in $G$, which depends additionally on the choice of a connection $A$ on $P$.
Moreover, Equation (\ref{eq-liehol}) of Lemma~\ref{lem-propholonom} implies that the holonomy is a $G$-equivariant map from $\ev_0^*P$ to $G$; thus, it may be viewed as a gauge transformation of $\ev_0^*P$.
Similarly, the parallel transport w.r.t.\ a chosen connection $A$ depends on a curve $\gamma$ in $M$, some $t\in\unint$, and two points in $P$, $(p,q)$, satisfying $\pi(p)=\gamma(0)$ and $\pi(q)=\gamma(t)$; hence, the parallel transport may be viewed as a function on the fibred product of $\ev_0^*P$ and $\ev^*P$ with values in $G$.
Moreover, Equation (\ref{eq-liepartr}) of Lemma~\ref{lem-proppartr} guarantees that the parallel transport is $G\times G$-equivariant w.r.t.\ the generalized conjugation introduced at the beginning of Section~\ref{sec-isombun}.
Thus, parallel transport defines a generalized gauge transformation between the bundle $\ev_0^*P$ and $\ev^*P$, depending on the choice of a connection $A$.

From now on, the holonomy $\holo{A}{\bullet}{\bullet}$, viewed as a gauge transformation of $\ev_0^*P$, is denoted simply by $\holo{A}{0}{1}$, whereas the parallel transport $\holonom{A}{\bullet}{\bullet}{\bullet}{\bullet}$ as a generalized gauge transformation is denoted by $\holo{A}{0}{\bullet}$.
Theorem~\ref{thm-isomequiv} of Section~\ref{sec-isombun} implies that the parallel transport $\holo{A}{0}{\bullet}$ induces a bundle (iso)morphism from $\ev_0^*P$ to $\ev^* P$, which is denoted by $\Phi_A$ in order to make explicit its dependence on the chosen connection $A$.

I now want to inspect now more carefully the dependence on the chosen connection $A$ of the isomorphism $\Phi_A$, in view of the computations done in Subsection~\ref{ssec-explcompgengauge} of Section~\ref{sec-groupoidconstr}.
First, write down the isomorphism $\Phi_A$ explicitly:
\begin{equation*}
\Phi_A\!\left(\gamma;t;p\right)=\left(\gamma;t;q \holonom{A}{\gamma}{t}{p}{q}\right),
\end{equation*}
where $\pi(p)=\gamma(0)$, and $q\in P$ is any point satisfying $\pi\left(q\right)=\gamma(t)$. 

It is well-known that the gauge group $\calG$ operates on the space of connections $\calA$ by pull-back. 
Considering now a gauge transformation $\sigma$ and taking the pull-back of the connection $A$ w.r.t.\ $\sigma$, the morphism $\Phi_A$ changes as follows:
\begin{equation}\label{eq-connisom}
\begin{aligned}
\Phi_{\sigma^* A}\left(\gamma;t;p\right)&=(\gamma;t;\tildep\ \holonom{\sigma^*A}{\gamma}{t}{p}{\tildep})=\\
&=(\gamma;t;\tildep\ g_{\sigma}\!\left(\tildep\right)^{-1}\! \holonom{A}{\gamma}{t}{p}{\tildep}\! g_{\sigma}(p))=\\
&=\left(\gamma;t;\sigma^{-1}\!\left(\tildep\right)\ \holonom{A}{\gamma}{t}{p\ g_{\sigma}(p)}{\tildep}\right)=\\
&=\left(\gamma;t;\sigma^{-1}\!\left(\tildep\ \holonom{A}{\gamma}{t}{\sigma(p)}{\tildep}\right)\right)=\\
&=\sigma^{-1}\left(\Phi_A \left(\gamma;t;\sigma(p)\right)\right).
\end{aligned}
\end{equation}
The main argument used in the proof of the above equation is Identity (\ref{eq-gaugepartr}) of Lemma~\ref{lem-proppartr}; this identity admits an interpretation in terms of equivariant maps between groupoid-spaces, as introduced in Subsection~\ref{ssec-equivgroupoid} of Section~\ref{sec-groupoidconstr}.

I need some technical arguments to make the statement more precise.
Assume $f$ to be a smooth map from some manifold $L$ to $M$; then
\begin{Lem}\label{lem-pullgauge}
Every gauge transformation of $P$ induces a gauge transformation of the pulled-back bundle $f^* P$.
\end{Lem}
\begin{proof}
Consider a general gauge transformation $\sigma\in \calG_P$.
By definition,
\[
f^* P=\left\{\left(l,p\right)\in L\times P\colon f(l)=\pi(p)\right\};
\]
consider then a general point $(l,p)$ in $f^* P$. 
Define a morphism $\sigma_f$ from $f^* P$ to itself as follows:
\begin{equation*}
\sigma_f\!\left(l,p\right)\colon=\left(l,\sigma(p)\right).
\end{equation*}
One has then to show that $\sigma_f$ is well-defined, that it respects the fibres of $f^*P$ and that it is $G$-equivariant; bijectivity, and hence invertibility, follows from Lemma~\ref{lem-gaugeisom} of Section~\ref{sec-isombun}.

Since $\pi\circ \sigma=\pi$, it follows that the image of a point in $f^* P$ still belongs to $f^* P$. 
By its very definition, $\sigma_f$ is fibre-preserving.
The $G$-equivariance of $\sigma_f$ is immediate.

Notice that the map taking $\sigma$ to $\sigma_f$ is clearly a homomorphism of groups from $\calG_P$ to $\calG_{f^* P}$.
\end{proof}
The gauge group of $P$, $\calG_P$, can be viewed as a groupoid, where the set of objects $X_{\calG_P}$ is simply a point, namely, one can consider this point as $P$, an object of the category $\mathsf{Bun}_{G,M}$; source map and target map are simply
\[
s_{\calG_P}(\sigma)=t_{\calG_P}(\sigma)=P;
\]     
identity is simply 
\[
\iota_P(P)=\id_P.
\]
Finally, the inversion is clearly
\[
j_{\calG_P}(\sigma)=\sigma^{-1}.
\]
Moreover, as already seen at the end of Section~\ref{sec-backgr}, the gauge group $\calG_P$ operates on the space of connections on $P$, $\calA_P$; it is easy to see that $\calA_P$ is a right $\calG_P$-space in the sense explained in Subsection~\ref{ssec-leftactgroupoid} of Section~\ref{sec-groupoidconstr}.

I want to construct a morphism of groupoids between $\calG_P$ and the product groupoid of the groupoid of generalized gauge transformations of the space of loops (or paths) in $M$ with itself. 
The groupoid of generalized gauge transformations, for a manifold $M$ and group $G$, is denoted by $\calG_{G,M}$ (hence, it follows $\calG_{G,M,P_1,P_2}=C^{\infty}!\left(P_1\odot P_2,G\right)^{G\times G}$, $P_1$ and $P_2$ being two principal $G$-bundles over $M$).
\begin{Lem}\label{lem-princmorphism}
The pair $\left(\Phi_P,\varphi_P\right)$, where
\[
\begin{cases}
\varphi_P(P)&\colon=\left(\ev_0^*P,\ev^*P\right)\\
\Phi_P(\sigma)&\colon=\left(\sigma_{\ev_0},\sigma_{\ev}\right)
\end{cases},
\]
is a morphism in the sense of Definition~\ref{def-morgroupoid} from the groupoid $\calG_{G,\lo M}^2$, for any manifold $M$ and any Lie group $G$..
\end{Lem}
\begin{proof}
According to Definition~\ref{def-morgroupoid}, one has to show first the commutativity of the three diagrams of (\ref{eq-commdiagmorph}); this is immediate, because of the definition of source map, target map and identity for the product groupoid of the groupoid of generalized gauge transformations of $\lo M$ (see Subsection~\ref{ssec-prodgroupoid} for more details).

Property (\ref{eq-homomorph}) follows immediately from Lemma~\ref{lem-pullgauge} and from the definition of product of a product groupoid.
\end{proof}

Using now all arguments sketched in Section~\ref{sec-isombun} and~\ref{sec-groupoidconstr}, one can restate Identity (\ref{eq-gaugepartr}) of Lemma~\ref{lem-proppartr} as follows:
\begin{Prop}\label{prop-equivpartr}
The triple $\left(\holo{\bullet}{0}{\bullet},\Phi_P,\varphi_P\right)$ defines an equivariant map from the space of connections $\calA_P$, viewed as a right $\calG_P$-space, to the groupoid $\calG_{G,\lo M}$, resp.\ $\calG_{G,\PM}$ of generalized gauge transformations of the space of loops (resp.\ paths) in $M$, viewed as a right $\calG_{G,\lo M}^2$-, resp.\ $\calG_{G,\PM}^2$-,space via the generalized conjugation introduced in Subsection~\ref{ssec-genconjgroupoid}.
(By $\holo{\bullet}{0}{\bullet}$, I mean the map associating to a connection $A$ the corresponding parallel transport as a generalized gauge transformation.)
\end{Prop}
\begin{proof}
First of all, the momentum $J_{\calA_P}$ is simply given by
\[
J_{\calA_P}\!(A)\colon=P,\quad \forall A\in\calA_P.
\]
With this definition in mind, and recalling the arguments of Subsection~\ref{ssec-genconjgroupoid}, one gets immediately the commutativity of the first diagram in Definition~\ref{def-equivgroupoid}. 

It remains to show the commutativity of the second diagram.
This is an immediate consequence of the explicit computation done in Subsection~\ref{ssec-explcompgengauge} together with Identity (\ref{eq-gaugepartr}) of Lemma~\ref{lem-proppartr}.
\end{proof}

An immediate consequence of Proposition~\ref{prop-equivpartr}, together with Identity (\ref{eq-gaugehol}) of Lemma~\ref{lem-propholonom}, is given by
\begin{Cor}\label{cor-equivhol}
The holonomy map $\holo{\bullet}{0}{1}$, assigning to any connection $A$ in $\calA_P$ the gauge transformation of the bundle $\ev_0^*P$ defined by the holonomy w.r.t.\ $A$, is an equivariant map from the space of connections on $P$ to the gauge group of $\ev_0^*P$ w.r.t.\ the conjugation.
\end{Cor}

\subsubsection{A consequence of Corollary~\ref{cor-equivhol}: the Wilson loop map}\label{sssec-invfunct}
In this subsubsection, I discuss an important consequence of Corollary~\ref{cor-equivhol}, namely the construction of the so-called {\em Wilson loop map}; this map plays an important r{\^o}le in the framework of Topological Quantum Field Theories, in particular Chern--Simons Theory.

To begin with, notice that the holonomy map $\holo{A}{0}{1}$, for any given connection $A$, defines immediately a section of the associated bundle $\Ad\!(\ev_0^*P)$ on $\lo M\times \unint$ w.r.t.\ conjugation.
For our purposes, it is better to view the holonomy map $\holo{A}{0}{1}$ as a map from the pulled-back bundle $\ev(0)^*P$ over $\lo M$, where $\ev(0)$ is defined via
\[
\ev(0)(\gamma)\colon=\gamma(0)
\]  
(see also the beginning of Subsection~\ref{ssec-restrboundary}).                       
Similarly, this holonomy map, denoted again by $\holo{A}{0}{1}$, descends to a section of the associated bundle $\Ad\!(\ev(0)^*P)$; this bundle is associated to $\ev(0)^*P$ via the conjugation of $G$. 
More generally, dropping out the explicit dependence on connections, the holonomy map $\holo{\bullet}{0}{1}$ may be viewed as a map from the cartesian product of $\calA_P$, the space of connections on $P$, with the principal bundle $\ev(0)^*P$ on the space of loops $\lo M$, to the group $G$; this map satisfies, by Corollary~\ref{cor-equivhol}, two kinds of equivariance, namely w.r.t.\ the $G$-action on $\ev(0)^*P$ and w.r.t.\ the $\calG_P$-action on $\calA_P$.
Consider now a representation $(V,\rho)$ of $G$, namely a (real or complex) vector space $V$ endowed with a group homomorphism $\rho$ from $G$ to the automorphism group of $V$.
Composing the holonomy map with $\rho$, one gets a map from the product $\calA_P\times \ev(0)^*P$ to the linear group $GL(V)$; the aforementioned equivariance properties hold again, since $\rho$ is a group homomorphism.
Finally, one can take the trace in the endomorphism ring of $V$ of the composite map $\rho\circ\holo{\bullet}{0}{1}$; the result of these operations is denoted by
\[
W_\rho(A;\gamma;p)\colon=\tr_V\!\left(\rho(\holon{A}{\gamma}{p})\right),
\] 
which is called {\em Wilson loop map $W_\rho$ w.r.t.\ the representation $\rho$}.
It is clearly a (real or complex, depending on the representation $V$) function on $\calA_P\times \ev(0)^*P$.

The cyclicity of the trace $\tr_V$, together with Corollary~\ref{cor-equivhol}, implies the two following invariance properties of the Wilson loop map: 
\begin{align}
\label{eq-wilslieinv}W_\rho(A;\gamma;pg)&=W_\rho(A;\gamma;p),\quad \forall g\in G,\\
\label{eq-wilsgaugeinv}W_\rho(A^\sigma;\gamma;p)&=W_\rho(A;\gamma;p),\quad \forall \sigma\in\calG.
\end{align}
Let me spend a few words on the space of connections on $P$, $\calA_P$: as already seen at the end of Section~\ref{sec-backgr}, the gauge group $\calG_P$ operates on the right on $\calA_P$.
It would naturally be interesting to analyse the properties of the gauge group action of $\calG_P$ on $\calA_P$ and to find suitable conditions for this action to be free; in fact, there could be in principle gauge transformations fixing connections, e.g.\ if the structure group $G$ possesses a nontrivial centre, there are constant gauge transformations fixing any given connection.
Therefore, in order to consider the {\em moduli space of connections on $P$}, i.e.\ the quotient space $\calA_P/ \calG_P$, one has to consider only the space $\calA_P^*$ of {\em irreducible connections}, which is the subset of $\calA_P$ on which the group $\calG_P^*$ (the quotient group of $\calG_P$ by the centre of the structure group $G$, which may be viewed as a normal subgroup of $\calG_P$).
Hence, it makes sense to consider the quotient space $\calA_P^*/ \calG_P^*$, which one can view as the base space of a principal $\calG_P^*$-bundle, whose total space is $\calA_P^*$; it is not my plan to discuss the local structure of this quotient space, thus I do not enter into the details of the construction of suitable topologies on it.

Equations (\ref{eq-wilslieinv}) and (\ref{eq-wilsgaugeinv}) imply immediately the following 
\begin{Prop}\label{prop-wilsloop}
For any linear representation $(V,\rho)$ of $G$, the structure group of $P$, the Wilson loop map
\[
\calA_P^*\times \ev(0)^*P\ni (A;\gamma;p)\mapsto W_\rho(A;\gamma;p)
\]
descends to a function $W_\rho$ on the product of the moduli space of connections $\calA_P^*/ \calG_P^*$ with the space of loops $\lo M$.
\end{Prop}

\subsection{Restriction to the boundary of $\unint$}\label{ssec-restrboundary}
I want to discuss how the restriction to the boundary of the ``cylinder'' $\lo M\times \unint$ affects the isomorphism between $\ev_0^*P$ and $\ev^*P$ given by the parallel transport $\holo{A}{0}{\bullet}$; in fact, I will prove later that both restrictions give well-known gauge transformations of a bundle over the space of loops $\lo M$.

First of all, one can define the evaluation map at some point $t\in\unint$, which is denoted by $\ev(t)$, as a map from the space of paths $\PM$ to $M$ by
\[
\ev(t)(\gamma)\colon=\gamma(t).
\]
It is clear that, taking restriction to $\lo M$, the evaluation maps at $0$ and at $1$ do coincide.
Restricting the bundles $\ev_0^* P$ and $\ev^* P$ to the subset $\lo M\times \{0\}$ (which is clearly diffeomorphic to $\lo M$), one gets immediately
\begin{equation*}
\begin{aligned}
\ev_0^*P\vert_{\lo M\times \{0\}}&\cong\left\{(\gamma;t;p)\in \ev_0^* P\colon t=0 \right\}\cong \ev(0)^* P;\\
\ev^* P\vert_{\lo M\times \{0\}}&\cong\left\{(\gamma;t;p)\in \ev^* P\colon t=0 \right\}\cong\\
&\cong\left\{(\gamma;0;p)\in \ev^* P\colon \gamma(0)=\pi(p) \right\}\cong\\
&\cong \ev(0)^* P.
\end{aligned}
\end{equation*}

The isomorphism $\Phi_A$ respects fibres, whence it follows that it maps the restriction to $\lo M\times \{0\}$ of $\ev_0^* P$ to the restriction to the same set of $\ev^* P$, hence mapping $\ev(0)^* P$ to itself.
By the very definition of $\Phi_A$, it follows:
\[
\Phi_A\vert_{\lo M\times \{0\}}\left(\gamma;0;p\right)=(\gamma;0;q \holonom{A}{\gamma}{0}{p}{q}),
\]
where in this case $\pi\left(q\right)=\gamma(0)=\pi(p)$. 
One can choose $q=p$, since the definition of $\Phi_A$ above does not depend on the choice of $q$ as long as it remains in the same fibre, and thus it follows
\[
\widetilde{\gamma}_{A,p}(0)=p,
\]
whence $\holonom{A}{\gamma}{0}{p}{p}=e$. 
Hence, the restriction to $\lo M\times \{0\}$ equals the identity gauge transformation of $\ev(0)^* P$.
On the other hand, taking restriction to the bundles $\ev_0^* P$ and $\ev^* P$ to $\lo M\times \{1\}$ (which is again diffeomorphic to $\lo M$); in this case, one gets
\begin{equation*}
\begin{aligned}
\ev_0^* P\vert_{\lo M\times \{1\}}&=\left\{\left(\gamma;1;p\right)\in \ev_0^* P\colon \pi(p)=\gamma(0)\right\}\cong \ev(0)^* P;\\
\ev^* P\vert_{\lo M\times \{1\}}&=\left\{\left(\gamma;1;p\right)\in \ev_0^* P\colon \pi(p)=\gamma(1)=\gamma(0)\right\}\cong \ev(0)^* P.
\end{aligned}
\end{equation*}
So, restricting the two bundles $\ev_0^* P$ and $\ev^* P$ to the subset $\lo M\times \{1\}$, one gets the same bundle as before, namely $\ev(0)^* P$. 
But the restriction of $\Phi_A$ is {\em not} the identity gauge transformation, as the following computation shows:
\[
\Phi_A\vert_{\lo M\times \{1\}}\left(\gamma;1;p\right)=\left(\gamma;1;q \holonom{A}{\gamma}{1}{p}{q}\right),
\]
where $q \in P$ obeys $\pi\left(q \right)=\gamma(1)=\gamma(0)=\pi(p)$. 
Hence, one can again choose $q=p$, and by definition
\[
\Phi_A\vert_{\lo M\times \{1\}}\left(\gamma;1;p\right)=\left(\gamma;1;p \holonom{A}{\gamma}{1}{p}{p}\right).
\]
The following identity holds:
\[
p \holonom{A}{\gamma}{1}{p}{p}=\widetilde{\gamma}_{A,p}(1)=p \holon{A}{\gamma}{p}.
\]
Since the action of $G$ on each fibre is free, it follows:
\[
\holon{A}{\gamma}{p}=\holonom{A}{\gamma}{1}{p}{p},\quad \forall \gamma\in \lo M,p\in P.
\]
Hence, the restriction of $\Phi_A$ to $\lo M\times \{1\}$ equals the gauge transformation of $\ev(0)^* P$ defined by the holonomy $\holon{A}{\bullet}{\bullet}$ (which is denoted by $\holo{A}{0}1$).

\section{Some consequences of flatness related to holonomy and parallel transport}\label{sec-flatness}
Motivated by Proposition~\ref{prop-equivpartr} of Section~\ref{sec-applpullbun}, which allows us to parameterize particular generalized gauge transformations of the space of loops or paths in $M$ by connections on a given background bundle $P$ over $M$ in an equivariant way, I am interested now in the restriction of the equivariant map of Proposition~\ref{prop-equivpartr} to the space of {\em flat connections}, on which the gauge group of $P$ also operates.
Before entering into the details, I need some technical facts about connections on the fibred product (see Section~\ref{sec-isombun} for more details) of any two principal bundles.

\subsection{Connections on the fibred product of two principal bundles}\label{ssec-connfiber}
Assume $P$, $\widetilde{P}$ to be objects (a priori distinct) of the category $\mathsf{Bun}_{G,M}$, endowed with respective connections $A$ and $\widetilde{A}$.
A general tangent vector to the fibred product $P\odot\widetilde{P}$ at a point $\left(p,\widetilde{p}\right)$ can be written as a pair $\left(X_p,X_{\widetilde{p}}\right)$, where $X_p$, resp.\ $X_{\widetilde{p}}$, is a tangent vector to $P$ at $p$, resp.\ to $\widetilde{P}$ at $\widetilde{p}$.
By the very definition of the projection $\overline{\pi}$ from $P\odot\widetilde{P}$ onto $M$, the tangent vector $\left(X_p,X_{\widetilde{p}}\right)$ is vertical if and only if
\[
\tangent{\left(p,\widetilde{p}\right)}{\overline{\pi}}\left(X_p,X_{\widetilde{p}}\right)=\tangent{p}{\pi}\left(X_p\right)=\tangent{\widetilde{p}}{\widetilde{\pi}}\left(X_{\widetilde{p}}\right)=0,
\]
i.e.\, if and only if both components of $\left(X_p,X_{\widetilde{p}}\right)$ are vertical.

As seen in Section~\ref{sec-backgr}, the connection $A$, resp.\ $\widetilde{A}$, specifies a smooth splitting of the tangent bundle of $P$, resp.\ $\widetilde{P}$, into vertical and $A$-horizontal bundle, resp.\ $\widetilde{A}$-horizontal bundle.
Hence, one can write any tangent vector $X_p$, resp.\ $X_{\widetilde{p}}$, to $P$ at $p$, resp.\ to $\widetilde{P}$ at $\widetilde{p}$, as
\[
X_p=X_p^v+X_p^h,\quad\text{resp.}\quad X_{\widetilde{p}}=X_{\widetilde{p}}^v+X_{\widetilde{p}}^h.
\]
Thus, any tangent vector $\left(X_p,X_{\widetilde{p}}\right)$ to $P\odot\widetilde{P}$ at $\left(p,\widetilde{p}\right)$ has a unique splitting
\begin{equation}\label{eq-horsplit}
\left(X_p,X_{\widetilde{p}}\right)=\left(X_p^v,X_{\widetilde{p}}^v\right)+\left(X_p^h,X_{\widetilde{p}}^h\right).
\end{equation}
It is clear from the definition of the action of $G\times G$ on $P\odot\widetilde{P}$ that the sets of all vectors of the form $\left(X_p^h,X_{\widetilde{p}}^h\right)$ in (\ref{eq-horsplit}), for any pair $(p,\tildep)\in P\odot\widetilde{P}$, specify a $G\times G$-invariant family of subspaces.
Therefore, by the arguments recalled in Section~\ref{sec-backgr}, to the splitting (\ref{eq-horsplit}) belongs obviously the connection $1$-form
\[
\left(A\odot\widetilde{A}\right)_{\!\left(p,\widetilde{p}\right)}\!\left(X_p,X_{\widetilde{p}}\right)\colon=\left(A_p\!\left(X_p\right),\widetilde{A}_{\widetilde{p}}\!\left(X_{\widetilde{p}}\right)\right),
\]
which is called the {\em fibred product connection} of the connections $A$ and $\widetilde{A}$; the fibred product connection $A\odot\widetilde{A}$ is a $\Lg\oplus\Lg$-valued form on $P\odot\widetilde{P}$.
Notice that the fibred product connection $A\odot\widetilde{A}$ can be written alternatively as the restriction of the connection on the product $P\times \widetilde{P}$ given by
\[
\left(\pr_{P}^*\!A,\pr_{\widetilde{P}}^*\!\widetilde{A}\right),
\]
where $\pr_P$, resp.\ $\pr_{\widetilde{P}}$, denotes the projection of $P\odot\widetilde{P}$ onto $P$, resp.\ $\widetilde{P}$, to $P\odot \widetilde{P}$.

Recall that a connection $A$ on $P$ is {\em flat}, if it satisfies the equation
\begin{equation}\label{eq-flatconn}
\dd A+\frac{1}2\Lie{A}{A}=0.
\end{equation}
An equivalent characterisation of flatness of a connection $A$, viewed as a way of splitting the tangent bundle of $P$ into vertical and $A$-horizontal bundle, is that the $A$-horizontal bundle specifies an {\em integrable distribution} over $P$, i.e.\ the Lie bracket of any two $A$-horizontal vectors is again horizontal.
From the latter characterisation of the fibred product connection $A\odot \widetilde{A}$, it follows that the fibred product connection is flat if and only if both $A$ and $\widetilde{A}$ are flat; it is a consequence of Equation (\ref{eq-flatconn}).

Given now a bundle morphism $\sigma$ from $P$ to $\widetilde{P}$ in the sense of Definition~\ref{def-bunisom}, one can construct a map $\sigma\odot\sigma^{-1}$ from $P\odot\widetilde{P}$ to $\widetilde{P}\odot P$ as follows:
\[
\left(\sigma\odot\sigma^{-1}\right)\left(p,\widetilde{p}\right)\colon=\left(\sigma(p),\sigma^{-1}\left(\widetilde{p}\right)\right),\quad \forall\ \left(p,\widetilde{p}\right)\in P\odot\widetilde{P}.
\]

\begin{Lem}\label{lem-fibprodmap}
The {\em fibred product map $\sigma\odot\sigma^{-1}$}, for any bundle morphism $\sigma$ from $P$ to $\widetilde{P}$, is a bundle morphism from the fibred product $P\odot\widetilde{P}$ to $\widetilde{P}\odot P$ (notice that both are $G\times G$-bundles).
\end{Lem}
\begin{proof}
First of all, since $\widetilde{\pi}\circ \sigma=\pi$, it follows immediately that $\sigma\odot\sigma^{-1}$ maps the fibred product $P\odot\widetilde{P}$ to $\widetilde{P}\odot P$.

It is clear that $\sigma\odot\sigma^{-1}$ is fibre-preserving; namely, the canonical projection from $\widetilde{P}\odot P$ to $M$ is e.g.\ the projection $\widetilde{\pi}$ from the first argument of any pair to $M$, whence it follows
\[
\left(\widetilde{\pi}\circ \sigma\right)(p)=\pi(p),
\]
which equals the canonical projection from $P\odot \widetilde{P}$.

Finally, a direct computation following from $G$-equivariance of $\sigma$, gives
\begin{align*}
\left(\sigma\odot\sigma^{-1}\right)\!\left((p,\tildep)(g,h)\right)&=\left(\sigma(pg),\sigma^{-1}(\tildep h)\right)=\\
&=\left(\sigma(p)g,\sigma^{-1}(\tildep)h\right)=\\
&=\left(\sigma(p),\sigma^{-1}(\tildep)\right)(g,h)=\\
&=\left(\left(\sigma\odot\sigma^{-1}\right)\!(p,\tildep)\right)(g,h),
\end{align*}
proving $G\times G$-equivariance.
\end{proof}
Recall that, by Lemma~\ref{lem-gaugeisom}, any bundle morphism $\sigma$ as above is invertible.  
Since $\sigma$ is an isomorphism, it is clear that $\sigma\odot\sigma^{-1}$ is also an isomorphism.
Moreover, the $G$-equivariance of $\sigma$ ensures that $\sigma\odot\sigma^{-1}$ is $G\times G$-equivariant.

There is a fibred product connection $\widetilde{A}\odot A$ on $\widetilde{P}\odot P$.
On the other hand, by pulling back the fibred product connection $\widetilde{A}\odot A$ w.r.t.\ $\sigma\odot\sigma^{-1}$, one gets a connection $1$-form on $P\odot\widetilde{P}$, by the $G\times G$-equivariance of $\sigma\odot\sigma^{-1}$.

\begin{Lem}\label{lem-fibprodconn}
The pull-back of a fibred product connection $\widetilde{A}\odot A$ on $\widetilde{P}\odot P$ w.r.t.\ the fibred product map $\sigma\odot\sigma^{-1}$ equals a fibred product connection $B\odot \widetilde{B}$ on $P\odot \widetilde{P}$, if and only if
\[
\begin{cases}
\sigma^*\widetilde{A}&=B\\
\sigma^*\widetilde{B}&=A.
\end{cases}
\]
\end{Lem}
\begin{proof}
An immediate computation gives the following identity
\begin{equation}\label{eq-tangfibprod}
\tangent{(p,\tildep)}{\sigma\odot\sigma^{-1}}\!(X_p,X_{\tildep})=\left(\tangent{p}{\sigma}(X_p),\tangent{\tildep}{\sigma^{-1}}(X_{\tildep})\right),
\end{equation}
for any tangent vector $(X_p,X_{\tildep})$ to $P\odot\widetilde{P}$ at $(p,\tildep)$.

It follows immediately from the definition of fibred product connection that
\begin{align*}
\left(\left(\sigma\odot\sigma^{-1}\right)^*\!\left(\widetilde{A}\odot A\right)\right)_{(p,\tildep)}\!\left(X_p,X_{\tildep}\right)&=\left(B\odot\widetilde{B}\right)_{(p,\tildep)}\!\left(X_p,X_{\tildep}\right)\Longleftrightarrow\\
\left(\widetilde{A}_{\sigma(p)}\!\left(\tangent{p}{\sigma}(X_p)\right),A_{\sigma^{-1}(\tildep)}\!\left(\tangent{\tildep}{\sigma^{-1}}(X_{\tildep})\right)\right)&=\left(B_p(X_p),\widetilde{B}_{\tildep}(X_{\tildep}\right).
\end{align*}
The surjectivity of $\pi$ and $\widetilde{\pi}$ yields then the claim.
\end{proof}
Hence, by Lemma~\ref{lem-fibprodconn}, if the pull-back of $\widetilde{A}\odot A$ w.r.t.\ $\sigma\odot \sigma^{-1}$ equals $A\odot\widetilde{A}$, the connection $\widetilde{A}$ is the pull-back of $A$ w.r.t.\ $\sigma$.

\subsection{Composition properties and inversion of the parallel transport}\label{ssec-compinvpartr}
In this subsection, I state and prove two technical lemmata, which will be used later, when proving the two theorems on the consequences of flatness of the reference connection $A$ on both holonomy and parallel transport, viewed respectively as a gauge transformation of the pulled-back bundle $\ev_0^*P$ on $\lo M\times\unint$ (or also $\ev(0)^*P$ on $\lo M$) and as a generalized gauge transformation on $\lo M$ or, more generally $\PM$.

In fact, I will prove later how the parallel transport behaves w.r.t.\ the composition and inversion of paths.

\begin{Lem}\label{lem-comphol}
Let $\gamma_1$ and $\gamma_2$ be two composable paths in $M$, in the sense that $\gamma_1(1)=\gamma_2(0)$.
Additionally, let us consider three points $p_1$, $p_2$ and $p_3$ in $P$, such that 
\[
\gamma_1(0)=\pi(p_1),\quad \gamma_1(1)=\gamma_2(0)=\pi(p_2),\quad \gamma_2(1)=\pi(p_3).
\]
Define the {\em composite curve of $\gamma_1$ and $\gamma_2$}, denoted by$\gamma_2\circ\gamma_1$, as
\[
\left(\gamma_2\circ \gamma_1\right)\!\left(t\right)\colon=\begin{cases}
\gamma_1(2t),& t\in \left[0;\frac{1}2\right]\\
\gamma_2\left(2t-1\right),& t\in\left[\frac{1}2;1\right]
\end{cases}.
\]
Clearly, the composite curve $\gamma_2\circ \gamma_1$ is piecewise smooth.

Then, the following identity holds
\begin{equation}\label{eq-comphol}
\holonom{A}{\gamma_2\circ\gamma_1}{1}{p_1}{p_3}=\holonom{A}{\gamma_2}{1}{p_2}{p_3}\ \holonom{A}{\gamma_1}{1}{p_1}{p_2}.
\end{equation}
\end{Lem}
\begin{proof}
By Definition~\ref{def-partransp}, one gets
\[
\widetilde{\gamma_2\circ\gamma_1}_{A,p_1}(1)=p_3\ \holonom{A}{\gamma_2\circ\gamma_1}{1}{p_1}{p_3},
\]
where $\widetilde{\gamma_2\circ\gamma_1}_{A,p_1}$ is the unique horizontal lift of $\gamma_2\circ \gamma_1$ based at $p_1$.

On the other hand, consider the composite curve $\widetilde{\gamma_2}_{A,\widetilde{\gamma_1}_{A,p_1}(1)}\circ\widetilde{\gamma_1}_{A,p_1}$ of horizontal curves:
\[
\left(\widetilde{\gamma_2}_{A,\widetilde{\gamma_1}_{A,p_1}(1)}\circ\widetilde{\gamma_1}_{A,p_1}\right)\left(t\right)\colon=\begin{cases}
\widetilde{\gamma_1}_{A,p_1}(2t),& t\in\left[0,\frac{1}2\right]\\
\widetilde{\gamma_2}_{A,\widetilde{\gamma_1}_{A,p_1}(1)}(2t-1),& t\in\left[\frac{1}2;1\right].
\end{cases}
\]
By its very definition, it is clear that $\widetilde{\gamma_2}_{A,\widetilde{\gamma_1}_{A,p_1}(1)}\circ\widetilde{\gamma_1}_{A,p_1}$ lies over $\gamma_2\circ\gamma_1$ and that it is based at $p_1$.
Since it is the composition of two horizontal curves, it is also horizontal, whence it follows by Theorem~\ref{thm-horlift} of Section~\ref{sec-holpartransp}
\[
\widetilde{\gamma_2}_{A,\widetilde{\gamma_1}_{A,p_1}(1)}\circ\widetilde{\gamma_1}_{A,p_1}=\widetilde{\gamma_2\circ\gamma_1}_{A,p_1}.
\]
Therefore, one gets
\begin{align*}
\left(\widetilde{\gamma_2}_{A,\widetilde{\gamma_1}_{A,p_1}(1)}\circ\widetilde{\gamma_1}_{A,p_1}\right)\!(1)&=\widetilde{\gamma_2}_{A,\widetilde{\gamma_1}_{A,p_1}(1)}(1)=\\
&=\widetilde{\gamma_2}_{A,p_2}(1)\ \holonom{A}{\gamma_1}{1}{p_1}{p_2}=\\
&=p_3\ \holonom{A}{\gamma_2}{1}{p_2}{p_3}\ \holonom{A}{\gamma_1}{1}{p}{p_1},
\end{align*}
whence the claim follows.
\end{proof}

\begin{Lem}\label{lem-invhol}
Let $\gamma$ be a path in $M$; additionally, let $p$, $q$ be two point in $P$, such that
\[
\gamma(0)=\pi(p),\quad \gamma(1)=\pi(q).
\]
Define the {\em inverse curve $\gamma^{-1}$ of $\gamma$} by
\[
\gamma^{-1}\!\left(t\right)\colon=\gamma\left(1-t\right).
\]

Then, the following identity holds:
\begin{equation}\label{eq-invhol}
\holonom{A}{\gamma^{-1}}{1}{q}{p}=\holonom{A}{\gamma}{1}{p}{q}^{-1}.
\end{equation}
\end{Lem}
\begin{proof}
Consider the following curve: 
\[
\widetilde{\gamma}_{A,p\holonom{A}{\gamma}{1}{p}{q}^{-1}}^{-1}(t)\colon=\widetilde{\gamma}_{A,p}(1-t)\ \holonom{A}{\gamma}{1}{p}{q}^{-1}.
\]
It is immediate to see that this curve lies over $\gamma^{-1}$.
It is clearly based at $q$, since
\begin{align*}
\widetilde{\gamma}_{A,p}(1)\ \holonom{A}{\gamma}{1}{p}{q}^{-1}&=q \holonom{A}{\gamma}{1}{p}{q}\ \holonom{A}{\gamma}{1}{p}{q}^{-1}=\\
&=\widetilde{p}.
\end{align*}
Since the connection $A$ defines a $G$-invariant distribution, $\widetilde{\gamma}_{A,p \holonom{A}{\gamma}{1}{p}{q}^{-1}}^{-1}$ is $A$-horizontal, and by the uniqueness part of Theorem~\ref{thm-horlift} of Section~\ref{sec-holpartransp}, one gets
\[
\widetilde{\gamma}_{A,p \holonom{A}{\gamma}{1}{p}{q}^{-1}}^{-1}=\widetilde{\gamma^{-1}}_{A,q},
\]
whence the claim follows.
\end{proof}

\subsection{Holonomy and flatness}\label{ssec-holflatness}
Let $A$ be a flat connection on $P$.
One of the most important features of flat connections is encoded in the following
\begin{Thm}\label{thm-holflat}
If the connection $A$ is flat, the holonomy $\holon{A}{\gamma}{p}$ of a loop $\gamma$, based at $p\in P$ over $\gamma(0)$, depends {\em only} on the homotopy class of the loop $\gamma$.
\end{Thm}
(See \cite{KN1}, Chapter 2, Sections 4, 7 and 8 for more details; let me just recall that it is a consequence of the Ambrose--Singer Theorem, which characterizes completely the Lie algebra of the holonomy group $\Phi(p)$ w.r.t.\ a given connection on a principal bundle in terms of its curvature form, $p$ being a reference point in $P$.)
On the other hand (see e.g.~\cite{Kob}, Chapter 1), it is well-known that a representation of the fundamental group $\pi_1(M)$ of $M$ in $G$ implies the existence of a principal bundle $P$ with a flat connection, whose holonomy representation of the fundamental group of $M$ in $G$ coincides with the representation giving rise to the bundle with flat connection; this is the well-known correspondence between (gauge-equivalence classes of) flat connections and (conjugacy classes in $G$ of) representations of the fundamental group of $M$.
Notice that $\pi_1(M)$ is defined by considering {\em continuous} loops, whereas the holonomy w.r.t.\ a connection is defined in terms of {\em piecewise smooth} loops; the connection between the two settings is encoded in the so-called {\em restricted holonomy group}, for whose precise definition see~\cite{KN1}, Chapter 2, Section 4. 

As a consequence of Corollary~\ref{cor-equivhol} of Section~\ref{sec-applpullbun}, one can view the holonomy $\holon{A}{0}{1}$ w.r.t.\ $A$ as a map from the pulled-back bundle $\ev(0)^*P$ to $G$, equivariant w.r.t.\ conjugation in $G$; the corresponding gauge-transformation of $\ev(0)^*P$ is denoted by $\Phi_A$.
There is a natural map from $\ev(0)^*P$ to $P$, denoted by $\widetilde{\ev(0)}$:
\[
\widetilde{\ev(0)}\left(\gamma;p\right)\colon=p.
\]
Clearly, $\widetilde{\ev(0)}$ is $G$-equivariant; thus, the pull-back of $A$ w.r.t.\ $\widetilde{\ev(0)}$ is a connection on $\ev(0)^*P$, which is denoted by $\ev(0)^*A$. 

Now, I have all elements needed to state the following
\begin{Thm}\label{thm-flatness}
The gauge transformation $\Phi_A$ of $\ev(0)^*P$ induced by the holonomy w.r.t.\ $A$ stabilizes the connection $\ev(0)^*A$, i.e.\ $\Phi_A^* \ev(0)^*A=\ev(0)^*A$, if and only if $A$ is flat.
\end{Thm}

Before proving Theorem~\ref{thm-flatness}, I need a few preliminary facts.
A general tangent vector at a point $\left(\gamma;p\right)\in \ev(0)^*P$ can be written as a pair of tangent vectors $\left(X_\gamma;X_p\right)$, where $X_\gamma\in \tangent{\gamma}{\lo M}$ (i.e.\ it is tangent to $\lo M$ at $\gamma$, which is equivalent to saying that $X_\gamma$ is a section of the tangent bundle of $M$ w.r.t.\ $\gamma$) and $X_p\in \tangent{p}P$ (i.e.\ it is tangent to $P$ at $p$).
The condition for the tangent vector $\left(X_\gamma;X_p\right)$ to $\ev(0)^*P$ at $\left(\gamma;p\right)$ to be vertical is that
\[
X_\gamma=0\Rightarrow \tangent{p}{\pi}(X_p)=\tangent{\gamma}{\ev(0)}(X_\gamma)=0,
\]
since $\left(X_\gamma;X_p\right)$ is tangent to $\ev(0)^*P$.
Therefore, any vertical vector to $\ev(0)^*P$ at $\left(\gamma;p\right)$ can be uniquely written as $\left(0;X_p\right)$, where $X_p$ is vertical at $p$.

On the other hand, the tangent vector $\left(X_\gamma;X_p\right)$ at $\left(\gamma;p\right)\in \ev(0)^*P$ is $\ev(0)^*A$-horizontal if and only if
\[
\left(\ev(0)^*A\right)_{\left(\gamma;p\right)}\left(X_\gamma;X_p\right)\overset{!}=A_p\!\left(X_p\right)=0.
\]
Hence, a tangent vector $\left(X_\gamma;X_p\right)$ at $\left(\gamma;p\right)$ is $\ev(0)^*A$-horizontal if and only if its $P$-component is $A$-horizontal.

\fbox{\parbox{12cm}{\bf Since the connection $A$ specifies a splitting of the tangent bundle of $P$ into the vertical bundle and the $A$-horizontal bundle, it is possible to split any tangent vector $\left(X_\gamma;X_p\right)$ to $\ev(0)^*P$ at $\left(\gamma;p\right)$ into a unique sum
\begin{equation}\label{eq-verthorsplit}
\left(X_\gamma;X_p\right)=\left(0;X_p^v\right)+\left(X_\gamma;X_p^h\right),
\end{equation}
where $X_p^v$, resp.\ $X_p^h$, denotes the vertical, resp.\ $A$-horizontal, part of the vector $X_p$.}}
The splitting (\ref{eq-verthorsplit}) plays a pivotal r{\^o}le in the proof of Theorem~\ref{thm-flatness}, to which I come now.
\begin{proof}[Proof of Theorem~\ref{thm-flatness}]
Let me first prove that the flatness of $A$ implies that $\ev(0)^*A$ is fixed by the gauge transformation $\Phi_A$; Lemma~\ref{lem-flathoriz} implies the converse.
Consider a general tangent vector $\left(X_\gamma;X_p\right)$ to $\ev(0)^*P$ at a general point $\left(\gamma;p\right)$, and assume $X_p$ to be $A$-horizontal at $p$.
One can assume that $\left(X_\gamma;X_p\right)$ is the tangent vector at $\left(\gamma;p\right)$ of a path $\left(\gamma_s;p_s\right)$ ($s\in \unint$), where $\gamma_s$ is a path in $\lo M$ and $p_s$ is an $A$-horizontal curve in $P$.
Moreover, since $\left(\gamma_s;p_s\right)$ belongs to $\ev(0)^*P$, the curve $p_s$ lies over $\gamma_s(0)$, for any $s\in \unint$.

Define a curve $\overline{\gamma}_s$ in $M$, for any $s\in\unint$, by the rule
\begin{equation}\label{eq-curvinit}
\overline{\gamma}_s\!\left(t\right)\colon=\gamma_{st}(0).
\end{equation}

Construct now a map $\Gamma$ on the unit square $\unint\times\unint$ by means of the family of curves (\ref{eq-curvinit}) as follows:
\begin{equation*}
\Gamma\left(t,s\right)\colon=\begin{cases}
\overline{\gamma}_s\!\left(3t\right),& t\in \left[0;\frac{1}3\right]\\
\gamma_s\!\left(3t-1\right),& t\in \left[\frac{1}3;\frac{2}3\right]\\
\overline{\gamma}_s\!\left(3-3t\right),& t\in\left[\frac{2}3,1\right].
\end{cases}
\end{equation*}
First of all, for any $s\in\unint$, $\Gamma(\bullet,s)$ is a closed curve: in fact,
\[
\Gamma(0,s)=\overline{\gamma}_s(0)=\gamma_{0}(0)=\gamma(0),\quad \Gamma(1,s)=\overline{\gamma}_s(0)=\gamma(0).
\]
It follows also that $\Gamma(\bullet,s)$ is based in $\gamma(0)$, for any $s$ in the unit interval.
Finally, it is clear that $\Gamma$ is a homotopy of $\gamma$.

Since any $\Gamma(\bullet,s):=\Gamma_s$ is homotopic to $\gamma$ and the connection $A$ is by assumption flat, one gets, by Theorem~\ref{thm-holflat}, 
\[
\holon{A}{\Gamma_s}{p}=\holon{A}{\gamma}{p},\quad \forall s\in \unint.
\]
On the other hand, for any $s\in\unint$, the curve $\Gamma_s$ may be written as the composite curve $\overline{\gamma}_s^{-1}\circ\gamma_s\circ\overline{\gamma}_s$.
Then, Lemma~\ref{lem-comphol} and~\ref{lem-invhol} imply together
\begin{align*}
\holon{A}{\gamma}{p}&=\holon{A}{\Gamma_s}{p}=\\
&=\holonom{A}{\overline{\gamma}_s}{1}{p}{p_s}^{-1}\ \holon{A}{\gamma_s}{p_s}\ \holonom{A}{\overline{\gamma}_s}{1}{p}{p_s}.
\end{align*}

Define, for any $s\in\unint$, the curve $\overline{p}_s(t)$ by the rule
\[
\overline{p}_s\!\left(t\right)\colon=p_{ts},
\]
where $t$ is also in the unit interval.

The curve $\overline{p}_s$ is clearly $A$-horizontal, as it is a reparametrization of an $A$-horizontal curve.
Moreover, it lies over $\overline{\gamma}_s$:
\[
\pi\left(\overline{p}_s\!(t)\right)=\pi\left(p_{st}\right)=\gamma_{st}(0)=\overline{\gamma}_s(t),\quad \forall t\in \unint,
\]
and is based at $p$, since $\overline{p}_s(0)=p_0=p$.

Hence, the curve $\overline{p}_s$, for any $s\in\unint$, is the unique $A$-horizontal lift of the curve $\overline{\gamma}_s$, thus it follows
\[ 
\holonom{A}{\overline{\gamma}_s}{1}{p}{p_s}=e,\quad \forall s\in\unint.
\] 
As a consequence, one gets the following identity:
\[
\holon{A}{\gamma}{p}=\holon{A}{\gamma_s}{p_s},\quad \forall s\in\unint,
\]
which is equivalent to the fact that the holonomy of the family of curves $\left(\gamma_s;p_s\right)$ in $\ev(0)^*P$ is constant, and this yields in turn
\[
\tangent{\left(\gamma;p\right)}{\holo{A}{0}{1}}\left(X_\gamma;X_p\right)=0,\quad \text{if $X_p$ is $A$-horizontal}.
\]

An explicit computation gives
\begin{equation}\label{eq-derholon}
\begin{aligned}
\tangent{\left(\gamma;p\right)}\Phi_A\left(X_\gamma;X_p\right)&=\left(X_\gamma;\tangent{p}{\raction_{\holon{A}{\gamma}{p}}}\left(X_p\right)\right.+\\
&\phantom{=}\left.+\tangent{\holon{A}{\gamma}{p}}{\laction_p}\left[\tangent{\left(\gamma;p\right)}{\holo{A}{0}{1}}\left(X_\gamma;X_p\right)\right]\right),
\end{aligned}
\end{equation}
for any tangent vector $\left(X_\gamma;X_p\right)$ on $\ev(0)^*P$ at $\left(\gamma;p\right)$.

Therefore, if $X_p$ is $A$-horizontal, the above equation simplifies to
\[
\tangent{\left(\gamma;p\right)}\Phi_A\left(X_\gamma;X_p\right)=\left(X_\gamma;\tangent{p}{\raction_{\holon{A}{\gamma}{p}}}\left(X_p\right)\right).
\]
Finally, one gets by direct computations
\begin{align*}
\Phi_A^*\left(\ev(0)^*A\right)_{\left(\gamma;p\right)}\left(X_\gamma;X_p\right)&=\left(\ev(0)^*A\right)_{\left(\gamma;p\! \holon{A}{\gamma}{p}\right)}\left[\left(X_\gamma;\tangent{p}{\raction_{\holon{A}{\gamma}{p}}}\left(X_p\right)\right)\right]=\\
&=A_{p\! \holon{A}{\gamma}{p}}\left[\tangent{p}{\raction_{\holon{A}{\gamma}{p}}}\left(X_p\right)\right]=\\
&=\Ad\left(\holon{A}{\gamma}{p}^{-1}\right)\left[A_p\left(X_p\right)\right]=\\
&=0,
\end{align*}
if $\left(X_\gamma;X_p\right)$ is $\ev(0)^*A$-horizontal; in the third equality, I made use of the equivariance of $A$.

\fbox{\parbox{12cm}{\bf Hence, any $\ev(0)^*A$-horizontal vector at $\left(\gamma;p\right)$ is also $\Phi_A^*\left(\ev(0)^*A\right)$-horizontal.}}

It follows immediately that $\Phi_A^*\left(\ev(0)^*A\right)=\ev(0)^*A$, since at any point $\left(\gamma;p\right)\in \ev(0)^*P$ a tangent vector $\left(X_\gamma;X_p\right)$ can be decomposed in a unique way into a vertical piece and an $\ev(0)^*A$-horizontal piece; by the above computations, $\Phi_A^*\left(\ev(0)^*A\right)$ and $\ev(0)^*A$ agree on any tangent vector, and thus they must be identical.
\end{proof}

As mentioned before, to prove the converse claim one needs the following 
\begin{Lem}\label{lem-flathoriz}
The condition that the gauge transformation $\Phi_A$ associated to the holonomy w.r.t.\ $A$ fixes the connection $\ev(0)^*A$ can be alternatively encoded into the following horizontality equation:
\begin{equation*}
\dd_{\ev(0)^*A}\holo{A}{0}{1}=0.
\end{equation*}  
\end{Lem}
\begin{proof}
Theorem~\ref{thm-holflat} implies that, for any couple $(\gamma;p)\in\ev(0)^*P$ and any tangent vector $(X_\gamma;X_p)$ to $\ev(0)^*P$ at $(\gamma;p)$, the identity holds:
\[
\left(\Phi_A^*\!(\ev(0)^*A)\right)_{(\gamma;p)}\!(X_\gamma;X_p)=\ev(0)^*A_{(\gamma;p)}\!(X_\gamma;X_p).
\]

Equation (\ref{eq-derholon}) gives an explicit expression for the tangent map of the holonomy; one can split the right-hand side of Equation (\ref{eq-derholon}) into two pieces:
\begin{align*}
\tangent{\left(\gamma;p\right)}\Phi_A\left(X_\gamma;X_p\right)&=\left(X_\gamma;\tangent{p}{\raction_{\holon{A}{\gamma}{p}}}\left(X_p\right)\right)+\\
&\phantom{=}+\left(0;\tangent{\holon{A}{\gamma}{p}}{\laction_p}\left[\tangent{\left(\gamma;p\right)}{\holo{A}{0}{1}}\left(X_\gamma;X_p\right)\right]\right).
\end{align*}
The second term of the pair on the right-hand side may be also rewritten as
\begin{align*}
\tangent{\holon{A}{\gamma}{p}}{\laction_p}\left[\tangent{\left(\gamma;p\right)}{\holo{A}{0}{1}}\left(X_\gamma;X_p\right)\right]=&\tangent{e}{\laction_{p\holon{A}{\gamma}{p}}}\!\left\{\tangent{\holon{A}{\gamma}{p}}{\laction_{\holon{A}{\gamma}{p}^{-1}}}\right.\\
&\left.\left[\tangent{\left(\gamma;p\right)}{\holo{A}{0}{1}}\left(X_\gamma;X_p\right)\right]\right\},
\end{align*}
which turns out to be a vertical vector in $\ev(0)^*P$.

Thus, it follows by the very definition of pull-back:
\begin{align*}
\left(\Phi_A^*\!(\ev(0)^*A)\right)_{(\gamma;p)}\!(X_\gamma;X_p)&=\ev(0)^*A_{(\gamma;p\holon{A}{\gamma}{p})}\!\left(X_\gamma;\tangent{p}{\raction_{\holon{A}{\gamma}{p}}}\left(X_p\right)\right)+\\
&\phantom{=}+\ev(0)^*A_{(\gamma;p\holon{A}{\gamma}{p})}\!\left(0;\tangent{\holon{A}{\gamma}{p}}{\laction_p}\right.\\
&\phantom{=+}\left.\left[\tangent{\left(\gamma;p\right)}{\holo{A}{0}{1}}\left(X_\gamma;X_p\right)\right]\right)=\\
&=\Ad\left(\holon{A}{\gamma}{p}^{-1}\right)\ev(0)^*A_{(\gamma;p)}\!(X_\gamma;X_p)+\\
&\phantom{=}+\tangent{\holon{A}{\gamma}{p}}{\laction_{\holon{A}{\gamma}{p}^{-1}}}\left[\tangent{\left(\gamma;p\right)}{\holo{A}{0}{1}}\left(X_\gamma;X_p\right)\right],
\end{align*}
where in the third equality, I made use of the $G$-equivariance of $\ev(0)^*A$ to get the first term and of verticality of $\ev(0)^*A$ to get the second term.

By writing down explicitly the adjoint action, the equation at the beginning of the proof, stating the fact that $\Phi_A$ fixes the connection $\ev(0)^*A$, takes the form
\begin{equation}\label{eq-horizcond}
\tangent{\left(\gamma;p\right)}{\holo{A}{0}{1}}\left(X_\gamma;X_p\right)+\tangent{e}{\raction_{\holon{A}{\gamma}{p}}}\!\left[A_p(X_p)\right]-\tangent{e}{\laction_{\holon{A}{\gamma}{p}}}\!\left[A_p(X_p)\right]=0,
\end{equation}
and recall that, by the very definition of $\ev(0)^*A$, one has
\[
A_p(X_p)=\ev(0)^*A_{(\gamma;p)}\!\left(X_\gamma;X_p\right).
\]
This is exactly the horizontality equation for the holonomy.
\end{proof}
Lemma~\ref{lem-flathoriz} implies that the gauge transformation $\Phi_A$, viewed as a $G$-equivariant map from $\ev(0)^*P$ to $G$, is constant along $\ev(0)^*A$-horizontal curves in $\ev(0)^*P$: namely, consider a $\ev(0)^*A$-horizontal curve $\widetilde{\gamma}(s)=\left(\Gamma_s,p_s\right)$ in $\ev(0)^*P$, where $\Gamma_s$ is a curve in $\lo M$, i.e.\ a curve of loops, whereas $p_s$ denotes a $A$-horizontal curve in $P$, such that the projection of $p_s$ onto $M$ equals $\Gamma_s(0)$, the initial (and final) point of the loop $\Gamma_s$.
Since $\widetilde{\gamma}$ is $\ev(0)^*A$-horizontal, it follows that
\[
\ev(0)^*A_{\widetilde{\gamma}}\!\left(\ddss \widetilde{\gamma}\right)=0,
\]
whence it follows, by the very definition of the tangent map and by (\ref{eq-horizcond}),
\[
\tangent{\widetilde{\gamma}}{\holo{A}{0}{1}}\!\left(\ddss \widetilde{\gamma}\right)=\ddss \holo{A}{0}{1}(\widetilde{\gamma})\overset{!}=0.
\]
This is equivalent to the condition that $\holo{A}{0}1$ is constant along the path $\widetilde{\gamma}$; this is the main argument in the proof of the ``only if'' part of Theorem~\ref{thm-flatness}.
\begin{proof}[End of the proof of Theorem~\ref{thm-flatness}]
By what was recalled shortly after the statement of Theorem~\ref{thm-holflat}, given a connection $A$, it suffices to prove that the holonomy representation w.r.t.\ $A$ based at some point $p\in P$ over some $x\in M$,
\[
\lo M\ni \gamma\mapsto \holo{A}{0}{1}\!\left(\gamma,p\right)\in G,
\]
descends to a representation of the fundamental group $\pi_1(M)$.
To be more precise (see~\cite{KN1}, Chapter 2, Section 4, for more details), one has to consider the holonomy as a representation of the fundamental group $\pi_1(M)$ (consisting of all {\em continuous loops in $M$}) into the quotient Lie group $\Phi(p)/ \Phi^0(p)$, where, borrowing the notations from~\cite{KN1}, $\Phi(p)$, resp.\ $\Phi^0(p)$, denotes the {\em holonomy group w.r.t.\ $A$ with reference point $p$}, resp.\ {\em the restricted holonomy group w.r.t\ $A$ with reference point $p$}; they are respectively subgroups of $G$ consisting of all elements which can be written in the form $\holo{A}{0}1(\gamma,p)$, for $\gamma$ a (piecewise smooth) loop in $M$, resp.\ a (piecewise smooth) loop in $M$ homotopic to the constant loop.
Hence, in this framework, to prove that the holonomy gives truly a representation of $\pi_1(M)$ to $G$, one has to show that the restricted holonomy group $\Phi^0(p)$ with reference point $p$ is trivial.   
To see this, consider a (piecewise smooth) loop $\gamma$ based at $x$, homotopic to the constant loop $e_x$ based at $x$; a base point $p$ over $x$ is chosen.
A homotopy between $\gamma$ and $e_x$ is denoted by $\Gamma$, or, in other words,
\[
\Gamma(t,0)=\gamma(t),\quad \Gamma(t,1)=e_x(t)=x,\quad \Gamma(0,s)=x,\quad \forall (s,t)\in \unint\times\unint, 
\]
whence it follows that
\[
\pi(p)=\Gamma(0,s),\quad\forall s\in \unint.
\]
Using the so called {\em Factorization Lemma} (see~\cite{KN1}, Appendix 7, for the proof), one can assume that the homotopy $\Gamma$ is also piecewise smooth.
Thus, one can consider the curve
\[
\unint\ni s\mapsto \widetilde{\gamma}(s)\colon=\left(\Gamma(\bullet,s),p\right);
\]
it is obvious that $\widetilde{\gamma}$ lies in $\ev(0)^*P$ and, by its very construction, $\widetilde{\gamma}$ is $\ev(0)^*A$-horizontal.
Therefore, the gauge transformation $\Phi_A$ is constant along the curve $\widetilde{\gamma}$, i.e.
\[
\holo{A}{0}1\!\left(\Gamma(\bullet,0),p\right)=\holo{A}{0}{1}\!\left(\gamma,p\right)\overset{!}=\holo{A}{0}1\!\left(\Gamma(\bullet,1),p\right)=\holo{A}{0}{1}\!\left(e_x,p\right)=e,
\]
whence the claim follows.

Thus, the quotient group $\Phi(p)/ \Phi^0(p)$ equals $\Phi(p)$, which is a subgroup of $G$; it follows that $A$ is flat, since the corresponding holonomy gives a representation of the fundamental group of $M$ in $G$.
\end{proof}

\subsubsection{Flatness and the Wilson loop map}\label{sssec-flatwils}
In this subsubsection, I want to discuss the consequences of Theorem~\ref{thm-holflat} and Lemma~\ref{lem-flathoriz} for the Wilson loop map, introduced and discussed in Subsubsection~\ref{sssec-invfunct}.

Recall that the Wilson loop map $W_\rho$, defined upon the choice of a linear representation $(V,\rho)$ of the structure group $G$, is a function on the product of the moduli space of (irreducible) connections $\calA_P^*/ \calG_P^*$ and the space of loops $\lo M$.
Moreover, the Wilson loop function $W_\rho$ comes from a $\calG_P^*\times G$-invariant function on the product bundle $\calA_P^*\times\ev(0)^*P$, which is simply the trace of the holonomy map $\holo{\bullet}{0}{1}$ in the representation $\rho$.

Notice that, by the very definition of flatness, the (restricted) gauge group $\calG_P^*$ operates on the subspace of $\calA_P^*$ consisting of flat connections; hence, one can consider the Wilson loop map $W_\rho$ on the {\em moduli space of flat connections $\calM_P$ on $P$}, which is simply
\[
\calM_P\colon=\left\{A\in\calA_P^*\colon F_A=0\right\}/ \calG_P^*.
\]
Lemma~\ref{lem-flathoriz} yields the following
\begin{Prop}\label{prop-derhamhol}
The Wilson loop map $W_\rho$ on the product $\calM_P\times \lo M$ is closed w.r.t.\ the De Rham differential, for any class of flat connections; equivalently, the map
\[
\calM_P\ni[A]\mapsto W_\rho(A),
\]
assigning to every class of flat connections the Wilson loop map w.r.t.\ $[A]$ on $\lo M$, is map from the moduli space of (irreducible) flat connections to the $0$-th De Rham cohomology group of $\lo M$.
\end{Prop}
\begin{proof}
It suffices to prove that, for any (irreducible) flat connection $A$ on $P$, the function
\[
\ev(0)^*P\ni(\gamma;p)\mapsto W_\rho(A)(\gamma;p)\colon=W_\rho(A;\gamma;p)
\]
on the bundle $\ev(0)^*P$ is $\dd$-closed.

The horizontality equation of Lemma~\ref{lem-flathoriz} takes the following form in the linear group $\GL(V)$:
\[
\tangent{(\gamma;p)}{\rho(\holo{A}{0}{1})}\!(X_\gamma;X_p)+\overline{\rho}(A_p(X_p))\rho(\holo{A}{0}{1})(\gamma;p)-\rho(\holo{A}{0}{1})(\gamma;p)\overline{\rho}(A_p(X_p))=0,
\]
where $\overline{\rho}\colon=\tangent{e}{\rho}$ is the derived Lie algebra-representation of $\Lg$ coming from the $G$-representation $\rho$.

By its very definition, the following equation for the exterior derivative of the Wilson loop map $W_\rho(A)$ holds:
\begin{align*}
\dd W_\rho(A)_{(\gamma;p)}\!(X_\gamma;X_p)&=\dd \tr_V\left(\rho(\holo{A}{0}{1}\right)_{(\gamma;p)}\!(X_\gamma;X_p)=\\
&=\tr_V\left(\tangent{(\gamma;p)}{\rho(\holo{A}{0}{1})}\!(X_\gamma;X_p)\right)=\\
&=\tr_V\left(\tangent{(\gamma;p)}{\rho(\holo{A}{0}{1})}\!(X_\gamma;X_p)\right)+\\
&\phantom{=}+\tr_V\left(\overline{\rho}(A_p(X_p))\rho(\holo{A}{0}{1})(\gamma;p)\right)-\\
&\phantom{=}-\tr_V\left(\rho(\holo{A}{0}{1})(\gamma;p)\overline{\rho}(A_p(X_p))\right)=\\
&=0,
\end{align*} 
where the third identity is a consequence of the cyclicity of the trace, and the fourth identity follows directly from the horizontality equation of Lemma~\ref{lem-flathoriz}.
\end{proof}
Hence, to any orbit in the moduli space $\calM_P$, the Wilson loop map $W_\rho(A)$ associates an element of the $0$-th De Rham cohomology group, i.e. a locally constant function on the space of loops.

\subsection{Flatness and parallel transport}\label{ssec-holpartr}
Let me return to the principal bundles $\ev^*P$ and $\ev_0^*P$ on the cylinder $\lo M\times\unint$ or $\PM\times\unint$.
Denote by $\widetilde{\ev}$, resp.\ $\widetilde{\ev_0}$, the natural map from $\ev^*P$, resp.\ $\ev_0^*P$, to $P$, given by
\[
\widetilde{\ev}\left(\gamma;t;p\right)\colon=p,\quad \text{resp.}\quad \widetilde{\ev_0}\left(\gamma;t;q\right)\colon=q.
\]
Denote also by $\ev^*A$, resp.\ $\ev_0^*A$, the pull-back of $A$ w.r.t.\ $\widetilde{\ev}$, resp.\ $\widetilde{\ev_0}$; $\ev^*A$, resp.\ $\ev_0^*A$, is clearly a connection on $\ev^*P$, resp.\ $\ev_0^*P$, since both maps $\widetilde{\ev}$ and $\widetilde{\ev_0}$ are $G$-equivariant.
Consider the fibred product of $\ev_0^*P$ and $\ev^*P$, and a point $\left(\gamma;t;p,q\right)$ in it.
Denote by the $4$-tuple $\left(X_\gamma;X_t;X_p,X_{q}\right)$ a general tangent vector to $\ev_0^*P\odot\ev^*P$ at a point $\left(\gamma;t;p,q\right)$, where $i)$ $X_\gamma$ is tangent to $\lo M$ at the curve $\gamma$, $ii)$ $X_t$ is tangent to $\unint$ at $t$, and $iii)$ $X_p$, resp.\ $X_{q}$, is tangent to $P$ at $p$, resp.\ $q$. 

The condition on the tangent vector $\left(X_\gamma;X_t;X_p,X_{q}\right)$ to be vertical may be translated into the set of equations
\[
X_\gamma=0,\quad X_t=0\Rightarrow\ \begin{cases}
\tangent{p}{\pi}(X_p)&=\tangent{\left(\gamma,t\right)}{\ev_0}\left(X_\gamma;X_t\right)=0,\\
\tangent{q}{\pi}\left(X_{q}\right)&=\tangent{\left(\gamma,t\right)}{\ev}\left(X_\gamma;X_t\right)=0,
\end{cases}
\]
where the first identities follow from the definition of the bundle projection from $\ev_0^*P\odot\ev^*P$ onto $\lo M\times\unint$, and the remaining are consequences of $\left(X_\gamma;X_t;X_p,X_{q}\right)$ being tangent to $\ev_0^*P\odot\ev^*P$.Hence, a tangent vector $\left(X_\gamma;X_t;X_p,X_{q}\right)$ to the fibred product $\ev_0^*P\odot\ev^*P$ at a point $\left(\gamma;t;p,q\right)$ is vertical if and only if its $P$-pieces are vertical and its $\lo M\times\unint$-piece vanishes.
As already known, the connection $A$ specifies a splitting of the tangent bundle of $P$.
Therefore, any tangent vector $X_p$ to $P$ at $p$ can be written uniquely as $X_p=X_p^v+X_p^h$, where $X_p^v$, resp.\ $X_p^h$, denotes the vertical part, resp.\ the horizontal part, of $X_p$.

\fbox{\parbox{10.5cm}{\bf Thus, it is possible to write in a unique way
\begin{equation*}
\left(X_\gamma;X_t;X_p,X_{q}\right)=\left(0;0;X_p^v,X_{q}^v\right)+\left(X_\gamma;X_t;X_p^h,X_{q}^h\right),
\end{equation*}
for any tangent vector $\left(X_\gamma;X_t;X_p,X_{q}\right)$ to $\ev_0^*P\odot\ev^*P$ at $\left(\gamma;t;p,q\right)$.}}
The above splitting specifies in an obvious way a $G\times G$-invariant distribution on $\ev_0^*P\odot\ev^*P$; the corresponding connection $1$-form is simply the fibred product connection
\begin{equation*}
\left(\ev_0^*A\odot\ev^*A\right)_{\left(\gamma;t;p,q\right)}\left(X_\gamma;X_t;X_p,X_{q}\right):=\left(A_p\left(X_p\right),A_{q}\left(X_{q}\right)\right).
\end{equation*}
Recall by Subsection~\ref{ssec-connfiber}, that the fibred product connection $\ev_0^*A\odot\ev^*A$ is flat if and only if both its components are flat; hence, if $A$ is flat, then $\ev_0^*A\odot\ev^*A$ is also flat.
Analogous results hold for $\ev^*P\odot\ev_0^*P$; denote the connection constructed from $A$ on the fibred product $\ev^*P\odot\ev_0^*P$ by $\ev^*A\odot\ev_0^*A$.
Consider the fibred product map of $\Phi_A$, which gives by Lemma~\ref{lem-fibprodmap} of Subsection~\ref{ssec-connfiber}, a bundle (iso)morphism from $\ev_0^*P\odot\ev^*P$ to $\ev^*P\odot\ev_0^*P$:
\begin{equation*}
\Phi_A\odot\Phi_A^{-1}\left(\gamma;t;p,q\right):=\left(\gamma;t;q\ \holonom{A}{\gamma}{t}{p}{q},p\ \holonom{A}{\gamma}{t}{p}{q}^{-1}\right).
\end{equation*}

After this discussion, one can state and prove the following
\begin{Thm}\label{thm-flatpartransp}
If $A$ flat, the following identity holds:
\begin{equation}\label{eq-flatpartransp}
\left(\Phi_A\odot\Phi_A^{-1}\right)^*\left(\ev^*A\odot\ev_0^*A\right)=\ev_0^*A\odot\ev^*A.
\end{equation}
\end{Thm}

\begin{proof}
It suffices to show that any $\ev_0^*A\odot\ev^*A$-horizontal vector in $\ev_0^*P\odot\ev^*P$ is also horizontal w.r.t.\ the pull-back of $\ev^*A\odot\ev_0^*A$ by the fibred product map $\Phi_A\odot\Phi_A^{-1}$.
One can view a general $\ev_0^*A\odot\ev^*A$-horizontal vector $\left(X_\gamma;X_t;X_p;X_{q}\right)$, where the subscripts are related to the base points of the respective tangent vectors, as the initial tangent direction of an $\ev_0^*A\odot\ev^*A$-horizontal curve, which is denoted by $\left(\gamma_s;\overline{t}_s;p_s,q_s\right)$; one can also assume that $p_s$ and $q_s$ are $A$-horizontal curves in $P$.
Here, $\overline{t}_s$ denotes a curve in $\unint$ starting at $\overline{t}_0=t$.

Denote by $\overline{\gamma}_s$, resp.\ $\widehat{\gamma}_s$, the curve
\[
\overline{\gamma}_s(t)\colon=\gamma_{st}(0),\quad \text{resp.}\quad \widehat{\gamma}_{s}(t)\colon=\gamma_{st}(\overline{t}_{st}).
\]
Since the curve $\left(\gamma_s;\overline{t}_s;p_s,q_s\right)$ belongs by assumption to $\ev_0^*P\odot\ev^*P$, it follows
\[
\pi(p_s)=\gamma_s(0),\quad \pi(q_s)=\gamma_s(t_s),\quad \forall s\in\unint.
\]
Furthermore, define the (piecewise smooth) family of loops in $M$, which is denoted by $\Gamma$, via the assignment
\begin{equation*}
\Gamma(t,s)\colon=\begin{cases}
\overline{\gamma}_s(4t),& t\in \left[0;\frac{1}4\right]\\
\gamma_s\left((4t-1)\overline{t}_s\right),& t\in\left[\frac{1}4,\frac{1}2\right]\\
\widehat{\gamma}_s\left(3-4t\right),& t\in \left[\frac{1}2;\frac{3}4\right]\\
\gamma\left((4-4t)\overline{t}_0\right),& t\in \left[\frac{3}4,1\right].
\end{cases}
\end{equation*}
It is immediate to see that $\Gamma$ is a homotopy of the loop $\gamma_{\overline{t}_0}^{-1}\circ \gamma_{\overline{t}_0}$, based at $\gamma(0)$, where $\gamma_{\overline{t}_0}(s)\colon=\gamma(s\overline{t}_0)$.
Theorem~\ref{thm-holflat} implies then the identity
\[
\holon{A}{\Gamma_s}{p}=\holon{A}{\gamma_{\overline{t}_0}^{-1}\circ \gamma_{\overline{t}_0}}{p}.
\]

Viewing the curve $\Gamma(\bullet,s)$, for any $s\in\unint$, as a composite curve, it follows in virtue of Lemma~\ref{lem-comphol} and~\ref{lem-invhol}:
\begin{align*}
\holon{A}{\Gamma_s}{p}&=\holonom{A}{\gamma_{\overline{t}_0}^{-1}}{1}{q}{p}\ \holonom{A}{\widehat{\gamma}_s^{-1}}{1}{q_s}{q}\\
&\phantom{=}\ \ \holonom{A}{\gamma_s(\bullet \overline{t}_s)}{1}{p_s}{q_s}\ \holonom{A}{\overline{\gamma}_s}{1}{p}{p_s}=\\
&=\holonom{A}{\gamma}{t}{p}{q}^{-1}\holonom{A}{\gamma_s}{t_s}{p_s}{q_s}.
\end{align*}

I made use of the fact that $p_{ts}$, resp.\ $q_{ts}$, (viewed as a curve w.r.t.\ the parameter $t$) is an $A$-horizontal curve over $\gamma_{ts}(0)$, resp.\ $\gamma_{ts}(\overline{t}_{ts})$: in fact, 
\[
\pi(p_{ts})=\gamma_{ts}(0)=\overline{\gamma}_s(t),\quad \pi(q_{ts})=\gamma_{st}(\overline{t}_{ts})=\widehat{\gamma}_s(t).
\]
Moreover, the first curve is based at $p$, while the latter is based at $q$. 
Since both curves are reparametrizations of $A$-horizontal curves, they are also $A$-horizontal, whence it follows
\[
\holonom{A}{\widehat{\gamma}_s}{1}{q}{q_s}=\holonom{A}{\overline{\gamma}_s}{1}{p}{q}=e.
\]
On the other hand, again by Lemma~\ref{lem-comphol} and~\ref{lem-invhol}, one gets
\begin{align*}
\holon{A}{\gamma_{\overline{t}_0}^{-1}\circ \gamma_{\overline{t}_0}}{p}&=\holonom{A}{\gamma_{\overline{t}_0}^{-1}}{1}{q}{p}\ \holonom{A}{\gamma_{\overline{t}_0}}{1}{p}{q}=\\
&=\holonom{A}{\gamma_{\overline{t}_0}}{1}{p}{q}^{-1}\ \holonom{A}{\gamma_{\overline{t}_0}}{1}{p}{q}=\\
&=e.
\end{align*}
Thus, one gets the following result:
\begin{equation}\label{eq-flatpartr}
\holonom{A}{\gamma_s}{\overline{t}_s}{p_s}{q_s}=\holonom{A}{\gamma}{t}{p}{q}.
\end{equation}
Notice finally that the curve $\gamma_{\overline{t}_0}$, resp.\ $\gamma_s(\bullet \overline{t}_s)$, for any $s\in\unint$, is a reparametrization of $\gamma$, resp.\ $\gamma_s$.

Taking the derivative w.r.t.\ $s$ at $s=0$ of (\ref{eq-flatpartr}) yields
\[
\tangent{(\gamma;t;p;q)}{\holo{A}{0}{\bullet}}\!\left(X_\gamma;X_t;X_p;X_{q}\right)=0,
\]
if the tangent vector $\left(X_\gamma;X_t;X_p;X_{q}\right)$ is $\ev_0^*A\odot\ev^*A$-horizontal.

An explicit computation gives, for a general tangent vector $(X_\gamma;X_t;X_p;X_{q})$ to $\ev_0^*P\odot\ev^*P$ at a general point $(\gamma;t;p;q)$
\begin{align*}
&\tangent{(\gamma,t;p,q)}{\left(\Phi_A\odot\Phi_A^{-1}\right)}\!\left(X_\gamma;X_t;X_p;X_{q}\right)=\big(X_\gamma;X_t;\tangent{q}{\widetilde{\raction}_{\holonom{A}{\gamma}{t}{p}{q}}}\!(X_{q})+\big.\\
&\big.+\tangent{\holonom{A}{\gamma}{t}{p}{q}}{\laction_{q}}\!
\left(\tangent{(\gamma;t;p;q)}{\holo{A}{0}{\bullet}}\!\left(X_\gamma;X_t;X_p;X_{q}\right)\right);\tangent{p}{\raction_{\holonom{A}{\gamma}{t}{p}{q}^{-1}}}\!(X_{p})+\big.\\
&\big.+\tangent{\holonom{A}{\gamma}{t}{p}{q}^{-1}}{\laction_{p}}\!
\left(\tangent{(\gamma;t;p;q)}{\left(\holo{A}{0}{\bullet}\right)^{-1}}\!\left(X_\gamma;X_t;X_p;X_{q}\right)\right)\big).
\end{align*}

The previous expression simplifies remarkably, assuming the tangent vector $(X_\gamma;X_t;X_p;X_{q})$ to be $\ev_0^*A\odot\ev^*A$-horizontal:
\begin{align*}
\tangent{(\gamma,t;p,q)}{\left(\Phi_A\odot\Phi_A^{-1}\right)}\!\left(X_\gamma;X_t;X_p;X_{q}\right)&=\big(X_\gamma;X_t;\tangent{q}{\widetilde{\raction}_{\holonom{A}{\gamma}{t}{p}{q}}}\!(X_{q});\big.\\
&\big.\phantom{=}\quad\tangent{p}{\raction_{\holonom{A}{\gamma}{t}{p}{q}^{-1}}}\!(X_{p})\big)
\end{align*}
Hence, a straightforward computation similar in spirit to the final computation in the proof of Theorem~\ref{thm-holflat}, where I make use of the $G$-equivariance of the connection $A$, implies: \

\fbox{\parbox{12cm}{\bf Any $\ev_0^*A\odot \ev^*A$-horizontal vector in the fibred product $\ev_0^*P\odot\ev^*P$ is also $\left(\Phi_A\odot\Phi_A^{-1}\right)^*\!\left(\ev^*A\odot\ev_0^*A\right)$-horizontal.}}

The claim is then a consequence of arguments very similar to those used in the final steps of the proof of Theorem~\ref{thm-flatness}.
\end{proof}
By Lemma~\ref{lem-fibprodconn}, Equation (\ref{eq-flatpartransp}) in Theorem~\ref{thm-flatpartransp} may be rewritten as follows:
\[
\Phi_A^*\ev^*A=\ev_0^*A.
\]
On the other hand, one can rewrite Equation (\ref{eq-flatpartransp}) as a ``horizontality equation'', with help of the following
\begin{Cor}\label{cor-horpartr}
The condition that the isomorphism $\Phi_A\odot\Phi_A^{-1}$ associated to the parallel transport w.r.t.\ $A$ intertwines the fibred product connections $\ev_0^*A\odot\ev^*A$ and $\ev^*A\odot\ev_0^*A$ may be written as a horizontality equation for the parallel transport:
\[
\dd \holo{A}{0}{\bullet}+\ev^*A\ \holo{A}{0}{\bullet}-\holo{A}{0}{\bullet}\ \ev_0^*A=0,
\]
where, by abuse of notation, I denoted by $\ev^*_0A$ and $\ev^*A$ the respective pull-backs w.r.t.\ the projections $\pr_1$ and $\pr_2$ from $\ev_0^*P\odot\ev^*P$ onto $\ev_0^*P$ and $\ev^*P$, respectively.
By abuse of notation, I have also simply denoted as products in the above formula what I should have written as tangent maps of left- and right multiplication; the explicit formulae in the proof will point out to the precise notation I have skipped above. 
\end{Cor}
\begin{proof}
Theorem~\ref{thm-flatpartransp} implies that, for any point $(\gamma;t;p,q)$ in $\ev_0^*P\odot\ev^*P$ and any tangent vector $(X_\gamma;X_t;X_p,X_q)$ at this point to the fibred product bundle, the identity holds:
\begin{multline*}
\left((\Phi_A\odot\Phi_A^{-1})^*\!(\ev^*A\odot\ev_0^*A)\right)_{(\gamma;t;p,q)}\!(X_\gamma;X_t;X_p,X_q)=\\
=(\ev_0^*A\odot\ev^*A)_{(\gamma;t;p,q)}\!(X_\gamma;X_t;X_p,X_q).
\end{multline*}

Recall the expression for the tangent map of the isomorphism $\Phi_A\odot\Phi_A^{-1}$:
\begin{equation}\label{eq-tangisom}
\begin{aligned}
&\tangent{(\gamma,t;p,q)}{\left(\Phi_A\odot\Phi_A^{-1}\right)}\!\left(X_\gamma;X_t;X_p;X_{q}\right)=\big(X_\gamma;X_t;\tangent{q}{\widetilde{\raction}_{\holonom{A}{\gamma}{t}{p}{q}}}\!(X_{q})+\big.\\
&\big.+\tangent{\holonom{A}{\gamma}{t}{p}{q}}{\laction_{q}}\!
\left(\tangent{(\gamma;t;p;q)}{\holo{A}{0}{\bullet}}\!\left(X_\gamma;X_t;X_p;X_{q}\right)\right);\tangent{p}{\raction_{\holonom{A}{\gamma}{t}{p}{q}^{-1}}}\!(X_{p})+\big.\\
&\big.+\tangent{\holonom{A}{\gamma}{t}{p}{q}^{-1}}{\laction_{p}}\!
\left(\tangent{(\gamma;t;p;q)}{\left(\holo{A}{0}{\bullet}\right)^{-1}}\!\left(X_\gamma;X_t;X_p;X_{q}\right)\right)\big).
\end{aligned}
\end{equation}
Notice that the right-hand side of Equation (\ref{eq-tangisom}) can be splitted in two tangent vectors, namely
\[
\left(X_\gamma;X_t;\tangent{q}{\widetilde{\raction}_{\holonom{A}{\gamma}{t}{p}{q}}}\!(X_{q}),\tangent{p}{\raction_{\holonom{A}{\gamma}{t}{p}{q}^{-1}}}\!(X_{p})\right)
\]
and
\begin{multline*}
\bigg(0;0;\tangent{\holonom{A}{\gamma}{t}{p}{q}}{\laction_{q}}\!
\left(\tangent{(\gamma;t;p;q)}{\holo{A}{0}{\bullet}}\!\left(X_\gamma;X_t;X_p;X_{q}\right)\right);\bigg.\\
\bigg.\tangent{\holonom{A}{\gamma}{t}{p}{q}^{-1}}{\laction_{p}}\!
\left(\tangent{(\gamma;t;p;q)}{\left(\holo{A}{0}{\bullet}\right)^{-1}}\!\left(X_\gamma;X_t;X_p;X_{q}\right)\right)\bigg).
\end{multline*}

The tangent vector $\tangent{\holonom{A}{\gamma}{t}{p}{q}}{\laction_{q}}\!
\left(\tangent{(\gamma;t;p;q)}{\holo{A}{0}{\bullet}}\!\left(X_\gamma;X_t;X_p;X_{q}\right)\right)$ can be then rewritten as
\begin{align*}
&\tangent{\holonom{A}{\gamma}{t}{p}{q}}{\laction_{q}}\!
\left(\tangent{(\gamma;t;p;q)}{\holo{A}{0}{\bullet}}\!\left(X_\gamma;X_t;X_p;X_{q}\right)\right)=\\
&=\tangent{e}{\laction_{q\holonom{A}{\gamma}{t}{p}{q}}}\!\left\{\tangent{\holonom{A}{\gamma}{t}{p}{q}}{\laction_{\holonom{A}{\gamma}{t}{p}{q}^{-1}}}\right.\\
&\phantom{=}\left.\left[\tangent{(\gamma;t;p;q)}{\holo{A}{0}{\bullet}}\!\left(X_\gamma;X_t;X_p;X_{q}\right)\right]\right\};
\end{align*}
similarly, one has
\begin{align*}
&\tangent{\holonom{A}{\gamma}{t}{p}{q}^{-1}}{\laction_{p}}\!\left(\tangent{(\gamma;t;p;q)}{\left(\holo{A}{0}{\bullet}\right)^{-1}}\!\left(X_\gamma;X_t;X_p;X_{q}\right)\right)=\\
&=\tangent{e}{\laction_{q\holonom{A}{\gamma}{t}{p}{q}^{-1}}}\!\left\{\tangent{\holonom{A}{\gamma}{t}{p}{q}^{-1}}{\laction_{\holonom{A}{\gamma}{t}{p}{q}}}\right.\\
&\phantom{=}\left.\left[\tangent{(\gamma;t;p;q)}{(\holo{A}{0}{\bullet})^{-1}}\!\left(X_\gamma;X_t;X_p;X_{q}\right)\right]\right\}.
\end{align*}

Thus, the second vector in the previous splitting of Equation (\ref{eq-tangisom}) is vertical, whence one can deduce from the definition of fibred product connection and from the first equation at the beginning of the proof the pair of equations:
\begin{align*}
A_p(X_p)&=A_{q\holo{A}{0}{\bullet}\!(\gamma;t;p,q)}\!(\tangent{q}{\raction_{\holo{A}{0}{\bullet}\!(\gamma;t;p,q)}}\!(X_q))+\\
&\phantom{=}+A_{q\holo{A}{0}{\bullet}\!(\gamma;t;p,q)}\!\bigg(\tangent{e}{\laction_{q\holo{A}{0}{\bullet}\!(\gamma;t;p,q)}}\!\left\{\tangent{\holo{A}{0}{\bullet}\!(\gamma;t;p,q)}{\laction_{(\holo{A}{0}{\bullet})^{-1}\!(\gamma;t;p,q)}}\right.\bigg.\\
&\phantom{=+}\bigg.\left.\left[\tangent{(\gamma;t;p;q)}{\holo{A}{0}{\bullet}}\!\left(X_\gamma;X_t;X_p;X_{q}\right)\right]\right\}\bigg)=\\
&=\Ad\!\left((\holo{A}{0}{\bullet})^{-1}\!(\gamma;t;p,q)\right)A_q(X_q)+\\
&\phantom{=}+\tangent{\holo{A}{0}{\bullet}\!(\gamma;t;p,q)}{\laction_{(\holo{A}{0}{\bullet})^{-1}\!(\gamma;t;p,q)}}\!\left[\tangent{(\gamma;t;p;q)}{\holo{A}{0}{\bullet}}\!\left(X_\gamma;X_t;X_p;X_{q}\right)\right]
\end{align*}
and
\begin{align*}
A_q(X_q)&=A_{p(\holo{A}{0}{\bullet})^{-1}\!(\gamma;t;p,q)}\!(\tangent{p}{\raction_{(\holo{A}{0}{\bullet})^{-1}\!(\gamma;t;p,q)}}\!(X_p))+\\
&+A_{p(\holo{A}{0}{\bullet})^{-1}\!(\gamma;t;p,q)}\!\bigg(\tangent{e}{\laction_{p(\holo{A}{0}{\bullet})^{-1}\!(\gamma;t;p,q)}}\!\left\{\tangent{(\holo{A}{0}{\bullet})^{-1}\!(\gamma;t;p,q)}{\laction_{(\holo{A}{0}{\bullet})\!(\gamma;t;p,q)}}\right.\bigg.\\
&\phantom{=+}\bigg.\left.\left[\tangent{(\gamma;t;p;q)}{(\holo{A}{0}{\bullet})^{-1}}\!\left(X_\gamma;X_t;X_p;X_{q}\right)\right]\right\}\bigg)=\\
&=\Ad\!\left((\holo{A}{0}{\bullet})\!(\gamma;t;p,q)\right)A_p(X_p)+\\
&\phantom{=}+\tangent{(\holo{A}{0}{\bullet})^{-1}\!(\gamma;t;p,q)}{\laction_{(\holo{A}{0}{\bullet})\!(\gamma;t;p,q)}}\left[\tangent{(\gamma;t;p;q)}{(\holo{A}{0}{\bullet})^{-1}}\!\left(X_\gamma;X_t;X_p;X_{q}\right)\right].
\end{align*}
Observe that I used the $G$-equivariance and the verticality of $A$ in the second identities.

It is not difficult to see that both equations are equivalent, and each one is in turn equivalent to the horizontality equation for the parallel transport, thus yielding the claim.
\end{proof}

\section{Applications of the generalized gauge transformation associated to the parallel transport: iterated Chen integrals}\label{sec-chenint}
In this section, I want to discuss an application of the computations performed in Section~\ref{sec-flatness}, which was already mentioned in the Introduction and which was the main motivation for this paper, namely iterated Chen integrals in the spirit of~\cite{CCL} and~\cite{CR}, where I discussed so-called {\em generalized holonomies} as differential forms on the space of loops $\lo M$.
These differential forms are constructed via iterated Chen integrals in the sense of~\cite{Ch}; they are obtained by taking push-forward (i.e.\ integration along the fibre of trivial fibre bundles, whose fibres are simplices $\triangle_n$, for any $n\in\mathbb{N}$) of traces in some representation of the Lie algebra $\Lg$ of the structure group $G$ of products of differential forms on $\lo M\times \triangle_n$ with values in $\Lg$ (or, to be more precise, in the universal enveloping algebra of $\Lg$), which are in turn pull-backs of differential forms on $M$ with values in $\Lg$ w.r.t.\ evaluation maps at different points of the simplex $\triangle_n$.
As it is well-known, differential forms on $M$ with values in the Lie algebra $\Lg$ of $G$ are in one-to-one correspondence with differential forms on $M$ with values in the vector bundle associated to the trivial principal $G$-bundle over $M$ via the adjoint representation of $G$.

The main question is: is it possible to consider analogues of the aforementioned Chen-type integrals, when considering a {\em nontrivial} principal $G$-bundle $P$, and consequently differential forms on $M$ with values in the associated bundle $\ad P$ to $P$ via the adjoint representation of $G$?  
The answer is positive, 

\subsection{Parallel transport and simplices}\label{ssec-simplpartr}
I am therefore interested in the following more general situation: instead of considering the ``cylinder'' $\lo M\times \unint$, consider the more general space $\lo M\times \triangle_n$, where $\triangle_n$ denotes the standard $n$-dimensional simplex
\[
\triangle_n\colon=\left\{\left(t_1,\dots,t_n\right)\in \unint^n\colon 0\leq t_1\leq \cdots\leq t_n\leq 1\right\}.
\]
Notice that $\unint=\triangle_1$.

There are $n+1$ natural maps from $\lo M\times \triangle_n$ into $M$, given by
\[
\ev_{i,n}\colon=\ev\circ \pi_{i,n},\quad 1\leq i\leq n,\quad \ev(0)_n\colon=\ev(0)\circ \pi_n,
\]
where $\pi_{i,n}\left(\gamma;t_1,\dots,t_n\right)\colon=\left(\gamma;t_i\right)$ and $\pi_n\left(\gamma;t_1,\dots,t_n\right)\colon=\gamma$.

Accordingly to the definitions of Section~\ref{sec-applpullbun}, define the pulled-back bundles $\ev_{i,n}^*P$ and $\ev(0)_n^*P$ over $\lo M\times \triangle_n$, as well as their fibred product $\ev(0)_n^*P\odot\ev_{i,n}^*P$ and $\ev_{i,n}^*P\odot \ev_{j,n}^*P$, for any two indices $1\leq i<j\leq n$; the natural maps from $\ev(0)_n^*P$ and $\ev_{i,n}^*P$ to $P$ are denoted by $\widetilde{\ev(0)_n}$ and $\widetilde{\ev_{i,n}}$ respectively. 

Given a connection $A$ on $P$, natural connections on the fibred products $\ev(0)_n^*P\odot\ev_{i,n}^*P$ and $\ev_{i,n}^*P\odot\ev_{j,n}^*P$ are given by the fibred product connections $\ev(0)_n^*A\odot\ev_{i,n}^*A$ and $\ev_{i,n}^*A\odot\ev_{j,n}^*A$; denote by $\ev_{i,n}^*A\odot\ev(0)_n^*A$ the corresponding fibred product connection on $\ev_{i,n}^*P\odot\ev(0)_n^*P$. 
(As before, I changed slightly the notation: in fact, I pulled the connection $A$ w.r.t.\ the corresponding bundle maps $\widetilde{\ev_{i,n}}$. resp.\ $\widetilde{\ev(0)_n}$.)
\begin{Prop}\label{prop-simplpartr}
For any connection $A$ on $P$, the {\em parallel transport map $\holo{A}{0}{i}$} defined by
\[
\ev(0)_n^*P\odot \ev_{i,n}^*P\ni\left(\gamma;t_1,\dots,t_n;p_0,p_i\right)\mapsto \holonom{A}{\gamma}{t_i}{p_0}{p_i}\in G,
\]
is a generalized gauge transformation between the bundle $\ev(0)_n^*P$ and $\ev_{i,n}^*P$. 
Moreover, the triple $\left(\holo{\bullet}{0}{i},\Phi_{P,i},\varphi_{P,i}\right)$, where $\holo{\bullet}{0}{i}$ is viewed as a map from the space of connections on $P$ to the generalized gauge transformations between $\ev(0)_n^*P$ and $\ev_{i,n}^*P$, is an equivariant map from the space of connections, viewed as a $\calG_P$-space, to the groupoid $\calG_{G,\lo M\times \triangle_n}$ of generalized gauge transformations on $\lo M\times\triangle_n$.
Here, the morphism $\left(\Phi_{P,i},\varphi_{P,i}\right)$ is defined as follows
\[
\Phi_{P,i}(\sigma)\colon=\left(\sigma_{\ev(0)_n},\sigma_{\ev_{i,n}}\right),\quad \varphi_{P,i}(P)\colon=\left(\ev(0)_n^*P,\ev_{i,n}^*P\right),
\] 
where I used the same notations of Lemma~\ref{lem-pullgauge} and Lemma~\ref{lem-princmorphism} of Section~\ref{sec-applpullbun}.
\end{Prop}
\begin{proof}
The proof is a consequence of the very same arguments used in the proof of Proposition~\ref{prop-equivpartr}.
\end{proof}
As a consequence of Proposition~\ref{prop-simplpartr} together with Theorem~\ref{thm-isomequiv}, the bundle $\ev(0)_n^*P$ is isomorphic to $\ev_{i,n}^*P$, for any index $i$, via  
\[
\Phi_{A,i,n}\left(\gamma;t_1,\dots,t_n;p\right)\colon=\left(\gamma;t_1,\dots,t_n;q\ \holonom{A}{\gamma}{t_i}{p}{q}\right),
\]
where $\pi(p)=\gamma(0)$, and $q$ is any point in $P$ over $\gamma(t_i)$.

Consider now two indices $1\leq i<j\leq n$ and the corresponding bundles $\ev_{i,n}^*P$, $\ev_{j,n}^*P$, and their fibred products $\ev_{i,n}^*P\odot\ev_{j,n}^*P$ and $\ev_{j,n}^*P\odot\ev_{i,n}^*P$.  

For a connection $A$ on $P$, define the following parallel transport map, which is denoted by $\holo{A}{i}{j}$,
\begin{equation}\label{eq-partrij}
\ev_{i,n}^*P\odot\ev_{j,n}^*P\ni\left(\gamma;t_1,\dots,t_n;p_i,p_j\right)\mapsto\holonomy{A}{\gamma}{t_i}{t_j}{p_i}{p_j}\in G. 
\end{equation}
Corollary~\ref{cor-proppartr} implies the following
\begin{Prop}\label{prop-simplpartr2}
For any connection $A$ on $P$, the parallel transport map $\holo{A}{i}{j}$ defined by Equation (\ref{eq-partrij}), for any two indices $1\leq i<j\leq n$, is a generalized gauge transformation between the bundles $\ev_{i,n}^*P$ and $\ev_{j,n}^*P$.
Moreover, the triple $\left(\holo{\bullet}{i}{j},\Phi_{P,i,j},\varphi_{P,i,j}\right)$, where $\holo{\bullet}{i}{j}$ is viewed as a map from the space of connections on $P$ to the generalized gauge transformations between $\ev_{i,n}^*P$ and $\ev_{j,n}^*P$, is an equivariant map from the space of connections as a $\calG_P$-space, to the groupoid $\calG_{G,\lo M\times \triangle_N}$ of generalized gauge transformations on $\lo M\times\triangle_n$ w.r.t.\ generalized conjugation.
The morphism $\left(\Phi_{P,i,j},\varphi_{P,i,j}\right)$ is defined in this case via
\[
\Phi_{P,i,j}(\sigma)\colon=\left(\sigma_{\ev_{i,n}},\sigma_{\ev_{j,n}}\right),\quad \varphi_{P,i,j}(P)\colon=\left(\ev_{i,n}^*P,\ev_{j,n}^*P\right),
\]
where I used the same notations of Lemma~\ref{lem-pullgauge} and Lemma~\ref{lem-princmorphism} of Section~\ref{sec-applpullbun}.
\end{Prop}
According to Theorem~\ref{thm-isomequiv}, there is a bundle (iso)morphism $\Phi_{A,i,j,n}$ from $\ev_{i,n}^*P$ to $\ev_{j,n}^*P$, defined via
\[
\ev_{i,n}^*P\ni\left(\gamma;t_1,\dots,t_n;p\right)\overset{\Phi_{A,i,j,n}}\mapsto\left(\gamma;t_1,\dots,t_n;q\holonomy{A}{\gamma}{t_i}{t_j}{p}{q}\right)\in \ev_{j,n}^*P,
\]
where $q$ is any point obeying $\pi(q)=\gamma(t_j)$.

\subsection{Flatness and parallel transport on simplices}\label{ssec-holpartrsimpl} 
In this subsection, I want to discuss analogues of Theorem~\ref{thm-flatpartransp} of Subsection~\ref{ssec-holpartr} for the isomorphisms $\Phi_{A,i,n}$ and $\Phi_{A,i,j,n}$ introduced and discussed in the previous subsection.
I begin therefore  with a discussion of verticality and $\ev(0)_n^*A\odot\ev_{i,n}^*A$-horizontality of tangent vectors to the fibred product $\ev(0)_n^*P\odot\ev_{i,n}^*P$, resp.\ $\ev_{i,n}^*P\odot\ev_{j,n}^*P$, for any two indices $1\leq i<j\leq n$.

A tangent vector to $\ev(0)_n^*P\odot\ev_{i,n}^*P$ at the point $\left(\gamma;t_1,\dots,t_n;p,q\right)$ is given by a $n+3$-tuple $\left(X_\gamma;X_1,\dots,X_n;X_p,X_q\right)$, where $i)$ $X_\gamma$ is a tangent vector to $\lo M$ at $\gamma$ (hence, a vector field on the pull-back of the tangent bundle of $M$ w.r.t.\ $\gamma$), $ii)$ $X_i$ a tangent vector to $\unint$ at $t_i$, which is the projection on the $i$-th factor of the $n$-simplex $\triangle_n$, and $iii)$ $X_p$, resp.\ $X_q$, is a tangent vector to $P$ at $p$, resp.\ $q$.

The condition on $\left(X_\gamma;X_1,\dots,X_n;X_p,X_q\right)$ to be vertical, by the very definition of the bundle projection onto $\lo M\times \triangle_n$, is encoded in the following equation:
\[
X_\gamma=0,\ X_i=0,\quad \forall 1\leq i\leq n\Rightarrow\begin{cases}
\tangent{p}{\pi}\!(X_p)&=\tangent{(\gamma;t_1,\dots)}{\ev(0)_n}\!(X_\gamma;X_1,\dots)=0,\\
\tangent{q}{\pi}\!(X_q)&=\tangent{(\gamma;t_1,\dots)}{\ev_{i,n}}\!(X_\gamma;X_1,\dots)=0.
\end{cases}
\]  
Hence, a general tangent vector $\left(X_\gamma;X_1,\dots,X_n;X_p,X_q\right)$ to $\ev(0)_n^*P\odot\ev_{i,n}^*P$ at the point $\left(\gamma;t_1,\dots,t_n;p,q\right)$ is vertical if and only if its $\lo M\times \triangle_n$-\-piece vanishes and both its $P$-pieces are vertical in $P$.

The choice of a connection $A$ on $P$ determines a unique splitting of any tangent vector to $P$ in its vertical part and $A$-horizontal part; hence, the two $P$-pieces of any tangent vector $\left(X_\gamma;X_1,\dots,X_n;X_p,X_q\right)$ to $\ev(0)_n^*P\odot\ev_{i,n}^*P$ at the point $\left(\gamma;t_1,\dots,t_n;p,q\right)$ may be decomposed in a unique way as
\[
X_p=X_p^v+X_p^h,\quad X_q=X_q^v+X_q^h.
\]

\fbox{\parbox{12cm}{\bf Thus, one has the unique decomposition of $\left(X_\gamma;X_1,\dots,X_n;X_p,X_q\right)$:
\[
\left(X_\gamma;X_1,\dots,X_n;X_p,X_q\right)=\left(0;0,\dots,0;X_p^v,X_q^v\right)+\left(X_\gamma;X_1,\dots,X_n;X_p^h,X_q^h\right).
\]
It is easy to verify, by the very definition of fibred product connection $\ev(0)_n^*A\odot\ev_{i,n}^*A$, that the tangent vector $\left(X_\gamma;X_1,\dots,X_n;X_p^h,X_q^h\right)$ is $\ev(0)_n^*A\odot\ev_{i,n}^*A$-horizontal.}}

Analogously, when considering the fibred product $\ev_{i,n}^*P\odot\ev_{j,n}^*P$, one gets \ 

\fbox{\parbox{12cm}{\bf Given a connection $A$ on $P$, there is a unique splitting of any tangent vector $\left(X_\gamma;X_1,\dots,X_n;X_p,X_q\right)$ to $\ev_{i,n}^*P\odot\ev_{j,n}^*P$:
\[
\left(X_\gamma;X_1,\dots,X_n;X_p,X_q\right)=\left(0;0,\dots,0;X_p^v,X_q^v\right)+\left(X_\gamma;X_1,\dots,X_n;X_p^h,X_q^h\right),
\]
where $X_p^h$, resp.\ $X_q^h$, is an $A$-horizontal vector to $P$ at $p$, resp.\ $q$.
The second piece in the above identity is obviously $\ev_{i,n}^*A\odot\ev_{j,n}^*A$-horizontal, by the very definition of fibred product connection.}}

Consider the fibred product map $\Phi_{A,i,n}\odot\Phi_{A,i,n}^{-1}$, resp.\ $\Phi_{A,i,j,n}\odot\Phi_{A,i,j,n}^{-1}$, which is an isomorphism from $\ev(0)_n^*P\odot\ev_{i,n}^*P$ to $\ev_{i,n}^*P\odot\ev(0)_n^*P$, resp.\ from $\ev_{i,n}^*\odot\ev_{j,n}^*P$ to $\ev_{j,n}^*P\odot\ev_{i,n}^*P$; explicitly, these maps take the form
\begin{align*}
\left(\Phi_{A,i,n}\odot\Phi_{A,i,n}^{-1}\right)\!\left(\gamma;t_1,\dots;p,q\right)&=\left(\gamma;t_1,\dots;q\holonom{A}{\gamma}{t_i}{p}{q},p\holonom{A}{\gamma}{t_i}{p}{q}^{-1}\right)\\
\left(\Phi_{A,i,j,n}\odot\Phi_{A,i,j,n}^{-1}\right)\!\left(\gamma;t_1,\dots;p,q\right)&=\Big(\gamma;t_1,\dots;q\holonomy{A}{\gamma}{t_i}{t_j}{p}{q},\Big.\\
&\qquad \Big.p\holonomy{A}{\gamma}{t_i}{t_j}{p}{q}^{-1}\Big). 
\end{align*}

The same arguments used in the proof of Theorem~\ref{thm-flatpartransp} yield the following
\begin{Thm}\label{thm-flatpartransp1}
If $A$ is flat, the following identities hold:
\begin{align*}
\left(\Phi_{A,i,n}\odot\Phi_{A,i,n}^{-1}\right)^*\!\left(\ev_{i,n}^*A\odot\ev(0)_n^*A\right)&=\ev(0)_n^*A\odot\ev_{i,n}^*A,\\
\left(\Phi_{A,i,j,n}\odot\Phi_{A,i,j,n}^{-1}\right)^*\!\left(\ev_{j,n}^*A\odot\ev_{i,n}^*A\right)&=\ev_{i,n}^*A\odot\ev_{j,n}^*A.
\end{align*}
\end{Thm} 

Theorem~\ref{thm-flatpartransp1} together with Lemma~\ref{lem-fibprodconn} implies the identities:
\begin{align*}
\Phi_{A,i,n}^*\left(\ev_{i,n}^*A\right)&=\ev(0)_n^*A,\\
\Phi_{A,i,j,n}^*\left(\ev_{j,n}^*A\right)&=\ev_{i,n}^*A,
\end{align*}
if $A$ is flat.
Hence, Theorem~\ref{thm-flatpartransp1} shows that the isomorphism $\Phi_{A,i,n}$, resp.\ $\Phi_{A,i,j,n}$, intertwines the connections $\ev_{i,n}^*A$ and $\ev(0)_n^*A$, resp.\ $\ev_{j,n}^*P$ and $\ev_{i,n}^*P$.

\subsection{Iterated Chen-type integrals}\label{ssec-chenint}
I consider now a differential form $\omega$ on $M$ with values in the associated bundle $\ad P=P\times_G\Lg$ of $P$, where $\Lg$, the Lie algebra of $G$, becomes a representation of $G$ w.r.t. the adjoint action; by the same label I denote the unique basic form on $P$ with values in $\Lg$ corresponding to $\omega$.
One can take the pull-back of $\omega$ w.r.t.\ the bundle map $\widetilde{\ev_{i,n}}$ and, as a result, one gets a basic form on $\ev_{i,n}^*P$ with values in the Lie algebra $\Lg$, which is denoted (by abuse of notations) by $\ev_{i,n}^*\omega$.
Taking then the pull-back of $\ev_{i,n}^*\omega$ w.r.t.\ $\Phi_{A,i,n}$, one gets a basic form on $\ev(0)_n^*P$ with values in $\Lg$, which in turn descends to a form on $\lo M\times \triangle_n$ with values in the associated bundle $\ad \ev(0)_n^*P$.
Denote the result of these three operations by $\widehat{\omega}_{A,i,n}$, for $1\leq i\leq n$.
In other words, as a basic form on $\ev(0)_n^*P$ with values in $\Lg$, one has
\[
\widehat{\omega}_{A,i,n}\colon=\left(\widetilde{\ev_{i,n}}\circ\Phi_{A,i,n}\right)^*\omega.
\]
I want now to analyze in detail, using the tools developed in the previous sections about the isomorphisms induced from parallel transport, the properties of the pulled-back forms $\widehat{\omega}_{A,i,n}$, i.e.\ its dependence on the space of connections $\calA_P$, which may be summarised in the following
\begin{Prop}\label{prop-covgauge}
\begin{itemize}
\item[i)] The assignment
\begin{align*}
&\calA_P\times\Omega^*\!\left(M,\ad P\right)\ni (A,\omega)\overset{\Psi_{0,i,n}}\mapsto\\
\overset{\Psi_{0,i,n}}\mapsto&(\ev(0)_n^*A,\widehat{\omega}_{A,i,n})\in\calA_{\ev(0)_n^*P}\times\Omega^*\!\left(\lo M\times \triangle_n,\ad \ev(0)_n^*P\right)
\end{align*}
is $\calG_P$-equivariant; the gauge group $\calG_P$ operates on the space of forms on $\lo M\times\triangle_n$ with values in the associated bundle $\ad \ev(0)_n^*P$ by composing the action of $\calG_{\ev(0)_n^*P}$ with the group homomorphism from $\calG_P$ to $\calG_{\ev(0)_n^*P}$ described in Lemma~\ref{lem-pullgauge} of Section~\ref{sec-applpullbun}.

\item[ii)] There is a natural map from $\calA_P\times\Omega^*\!\left(M,\ad P\right)$ to $\Omega^*\!\left(M,\ad P\right)$, resp.\ from the product $\calA_{\ev(0)_n^*P}\times\Omega^*\!\left(\lo M\times \triangle_n,\ad \ev(0)_n^*P\right)$ to $\Omega^*\!\left(\lo M\times \triangle_n,\ad \ev(0)_n^*P\right)$, given by
\[
(A,\omega)\overset{\mathnormal{cov}_P}\mapsto \dd_A\omega,\quad\text{resp.}\quad (B,\eta)\overset{\mathnormal{cov}_{\ev(0)_n^*P}}\mapsto \dd_B\eta.
\] 
Restricting to the subspace $F^{-1}(\left\{0\right\})$ of flat connections, the following diagram commutes:
\[
\begin{CD} 
F^{-1}(\left\{0\right\})\times\Omega^*\!\left(M,\ad P\right)  @>\Psi_{0,i,n}>>F^{-1}(\left\{0\right\})\times\Omega^*\!\left(\lo M\times \triangle_n,\ad \ev(0)_n^*P\right)\\
@V\pr_{\calA_P}\times\mathnormal{cov}_P VV              @VV\pr_{\calA_{\ev(0)_n^*P}}\times\mathnormal{cov}_{\ev(0)_n^*P} V\\
F^{-1}(\left\{0\right\})\times \Omega^*\!\left(M,\ad P\right) @>\Psi_{0,i,n}>>F^{-1}(\left\{0\right\}) \times\Omega^*\!\left(\lo M\times\triangle_n,\ad \ev(0)_n^*P\right)  
\end{CD}\quad.
\] 

As a consequence of the commutativity of the previous diagram, one can consider the ``cohomology bundle'' $\coh\!(P,\Lg)$ over the moduli space of (irreducible) flat connections, whose fibre at a class $[A]$ consists of all cohomology classes with values in the associated bundle $\ad P$ at any degree of the covariant differential $\dd_A$, and the ``cohomology bundle'' $\coh\!(\ev(0)_n^*P,\Lg)$, with fibre at $[B]$ consisting of cohomology classes with values in the associated bundle $\ad \ev(0)_n^*P$ w.r.t.\ the covariant differential $\dd_B$; then, the map $\Psi_{0,i,n}$ descends to a bundle morphism between the cohomology bundles. 
(To be more precise, the cohomology bundles may be defined over any open subset of $\calM_P$, where the dimensions of the aforementioned ``covariant'' cohomology groups are constant; therefore, one has to consider a family of cohomology bundles over all such open sets, where they make sense.)
\end{itemize}
\end{Prop}
\begin{proof}
The gauge group $\calG_P$ of $P$ defines by pull-back a right action on $\Omega^*\left(M,\ad P\right)$, viewed as the space of basic differential forms on $P$ with values in $\Lg$:
\[
\omega^\sigma\colon=\sigma^*\!\left(\omega\right),\quad \sigma\in\calG
\]
Proposition~\ref{prop-equivpartr} implies, using again the same notations of Section~\ref{sec-applpullbun},
\begin{align*}
\Psi_{0,i,n}\left((A,\omega)^\sigma\right)&=\widehat{\omega^\sigma}_{A^\sigma,i,n}=\\
&=\left(\sigma\circ\widetilde{\ev_{i,n}}\circ\Phi_{A^\sigma,i,n}\right)^*\omega=\\
&=\left(\widetilde{\ev_{i,n}}\circ\sigma_{\ev_{i,n}}\circ\Phi_{A^\sigma,i,n}\right)^*\omega=\\
&=\left(\widetilde{\ev_{i,n}}\circ\Phi_{A,i,n}\circ\sigma_{\ev(0)_n}\right)^*\omega=\\
&=\widehat{\omega}_{A,i,n}^{\sigma_{\ev(0)_n}}=\\
&=\left(\Psi_{0,i,n}(A,\omega)\right)^{\sigma_{\ev(0)_n}},
\end{align*}
from which the first assertion of the Proposition follows.

As for the second assertion, assume the connection $A$ on $P$ to be flat; then one gets:
\begin{align*}
\left(\mathnormal{cov}_{\ev(0)_n^*P}\circ\Psi_{0,i,n}\right)\!(A,\omega)&=\dd_{\ev(0)_n^*A}\widehat{\omega}_{A,i,n}=\\
&=\dd \widehat{\omega}_{A,i,n}+\Lie{\left(\ev(0)_n^*A\right)}{\widehat{\omega}_{A,i,n}}=\\
&=\widehat{\dd \omega}_{A,i,n}+\Lie{\left(\ev(0)_n^*A\right)}{\left(\Phi_{A,i,n}\circ\widetilde{\ev_{i,n}}\right)^*\omega}=\\
&=\widehat{\dd \omega}_{A,i,n}+\left(\Phi_{A,i,n}\circ\widetilde{\ev_{i,n}}\right)^*\Lie{A}{\omega}=\\
&=\widehat{\dd_A\omega}_{A,i,n}=\\
&=\left(\pr_2\circ\Psi_{0,i,n}\circ (\pr_{\calA_P}\times \mathnormal{cov}_P)\right)\!(A,\omega),
\end{align*}
where $\pr_2$ is the projection onto the second factor of $\calA_{\ev(0)_n^*P}$ with the forms on $\lo M\times \triangle_n$ with values in the associated bundle $\End(V)$.
The third equality is a direct consequence of Theorem~\ref{thm-flatpartransp1}, if $A$ is flat.
\end{proof}

Let me state and prove a Lemma useful for explicit computations.
\begin{Lem}\label{lem-explpullhol}
For any connection $A$ in $\calA_P$, for any basic form $\omega$ on $P$ with values in $\Lg$ and for any index $1\leq i\leq n$, the following formula holds:
\[
\pr_1^*\!\left(\widehat{\omega}_{A,i,n}\right)=\Ad\!\left((\holo{A}{0}{i})^{-1}\right)\pr_2^*\!\left(\ev_{i,n}^*\omega\right),
\]
where $\pr_j$ denotes the $j$-projection from $\ev(0)_n^*P\odot\ev_{i,n}^*P$ onto its $j$-factor, for $j=1,2$.
\end{Lem}
\begin{proof}
Assuming $\omega$ to be of degree $r$, consider $r$ tangent vectors $X^j\colon=\left(X_\gamma^j;X_1^j,\dots;X_p^j,X_q^j\right)$ to $\ev(0)_n^*P\odot \ev_{i,n}^*P$ at a general point.

Thus, by definition of pull-back, one has
\begin{align*}
\pr_1^*\!\left(\widehat{\omega}_{A,i,n}\right)_{\left(\gamma;t_1,\dots;p,q\right)}\!\left(X^1,\dots\right)&=\left(\widehat{\omega}_{A,i,n}\right)_{\left(\gamma;t_1,\dots;p\right)}\!\left(\tangent{\left(\gamma;t_1,\dots;p,q\right)}{\pr_1}\!\left(X^1\right),\dots\right)=\\
&=\omega_{q\holonom{A}{\gamma}{t_i}{p}{q}}\!\left(\tangent{\left(\gamma;t_1,\dots;p,q\right)}{\left(\widetilde{\ev_{i,n}}\circ \Phi_{A,i,n}\circ \pr_1\right)}\!\left(X^1\right),\right.\\
&\phantom{=\quad}\left.\dots,\tangent{\left(\gamma;t_1,\dots;p,q\right)}{\left(\widetilde{\ev_{i,n}}\circ \Phi_{A,i,n}\circ \pr_1\right)}\!\left(X^r\right)\right).
\end{align*}

An easy computation gives the following expression for the tangent map of the composite map $\widetilde{\ev_{i,n}}\circ \Phi_{A,i,n}\circ \pr_1$ at $(\gamma;t_1,\dots;p,q)$ applied to any tangent vector $X^j$:
\begin{multline*}
\tangent{\left(\gamma;t_1,\dots;p,q\right)}{\left(\widetilde{\ev_{i,n}}\circ \Phi_{A,i,n}\circ \pr_1\right)}\!\left(X^j\right)=\tangent{q}{\raction_{\holo{A}{0}{i}\!\left(\gamma;t_1,\dots;p,q\right)}}\!\left(X_q^j\right)+\\
+\tangent{\holo{A}{0}{i}\!\left(\gamma;t_1,\dots;p,q\right)}{\laction_{q}}\!\left(\tangent{\left(\gamma;t_1,\dots;p,q\right)}{\holo{A}{0}{i}}\!\left(X^j\right)\right).
\end{multline*}

The second term of the previous equation may be also rewritten as
\begin{align*}
&\tangent{\holo{A}{0}{i}\!\left(\gamma;t_1,\dots;p,q\right)}{\laction_{q}}\!\left(\tangent{\left(\gamma;t_1,\dots;p,q\right)}{\holo{A}{0}{i}}\!\left(X^j\right)\right)=\\
&\tangent{e}{\laction_{q\holo{A}{0}{i}\!\left(\gamma;t_1,\dots;p,q\right)}}\!\bigg\{\tangent{\holo{A}{0}{i}\!\left(\gamma;t_1,\dots;p,q\right)}{\laction_{\left(\holo{A}{0}{i}\!\left(\gamma;t_1,\dots;p,q\right)\right)^{-1}}}\bigg.\\
&\bigg.\left[\tangent{\left(\gamma;t_1,\dots;p,q\right)}{\holo{A}{0}{i}}\!\left(X^j\right)\right]\bigg\}.
\end{align*}
Thus, the tangent map of $\widetilde{\ev_{i,n}}\circ \Phi_{A,i,n}\circ \pr_1$ at a general point $\left(\gamma;t_1,\dots;p,q\right)$ of $\ev(0)_n^*P\odot \ev_{i,n}^*P$ applied to any tangent vector $X^j$ splits into a sum of two vectors, one of which being vertical in $P$.

Since $\omega$ is basic, it follows
\begin{align*}
\pr_1^*\!\left(\widehat{\omega}_{A,i,n}\right)_{\left(\gamma;t_1,\dots;p,q\right)}\!\left(X^1,\dots\right)&=\omega_{q\holo{A}{0}{i}\!\left(\gamma;t_1,\dots;p,q\right)}\!\left(\tangent{q}{\raction_{\holo{A}{0}{i}\!\left(\gamma;t_1,\dots;p,q\right)}}\!\left(X_q^j\right),\dots\right)=\\
&=\Ad\!\left(\left(\holo{A}{0}{i}\right)^{-1}\!\left(\gamma;t_1,\dots;p,q\right)\right)\omega_q\!\left(X_q^1,\dots\right)=\\
&=\Ad\!\left(\left(\holo{A}{0}{i}\right)^{-1}\!\left(\gamma;t_1,\dots;p,q\right)\right)\\
&\phantom{=\quad}\left(\widetilde{\ev_{i,n}}\circ\pr_2\right)^*\omega_{\left(\gamma;t_1,\dots;p,q\right)}\!\left(X^1,\dots\right),
\end{align*}
where in the first identity I used the horizontality of $\omega$ to get rid of the vertical parts of any argument, and in the second identity I used the $G$-equivariance of $\omega$.
\end{proof}
\begin{Rem}
The formula displayed in Lemma~\ref{lem-explpullhol} simplifies remarkably under the assumption $P$ trivial; it gives in this case the well-known formula used in~\cite{CR} for the construction of the generalized holonomy.
\end{Rem}
On the other hand, as was already proved in Proposition~\ref{prop-simplpartr2}, the parallel transport map $\holo{A}{i}{j}$, for a given connection $A$ on $P$, defines an isomorphism between the bundles $\ev_{i,n}^*P$ and $\ev_{j,n}^*P$; thus, taking a basic differential form $\omega$ on $P$ with values in $\Lg$, representing a form on $M$ with values in the bundle $\ad P$, one can take its pull-backs w.r.t.\ $\widetilde{\ev_{i,n}}$ and $\widetilde{\ev_{j,n}}$ respectively, which is denoted by $\ev_{i,n}^*\omega$ and $\ev_{j,n}^*\omega$.
Taking now $\ev_{j,n}^*\omega$ as a basic form on $\ev_{j,n}^*P$ with values in $\Lg$, one can consider its pull-back w.r.t.\ the isomorphism $\Phi_{A,i,j,n}$, induced by the parallel transport map $\holo{A}{i}{j}$, thus one has a basic form on $\ev_{i,n}^*P$ with values in $\Lg$, which is denoted by
\[
\widehat{\omega}_{A,i,j,n}\colon=\left(\widetilde{\ev_{j,n}}\circ\Phi_{A,i,j,n}\right)^*\omega.
\]
Is there some formula, similar to that displayed in Lemma~\ref{lem-explpullhol}, relating $\widehat{\omega}_{A,i,j,n}$ to $\ev_{j,n}^*\omega$?
This is encoded in the following
\begin{Lem}
For any connection $A$ in $\calA_P$, for any basic form $\omega$ on $P$ with values in $\Lg$ and for any two indices $1\leq i<j\leq n$, the following formula holds:
\[
\pr_1^*\!\left(\widehat{\omega}_{A,i,j,n}\right)=\Ad\!\left((\holo{A}{i}{j})^{-1}\right)\pr_2^*\!\left(\ev_{j,n}^*\omega\right),
\]
where $\pr_k$ denotes the $k$-projection from $\ev_{i,n}^*P\odot\ev_{j,n}^*P$ onto its $k$-factor, for $k=1,2$. \end{Lem}
\begin{proof}
The proof is similar to that of Lemma~\ref{lem-explpullhol}; one only has to change the notations accordingly.
\end{proof}

\subsection{Restrictions to boundary faces of simplices}
I need now some brief comments about the isomorphisms $\Phi_{A,i,n}$ and their restrictions on the faces of the simplices $\triangle_n$; in fact, the following computations will be used later.

Denoting by $\iota_{\alpha,n}$, for $0\leq\alpha\leq n$, the inclusion of $\lo M\times \triangle_{n-1}$ into $\lo M\times\triangle_n$ given by
\[
\iota_{\alpha,n}\left(\gamma;t_1,\dots,t_{n-1}\right)\colon=\begin{cases}
\left(\gamma;0,t_1,\dots,t_{n-1}\right),& \alpha=0\\
\left(\gamma;t_1,\dots,t_\alpha,t_\alpha,\dots,t_{n-1}\right),& 1\leq\alpha\leq n-1\\
\left(\gamma;t_1,\dots,t_{n-1},1\right),& \alpha=n,
\end{cases}
\]
let me also denote by $\widetilde{\iota}_{\alpha,n}$ the natural map from $\ev_{i,n-1}^*P$ into $\ev_{i,n}^*P$.

Now, given two manifolds $M$ and $N$, a principal $G$-bundle $P$ over $N$ and a smooth map $f$ from $M$ to $N$, if one takes the pull-back of any form on $N$ with values in the associated bundle $P\times_G V$ of $P$, for some representation $\left(\rho,V\right)$ of $G$, the resulting form on $M$ takes its values in the pull-back bundle $f^*\left(P\times_G V\right)$, which is canonically isomorphic to $f^*P\times_G V$.
Denoting by $\widetilde{f}$ the natural map from $f^*P$ to $P$, by $\omega$ a general form on $N$ with values in the associated bundle $P\times_G V$, it is well-known that the basic form on $f^*P$ with values in $V$ corresponding to $f^*\omega$ is simply the pull-back w.r.t.\ $\widetilde{f}$ of the basic form on $P$ with values in $V$ corresponding to $\omega$.

Therefore, in order to compute the result of the restriction by $\iota_{\alpha,n}$ of the form $\widehat{\omega}_{A,i,n}$ as a form on $\lo M\times \triangle_n$ with values in $\ev(0)_n^*P\times_G \Lg$, it suffices to compute the result of the following compositions of bundle maps, which is very easy; one can therefore write down only the results:
\begin{equation}\label{eq-boundary}
\begin{aligned}
\widetilde{\ev_{i,n}}\circ\Phi_{A,i,n}\circ\widetilde{\iota}_{0,n}&=\begin{cases}
\widetilde{\ev(0)_{n-1}},& i=1\\
\widetilde{\ev_{i-1,n-1}}\circ\Phi_{A,i-1,n-1},& 2\leq i\leq n
\end{cases},\\
\widetilde{\ev_{i,n}}\circ\Phi_{A,i,n}\circ\widetilde{\iota}_{\alpha,n}&=\begin{cases}
\widetilde{\ev_{i,n-1}}\circ\Phi_{A,i,n-1},& 1\leq i\leq\alpha\\
\widetilde{\ev_{\alpha,n-1}}\circ\Phi_{A,\alpha,n-1},& i=\alpha+1\\
\widetilde{\ev_{i-1,n-1}}\circ\Phi_{A,i-1,n-1},& \alpha+2\leq i\leq n\\
\end{cases},\\
\widetilde{\ev_{i,n}}\circ\Phi_{A,i,n}\circ\widetilde{\iota}_{n,n}&=\begin{cases}
\widetilde{\ev_{i,n-1}}\circ\Phi_{A,i,n-1},& 1<i\leq n-1\\
\widetilde{\ev(0)}\circ \Phi_A\circ\widetilde{\pi}_n,& i=1.
\end{cases}
\end{aligned}
\end{equation}

\subsection{The generalized holonomy and Chen's iterated integrals}\label{ssec-genhol}
Finally, let me define, for a general principal bundle $P$ over $M$, the so-called {\em generalized holonomy} via Chen's iterated integrals.

First of all, consider so-called {\em $BF$-couples}, i.e.\, assuming the manifold $M$ to be of dimension $m$, consider a general couple $(A,B)$, where $i)$ $A$ is a connection on the principal bundle $P$ and $ii)$ $B$ is a basic form on $P$ with values in $\Lg$ of any degree; equivalently, $B$ is a form on $M$ of degree with values in the adjoint bundle $\ad P$.
(For applications in TQFT, especially higher-dimensional $BF$ theories, I will restrict the degree of $B$ to $m-2$, $m$ being the dimension of the manifold $M$.)
Take additionally some representation $(V,\rho)$ of the structure group $G$; the derivative at the identity of the Lie group morphism $\rho$ defines a Lie algebra representation of $\Lg$ in $V$, which is denoted by $\widehat{\rho}$.

For any positive integer $n$ and $1\leq i\leq n$, define the following differential form:  
\[
\widehat{B}_{A,i,n,\rho}:=\widehat{\rho}\!\left(\widehat{B}_{A,i,n}\right),
\]
borrowing the notations from Subsection~\ref{ssec-chenint}.
Hence, for any choice of indices $n$ (and $i$), $\widehat{B}_{A,i,n,\rho}$ defines a basic differential form on the pulled-back bundle $\ev(0)_n^*P$ with values in the Lie algebra of endomorphisms of $V$; equivalently, it may be viewed as a differential form on $\lo M\times \triangle_n$ with values in the associated bundle $\ev(0)_n^*P\times_G \End(V)$, where $G$ acts on $\End(V)$ by a twisted conjugation w.r.t.\ $\rho$.
If one considers then the expression on $\ev(0)_n^*P$ with values in $\End(V)$ given by
\[
\left(\widehat{B}_{A,1,n,\rho}\land\cdots\land \widehat{B}_{A,n,n,\rho}\right)\rho(\holo{A}{0}{1})^{-1},
\]
where the wedge product means not only the usual wedge product of the form-piece but also multiplication in the ring $\End(V)$, it is easy to verify that it defines also a basic form on $\ev(0)_n^*P$ with values in $\End(V)$ of degree $n\deg B$; thus, its trace in the endomorphism ring $\End(V)$ defines a basic differential form on $\ev(0)_n^*P$ of degree $n(m-2)$, hence a differential form on $\lo M\times \triangle_n$.
Integrating this differential form along the fiber w.r.t.\ the projection from $\lo M\times\triangle_n$ onto $\lo M$, forgetting the simplex, one gets a differential form on the space of loops $\lo M$ of degree $n(m-3)$.

\begin{Def}\label{def-genholon}
Define the {\em generalized Wilson loop w.r.t.\ the $BF$-pair $(A,B)$ in the representation $\rho$} as the formal series
\begin{equation}\label{eq-genhol}
\calW_\rho(A,B):=W_\rho(A)+\sum_{n\geq 1} \pi_{n*} \tr_V\left[\left(\widehat{B}_{A,1,n,\rho}\land\cdots\land \widehat{B}_{A,n,n,\rho}\right)\rho(\holo{A}{0}{1})^{-1}\right],
\end{equation}
where $W_\rho(A)$ denotes the Wilson loop function introduced in Subsubsection~\ref{sssec-invfunct} and by $\pi_{n*}$ is denoted the push-forward (or integration along the fiber) w.r.t.\ the forgetful projection from $\lo M\times \triangle_n$ onto $\lo M$.
\end{Def} 

\begin{Rem}
Restricting to $BF$-pairs with $\deg B=m-2$, the generalized Wilson loop $\calW_\rho(A,B)$ is a formal series, whose terms are differential forms of increasing degree, with the exception of the case $m=1,2$, where all push-forwards vanish, because the degree of any integrand is strictly smaller than the dimension of the fiber along which one integrates the terms of the generalized Wilson loop up to order $1$, hence reducing it to the usual Wilson loop.
The dimension $m=3$ presents a different phenomenon: namely, since the degree of the summand indexed by the positive integer $n$ of the generalized Wilson loop is $n(m-3)$, if $m=3$, the degree of the aforementioned summand is $0$, hence for any summand.
Thus, the generalized Wilson loop $\calW_\rho(A,B)$ defines, for every $3$-dimensional manifold $M$, a formal function on the space of loops $\lo M$.
\end{Rem}


Now, recall that the gauge group $\calG_P$ operates from the right on the space of connections $\calA_P$ on $P$ and on basic forms on $P$ with values in $\Lg$ by pull-back; therefore, $\calG_P$ operates also on the space of $BF$-pairs $(A,B)$ from the right.
There is also an operation of $\Omega^{\deg B-1}(M,\ad P)$ on the space of $BF$-pairs, namely
\[
\left((A,B),\tau\right)\mapsto (A,B+\dd_A\tau).
\]
It turns out that both actions, of $\calG_P$ and of $\Omega^{\deg B-1}(M,\ad P)$, combine into a unique right action of the semidirect product group $\calG_P\ltimes \Omega^{\deg B-1}(M,\ad P)$:
\[
\left((A,B),(\sigma,\tau)\right)\mapsto (A^\sigma,B^\sigma+\dd_{A^\sigma}\tau). 
\] 
\begin{Rem}
Before entering into the details (and under the assumption $\deg B=m-2$, let me point out another peculiarity of the case $m=3$: the space of connections $\calA_P$ is an affine space modelled on the space of basic $1$-forms on $P$ with values in $\Lg$.
Hence, a $BF$-pair $(A,B)$ in dimension $3$ may be also viewed as a connection on $P$, decomposed w.r.t.\ to the chosen ``background'' connection $A$.
Notice that the connection $A+B$, corresponding to the $BF$-pair $(A,B)$, is flat, assuming $A$ flat, if and only if $B$ satisfies the $A$-covariant Maurer--Cartan Equation
\[
\dd_A B+\frac{1}2\Lie{B}B=0;
\]
this identity will play a key-r{\^o}le in subsequent computations.
\end{Rem}

The following theorem contains a list of the most important properties of the generalized Wilson loop $\calW_\rho(A,B)$.
\begin{Thm}\label{thm-propgenwils}
The following properties hold for the generalized Wilson loop function $\calW_\rho$:
\begin{itemize}
\item[a)] in any dimension $m$, the map $\calW_\rho$, assigning to any $BF$-pair $(A,B)$ the generalized Wilson loop $\calW_\rho(A,B)$, is invariant w.r.t.\ the action of the gauge group $\calG_P$:
\[
\calW_\rho\left((A,B)^\sigma\right)=\calW_\rho(A,B),\quad \forall \sigma\in \calG_P.
\] 
\item[b)] If one considers the subspace of $BF$-pairs obeying 
\begin{equation}\label{eq-flatBF}
F_A=0,\quad \dd_A B=0,
\end{equation}
then the corresponding generalized Wilson loop $\calW_\rho(A,B)$ is a closed form, assuming $\deg B$ to be odd and such that $2\deg B>m$; assuming $\deg B$ even and such that $2\deg B>m$, then the form made {\em only of the oddly-indexed summands in the generalized Wilson loop} is a closed form.
\item[c)] Assuming $\deg B=1$ and the background connection $A$ to be flat and $B$ obeys the $A$-covariant Maurer--Cartan Equation, then the generalized Wilson loop is a locally constant function on the space of loops in $M$.
\item[d)] If one considers the space of $BF$-pairs $(A,B)$ obeying Equation (\ref{eq-flatBF}), then the infinitesimal action of $\Omega^{\deg B-1}(M,\ad P)$ changes the generalized Wilson loop $\calW_\rho(A,B)$ by an exact form, assuming $2\deg B>m+1$.     
\end{itemize}
\end{Thm}
\begin{proof}
\begin{itemize}
\item[a)] It suffices to show the claim for any summand of the generalized Wilson loop; for this purpose, let me introduce the shorthand notation
\[
\calW_{\rho,n}(A,B)\colon=\pi_{n*} \tr_V\left[\left(\widehat{B}_{A,1,n,\rho}\land\cdots\land \widehat{B}_{A,n,n,\rho}\right)\rho(\holo{A}{0}{1})^{-1}\right].
\]
Now, recall Part $i)$ of Proposition~\ref{prop-covgauge} and that, in the notation used there, 
\[
\calW_{\rho,n}(A,B)=\pi_{n*}\tr_V\left[\Psi_{0,1,n}(A,B)_\rho\land \cdots \land \Psi_{0,n,n}(A,B)_\rho\rho(\holo{A}{0}{1})^{-1}\right],
\]
where the subscript $\rho$ denotes composition with the tangent map at the identity of the homomorphism $\rho$; by abuse of notation, I have written
\[
\Psi_{0,i,n}(A,B)
\]
instead of 
\[
\pr_2\!\left(\Psi_{0,i,n}(A,B)\right),
\]
where $\pr_2$ denotes the projection onto the second factor.
Therefore, it holds
\begin{align*}
\calW_{\rho,n}\!\left((A,B)^{\sigma}\right)&=\pi_{n*}\tr_V\left[\Psi_{0,1,n}\!\left((A,B)^\sigma\right)_\rho\land \cdots \land\rho(\holo{A^\sigma}{0}{1})^{-1}\right]=\\
&=\pi_{n*}\tr_V\left[\Psi_{0,1,n}\!\left((A,B)\right)_\rho^{\sigma_{\ev(0)_n}}\land \cdots\land\rho(\holo{A^\sigma}{0}{1})^{-1}\right]=\\
&=\calW_{\rho,n}(A,B),
\end{align*}
where I made use of the fact that a general gauge transformation operates on a basic form by conjugation, because of the horizontality and of the equivariance, of Proposition~\ref{prop-equivpartr} of Section~\ref{sec-applpullbun} and of the cyclicity of the trace.
\item[b)] First of all, notice that, under the assumption $2\deg B>m$, the form $B_\rho\land B_\rho$, which has obviously degree $2\deg B$ vanishes due to dimensional reason.

Now, I compute exemplarily the exterior derivative of the $n$-th summand $\calW_{\rho,n}(A,B)$ of the generalized Wilson loop; for this purpose, I need the generalized Stokes' Theorem (see~\cite{CR} for more details and to the Appendix for the orientation choices for the $n$-th simplex):
\begin{align*}
\dd \calW_{\rho,n}(A,B)&=(-1)^n \pi_{n*}\dd\tr_V\left[\Psi_{0,1,n}(A,B)_\rho\land\cdots \land \rho(\holo{A}{0}{1})^{-1}\right]-\\
&\phantom{=}-(-1)^n \pi_{\partial_n*}\iota_{\partial_n}^*\tr_V\left[\Psi_{0,1,n}(A,B)_\rho\land\cdots \land \rho(\holo{A}{0}{1})^{-1}\right]=\\
&=(-1)^n \pi_{n*}\tr_V\left\{\dd_{\ev(0)_n^*A_\rho}\left[\Psi_{0,1,n}(A,B)_\rho\land\cdots \land \rho(\holo{A}{0}{1})^{-1}\right]\right\}-\\
&\phantom{=}-(-1)^n \sum_{\alpha=0}^n (-1)^{\alpha+1}\pi_{n-1*}\iota_{\alpha,n}^*\tr_V\left[\Psi_{0,1,n}(A,B)_\rho\land\cdots \land\right.\\
&\phantom{=}\left.\land \rho(\holo{A}{0}{1})^{-1}\right]=\\
&=\sum_{i=1}^n (-1)^{n+(i-1)(m-2)}\pi_{n*}\tr_V\left[\cdots\land\dd_{\ev(0)_n^*A_\rho}\Psi_{0,i,n}(A,B)_\rho\land \cdots\land\right.\\
&\phantom{=}\left. \land \rho(\holo{A}{0}{1})^{-1}\right]+\\
&\phantom{=}+(-1)^{n+n(m-2)}\pi_{n*}\tr_V\left[\Psi_{0,1,n}(A,B)_\rho\land\cdots \land\dd_{\ev(0)_n^*A_\rho}\rho(\holo{A}{0}{1})^{-1}\right]
-\\
&\phantom{=}-\sum_{\alpha=0}^n(-1)^{n+\alpha+1}\pi_{n-1*}\iota_{\alpha,n}^*\tr_V\left[\Psi_{0,1,n}(A,B)_\rho\land\cdots \land \rho(\holo{A}{0}{1})^{-1}\right]=\\
&=\sum_{i=1}^n (-1)^{n+(i-1)(m-2)}\pi_{n*}\tr_V\left[\cdots\land\widehat{\dd_AB}_{A,i,n,\rho}\land \cdots\land\rho(\holo{A}{0}{1})^{-1}\right]-\\
&\phantom{=}-\sum_{\alpha=0}^n(-1)^{n+\alpha+1}\pi_{n-1*}\iota_{\alpha,n}^*\tr_V\left[\Psi_{0,1,n}(A,B)_\rho\land\cdots \land \rho(\holo{A}{0}{1})^{-1}\right]=\\
&=-\sum_{\alpha=0}^n(-1)^{n+\alpha+1}\pi_{n-1*}\iota_{\alpha,n}^*\tr_V\left[\Psi_{0,1,n}(A,B)_\rho\land\cdots \land \rho(\holo{A}{0}{1})^{-1}\right],
\end{align*}
where I made use of Part $ii)$ of Proposition~\ref{prop-equivpartr} of Section~\ref{sec-applpullbun} and of Lemma~\ref{lem-flathoriz} of Section~\ref{sec-flatness}, together with Equation (\ref{eq-flatBF}.
By the notation $\dd_{\ev(0)_n^*A_\rho}$ is meant the covariant derivative on forms on $\ev(0)_n^*P$ with values in the associated bundle $\ev(0)_n^*P\times_G \End(V)$ obtained by twisting the pull-back connection $\ev(0)_n^*A$ by the Lie-algebra representation $\widehat{\rho}$, the tangent map at the identity of $\rho$.

It therefore remains only to compute the boundary contributions; for this purpose, I will make use of the set of equations (\ref{eq-boundary}).
In fact, one gets the following three types of boundary contributions:
\begin{align*}
\iota_{0,n}^*\tr_V\left[\Psi_{0,1,n}(A,B)_\rho\land\cdots \land \rho(\holo{A}{0}{1})^{-1}\right]&=\tr_V\left[\widetilde{\ev(0)_n}^*B_\rho\land\Psi_{0,1,n-1}(A,B)_\rho\land\right.\\
&\phantom{=}\left.\land\cdots\land\Psi_{0,n-1,n-1}(A,B)_\rho\land\rho(\holo{A}{0}{1})^{-1}\right];\\
\iota_{\alpha,n}^*\tr_V\left[\Psi_{0,1,n}(A,B)_\rho\land\cdots \land \rho(\holo{A}{0}{1})^{-1}\right]&=\tr_V\left[\Psi_{0,1,n-1}(A,B)_\rho\land\cdots\land\right.\\
&\phantom{=}\left.\land \widehat{\left(B_\rho\land B_\rho\right)}_{A,\alpha,n-1}\land\cdots \land \rho(\holo{A}{0}{1})^{-1}\right],\\
&\phantom{=} 1\leq \alpha\leq n-1;\\ 
\iota_{n,n}^*\tr_V\left[\Psi_{0,1,n}(A,B)_\rho\land\cdots \land \rho(\holo{A}{0}{1})^{-1}\right]&=\tr_V\left[\Psi_{0,1,n-1}(A,B)_\rho\land\cdots\right.\\
&\phantom{=}\left.\land \Psi_{0,n-1,n-1}(A,B)_\rho\land\right.\\
&\phantom{=}\left.\land\rho(\holo{A}{0}{1})^{-1}\land \widetilde{\ev(0)_n}^*B_\rho\right].
\end{align*}
If $2\deg B>m$, the right-hand side of the second term in the above identities vanishes, for any index $1\leq\alpha\leq n-1$.
Thus, one has only to compare the right-hand sides of the first and the last term; for this purpose, I make use of the cyclicity of the trace, so as to get:
\begin{multline*}
-\sum_{\alpha=0}^n(-1)^{n+\alpha+1}\pi_{n-1*}\iota_{\alpha,n}^*\tr_V\left[\Psi_{0,1,n}(A,B)_\rho\land\cdots \land \rho(\holo{A}{0}{1})^{-1}\right]=\\
=(-1)^n \pi_{n-1*}\tr_V\left[\widetilde{\ev(0)_n}^*B_\rho\land\Psi_{0,1,n-1}(A,B)_\rho\land\right.\\
\phantom{=}\left.\land\cdots\land\Psi_{0,n-1,n-1}(A,B)_\rho\land\rho(\holo{A}{0}{1})^{-1}\right]+\\
+\tr_V\left[\Psi_{0,1,n-1}(A,B)_\rho\land\cdots\land \Psi_{0,n-1,n-1}(A,B)_\rho\land\right.\\
\left.\land \rho(\holo{A}{0}{1})^{-1}\land \widetilde{\ev(0)_n}^*B_\rho\right]=\\
=\left[(-1)^n+(-1)^{(n-1)(m-2)}\right] \pi_{n-1*}\tr_V\left[\widetilde{\ev(0)_n}^*B_\rho\land\Psi_{0,1,n-1}(A,B)_\rho\land\right.\\
\phantom{=}\left.\land\cdots\land\Psi_{0,n-1,n-1}(A,B)_\rho\land\rho(\holo{A}{0}{1})^{-1}\right].
\end{multline*}
Assuming additionally that $\deg B$ is odd, the sum of sign terms before the last expression in the above chain of equalities vanishes, hence yielding the claim for $\deg B$ odd.
On the other hand, if $m$ is even, the aforementioned sum of sign terms equals, for a general $n$,
\[
(-1)^n+1, 
\]
which vanishes if and only if $n$ is odd.
Thus, in even dimensions, the series made of the oddly-indexed summands of the generalized Wilson loop is a closed form, if the $BF$-pair obeys Equation (\ref{eq-flatBF}).
\item[c)] The proof of the claim goes in principle as the proof of Part $b)$, namely I make use of the generalized Stokes' Theorem.
Without repeating all the computations made in Part $b)$ and assuming the background connection $A$ to be flat, since $\deg B=1$ is odd, the remaining boundary terms for the $n$-th summand of the generalized Wilson loop are 
\begin{align*}
&\sum_{i=1}^{n-1}(-1)^{n+i}\pi_{n-1*}\tr_V\left[\Psi_{0,1,n-1}(A,B)_\rho\land\cdots\land \widehat{B_\rho\land B_\rho}_{A,i,n-1}\land\cdots\right.\\
&\left.\land\Psi_{0,n-1,n-1}(A,B)_\rho\land\rho(\holo{A}{0}{1})^{-1}\right]=\\
&=\sum_{i=1}^{n-1}(-1)^{n+i}\pi_{n-1*}\tr_V\left[\Psi_{0,1,n-1}(A,B)_\rho\land\cdots\land\Psi_{0,i,n-1}\!\left(A,\frac{1}2\Lie{B}B\right)_\rho\land\cdots\right.\\
&\left.\land\Psi_{0,n-1,n-1}(A,B)_\rho\land\rho(\holo{A}{0}{1})^{-1}\right],
\end{align*}
where I made use of the fact that $B$ is a $1$-form.

On the other hand, the exterior derivative of the $n-1$-th summand produces a term containing covariant derivatives w.r.t.\ $A$ of $B$, which takes the form:
\begin{align*}
&\sum_{i=1}^{n-1} (-1)^{n+i}\pi_{n-1*}\tr_V\left[\Psi_{0,1,n-1}(A,B)_\rho \land \cdots\land\Psi_{0,i,n-1}(A,\dd_AB)_\rho \land\cdots\land\right.\\
&\left.\land\Psi_{0,n-1,n-1}(A,B)_\rho\land \rho(\holo{A}{0}{1})^{-1}\right].
\end{align*}
It is immediate to see that the boundary term coming from the $n$-summand and the previous term sum up to give
\begin{align*}
&\sum_{i=1}^{n-1} (-1)^{n+i}\pi_{n-1*}\tr_V\left[\Psi_{0,1,n-1}(A,B)_\rho \land \cdots\land\Psi_{0,i,n-1}(A,\dd_AB+\frac{1}2\Lie{B}B)_\rho \land\cdots\land\right.\\
&\left.\land\Psi_{0,n-1,n-1}(A,B)_\rho\land \rho(\holo{A}{0}{1})^{-1}\right],
\end{align*}
which vanishes assuming the connection $A+B$ to be flat, whence the claim follows.

\item[d)] I will show the claim exemplarily for the second summand in the generalized Wilson loop, namely
\[
\calW_{2,\rho}(A,B)=\pi_{1*}\tr_V\left[\widehat{B}_{A,1,1,\rho}\land \rho(\holo{A}{0}1)^{-1}\right].
\] 
Assuming the $BF$-pair to satisfy Equation (\ref{eq-flatBF}), then one has
\[
\calW_{2,\rho}(A,B+\dd_A\tau)=\calW_{2,\rho}(A,B)+\calW_{2,\rho}(A,\dd_A\tau);
\] 
I want to show that the second term on the right-hand side is $dd$-exact.
Namely, again by the generalized Stokes' Theorem, one gets
\begin{align*}
\dd \pi_{1*}\tr_V\left[\Psi_{0,1,1}(A,\tau)\land\rho(\holo{A}{0}1)^{-1}\right]&=-\pi_{1*}\dd\tr_V\left[\Psi_{0,1,1}(A,\tau)\land\rho(\holo{A}{0}1)^{-1}\right]+\\
&\phantom{=}+\iota_{1,1}^*\tr_V\left[\Psi_{0,1,1}(A,\tau)\land\rho(\holo{A}{0}1)^{-1}\right]-\\
&\phantom{=}-\iota_{0,1}^*\tr_V\left[\Psi_{0,1,1}(A,\tau)\land\rho(\holo{A}{0}1)^{-1}\right]=\\
&=-\pi_{1*}\tr_V\left[\Psi_{0,1,1}(A,\dd_A\tau)\land \rho(\holo{A}{0}1)^{-1}\right]+\\
&\phantom{=}+\tr_V\left[\rho(\holo{A}{0}1)^{-1}\land\ev(0)^*\tau\right]-\\
&\phantom{=}-\tr_V\left[\ev(0)^*\tau\land\rho(\holo{A}{0}1)^{-1}\right]=\\
&=\calW_{2,\rho}(A,\dd_A\tau),
\end{align*}
by cyclicity of the trace again.
Thus, if the $BF$-pair obeys Equation (\ref{eq-flatBF}), a variation of the $B$-form by an $A$-covariant exact form produces a change of the second summand of the generalized Wilson loop, which is itself a closed form, {\em independently of the dimension of $m$} (as the previous computation shows immediately by repeating it almost verbatim), by a $\dd$-exact form:
\[
\calW_{2,\rho}(A,B+\dd_A\tau)=\calW_{2,\rho}(A,B)+\dd\widetilde{\calW}_{2,\rho}(A,\tau).
\]
The second term of the generalized Wilson loop, $\calW_{2,\rho}(A,B)$ is, by the above arguments, a ``toy'' example of a closed form on the space of loops in $M$, because there are no restrictions on the dimension of the manifold $M$.

The problems arise when considering the higher summands; the generalized Stokes' Theorem plays again a fundamental r{\^o}le, and one can repeat the arguments used in the proof of Part $b)$ almost verbatim, keeping in mind that, due to dimensional reason, the products
\[
B_\rho\land \tau_\rho\quad\text{or}\quad \tau_\rho\land B_\rho  
\]
vanish, whenever $2\deg B>m+1$.

Notice that the fact that I consider the infinitesimal action of $\Omega^{\deg B-1}(M,\ad P)$ simplifies considerably the computations, because the action on the $n$-th summand decomposes it into a sum of $n$ contributions, each containing at some place between $1$ and $n$ a term $\dd_A\tau$ instead of a form $B$.
Notice also that I have to restrict ourselves to odd dimensions, because of the two boundary contributions coming from the collapse of first coordinate of the simplex to $0$ or of the last one to $1$: namely, it arises a similar phenomenon as for the last part of the proof of Part $b)$.
Similarly, considering also the case of even dimensions, then one has to restrict the generalized Wilson loop only to the summands indexed by odd indices. 
\end{itemize}
\end{proof}

\appendix
\section{The orientation of the $n$th simplex $\triangle_n$}\label{app-simplex}
One of the fundamental ingredients in the construction of the generalized Wilson loop of Definition~\ref{def-genholon} of Section~\ref{sec-chenint} is the $n$-th simplex $\triangle_n$.

In order to compute the push-forward w.r.t.\ the forgetful projection $\pi_n$, I need to specify an orientation of the simplex.
This can be done by choosing an orientation form, which is the restriction to $\triangle_n$ of the standard orientation form on $\bbR^n$, i.e.
\begin{equation}\label{eq-orsimplex}
\dd\!\text{vol}_{\triangle_n}\colon=\dd t_1\land\cdots\land \dd t_n.
\end{equation}
Thus, one can integrate $n$-forms on $\triangle_n$ by the rule: if $\omega=f\dd t_1\land\cdots \land \dd t_n$ is a form of highest degree on $\triangle_n$, where $f$ is a smooth function on $\triangle_n$, its integral is given by
\[
\int_{\triangle_n} \omega\colon=\int_0^1\cdots \int_0^{t_2} f\left(t_1,\dots,t_n\right)\dd t_1\cdots \dd t_n.
\]
Clearly, by Fubini's theorem, one can also write e.g.
\[
\int_{\triangle_n} \omega=\int_0^1 \int_{t_i}^1\cdots\int_{t_i}^{t_{i+2}}\int_0^{t_i}\cdots \int_0^{t_2}f\left(t_1,\dots,t_n\right)\dd t_1\cdots\dd t_{i-1}\dd t_{i+1}\cdots \dd t_n \dd t_i,  
\]
for any $1\leq i\leq n$.

Now that a particular orientation of $\triangle_n$ is specified, one can define the push-forward w.r.t.\ the projection $\pi_n$ from $\lo M\times \triangle_n$ onto $\lo M$ as
\begin{equation}\label{eq-pushsimplex}
\pi_{n*}\left[f\left(\gamma;t_1,\dots\right)\dd\!\text{vol}_{\triangle_n}\land \pi_{n}^*\left(\omega\right)\right]\colon=\left[\int_{\triangle_n}f\left(\gamma;t_1,\dots\right)\dd\!\text{vol}_{\triangle_n}\right]\land\omega,
\end{equation}
otherwise, I set the push-forward to be $0$.
In equation (\ref{eq-pushsimplex}), $f$ is a smooth function on $\lo\times \triangle_n$, while $\omega$ is a smooth form on $\lo M$.
One can verify directly that equation (\ref{eq-pushsimplex}) is a special case of our definition of push-forward.

Now, it is not difficult to compute the orientations of $n+1$ boundary faces of the $n$-simplex $\triangle_n$, which I used explicitly on the proof of Theorem~\ref{thm-propgenwils} of Section~\ref{sec-chenint}.

The boundary $\partial \triangle_n$ of the simplex can be decomposed into a union of $n+1$ different sets as follows:
\[
\partial \triangle_n=\bigcup_{\alpha=0}^n \left(\partial \triangle_n\right)_\alpha,
\]
where 
\begin{align*}
\left(\partial \triangle_n\right)_0&\colon=\left\{\left(t_1,\dots,t_n\right)\in\triangle_n\colon t_1=0\right\};\\
\left(\partial \triangle_n\right)_\alpha&\colon=\left\{\left(t_1,\dots,t_n\right)\in\triangle_n\colon t_\alpha=t_{\alpha+1}\right\},\quad 1\leq \alpha\leq n-1;\\
\left(\partial \triangle_n\right)_n&\colon=\left\{\left(t_1,\dots,t_n\right)\in\triangle_n\colon t_n=1\right\}.
\end{align*}
It is clear that any $\left(\partial \triangle_n\right)_\alpha$ is isomorphic to the $n-1$-simplex $\triangle_{n-1}$; the orientation form on $\triangle_{n-1}$ is $\dd\!\text{vol}_{n-1}$.

Clearly, the $0$-th boundary face $\left(\partial \triangle_n\right)_0$ has $-\frac{\partial}{\partial t_1}$ as outward-directed (normalized) vector field (since $t_i\geq 0$, for all $1\leq i\leq n$).
Since the orientation form of the boundary $\partial M$ of a smooth, oriented manifold $M$ is induced from the orientation form of $M$ by the rule
\[
\dd\!\text{vol}_{\partial M}\colon=\iota_v \dd\!\text{vol}_M,
\]
where $\iota_v$ denotes contraction of the orientation form of $M$ by a given vector-field $v$, which has to be outward-directed from $\partial M$, it follows immediately that the $0$-th boundary face of $\triangle_n$ has orientation form $-\dd\!\text{vol}_{n-1}$, hence one can say that it has orientation $-1$.
In a similar way, by direct computations, it can be proved that any face $\left(\partial \triangle_n\right)_\alpha$ has orientation form $(-1)^{\alpha+1} \dd\!\text{vol}_{n-1}$, hence the orientation of the $n+1$ faces of $\triangle_n$ is given by the rule
\[
\orient\left(\partial \triangle_n\right)_\alpha=(-1)^{\alpha+1}.
\]

\end{document}